\documentclass[12pt]{article}
\usepackage{multirow}
\usepackage{graphicx}
\usepackage{latexsym,amsmath,amsfonts,amscd,amssymb}
\usepackage{epsfig}
\usepackage{epstopdf}
\usepackage{float}
\usepackage{changebar}
\usepackage{indentfirst}
\usepackage{verbatim}
\usepackage{color}
\usepackage{colordvi}

\usepackage{appendix}
\usepackage{cite}

\newcommand{\bx}{\mathbf{x}}
\newtheorem{remark}{Remark}[section]

\begin{document}

\baselineskip=2pc

\begin{center}
{\Large \bf High order finite difference Hermite WENO fast sweeping methods for static Hamilton-Jacobi equations}
\end{center}
\centerline{Yupeng Ren\footnote{School of Mathematical Sciences, Xiamen University,
Xiamen, Fujian 361005, P.R. China. E-mail: ypren@stu.xmu.edu.cn. This work was carried out when Y. Ren was visiting Department of Mathematics, the Ohio State University under the support of the China Scholarship Council (CSC NO. 201906310077).},
Yulong Xing\footnote{Department of Mathematics, The Ohio State University, Columbus, OH 43210, USA.
 E-mail: xing.205@osu.edu. The work of Y. Xing is partially supported by the NSF grant DMS-1753581.} and Jianxian Qiu\footnote{School of Mathematical Sciences and Fujian Provincial
Key Laboratory of Mathematical Modeling and High-Performance
Scientific Computing, Xiamen University,
Xiamen, Fujian 361005, P.R. China. E-mail: jxqiu@xmu.edu.cn.
The work of J. Qiu is partially supported by NSAF grant U1630247.}}
\baselineskip=1.8pc
\vspace{1cm}
\centerline{\bf Abstract}
\bigskip
In this paper, we propose a novel Hermite weighted essentially non-oscillatory (HWENO) fast sweeping method to solve the static Hamilton-Jacobi equations efficiently. During the HWENO reconstruction procedure, the proposed method is built upon a new finite difference fifth order HWENO scheme involving one big stencil and two small stencils. However, one major novelty and difference from the traditional HWENO framework lies in the fact that, we do not need to introduce and solve any additional equations to update the derivatives of the unknown function $\phi$. Instead, we use the current $\phi$ and the old spatial derivative of $\phi$ to update them. The traditional HWENO fast sweeping method is also introduced in this paper for comparison, where additional equations governing the spatial derivatives of $\phi$ are introduced. The novel HWENO fast sweeping methods are shown to yield great savings in both computational time and storage, which improves the computational efficiency of the traditional HWENO scheme. In addition, a hybrid strategy is also introduced to further reduce computational costs. Extensive numerical experiments are provided to validate the accuracy and efficiency of the proposed approaches.


\vfill {\bf Key Words:} finite difference; Hermite methods; weighted essentially non-oscillatory method; fast sweeping method; static Hamilton-Jacobi equations; Eikonal equation


\pagenumbering{arabic}

\newpage

\baselineskip=2pc

\section{Introduction}
\setcounter{equation}{0}
\setcounter{figure}{0}
\setcounter{table}{0}

In this paper, we design and validate high order accurate and efficient Hermite weighted essentially non-oscillatory (HWENO) fast sweeping methods for solving the static Hamilton-Jacobi (HJ) equation
\begin{equation}\label{Y1.0}
\begin{cases}
H(\nabla\phi,\mathbf{x})=0,    & \mathbf{x}\in\Omega \setminus \Gamma,\\
\phi(\mathbf{x})=g(\mathbf{x}),                         & \mathbf{x}\in\Gamma\subset\Omega,
\end{cases}
\end{equation}
where $\Omega$ is the computational domain in $\mathbb{R}^{d}$,
$\phi(\bx)$ is the unknown function in $\Omega$, the Hamiltonian $H$ is a nonlinear Lipschitz continuous function depending on
$\nabla\phi$ and $\bx$, and the boundary condition is given by $\phi(\mathbf{x})=g(\mathbf{x})$ on the subset $\Gamma\subset\Omega$.
One important example to be considered is the Eikonal equation, taking the form of
\begin{equation}\label{eq:eikonal}
\begin{cases}
|\nabla\phi|=f(\mathbf{x}),    & \mathbf{x}\in\Omega \setminus \Gamma,\\
\phi(\mathbf{x})=g(\mathbf{x}),                         & \mathbf{x}\in\Gamma\subset\Omega ,
\end{cases}
\end{equation}
where $f(\mathbf{x})>0$.

The HJ equations have extensive applications in many different fields, for instance in optimal control, computer vision, differential game and geometric optics, image processing and so on \cite{Huang1,Yxia}.
It is well-known that global $C^{1}$ solution does not exist for time-dependent HJ equations in
the generic situation, even if the initial condition is amply smooth. Singularities in the form of discontinuities would appear in the derivatives of the unknown function, hence it is necessary to define a ``weak solution'' for the HJ equations. The viscosity solutions of the HJ equations were first introduced by Crandall and Lions in \cite{viscosity}.


One popular way to numerically solve the static HJ equations is to treat the problem as a stationary boundary value problem, such as the fast marching method (FMM) \cite{Tsitsiklis,Sethian,Helmsen} and the fast sweeping method (FSM) \cite{KaoLF,osher,qianzhangzhao,Tsai,Zhaofsm1,Zhaofsm2} can be applied. Compared with FMM, FSM can be constructed to be high order accurate, and becomes a class of popular and effective methods for solving static HJ equations nowadays.
The FSM was first introduced in \cite{Boue} by Bou\'e and Dupuis, to solve a deterministic control problem with quadratic running cost using Markov chain approximation. Later, Zhao \cite{Zhaofsm1} applied the FSM to solve the Eikonal equations.
 Since then, many high order FSM have been developed to solving static HJ equations. In the framework of finite difference methods, Zhang et. al. \cite{zhang2006} combined the third order finite difference WENO-JP scheme \cite{WENOJP} with FSM, and Xiong et. al. \cite{xiong} studied fifth order WENO-JP FSM scheme. High order accurate boundary treatments (i.e., Richardson extrapolation and Lax-Wendroff type procedure), which are consistent with high order FSM, have been developed for the inflow boundary conditions in \cite{huangnum,xiong}. In \cite{stop}, a competent stopping criteria was recommended for high order FSM. In addition, high order FSM was also investigated in the framework of discontinuous Galerkin (DG) finite element method to solve Eikonal equation, and their numerical performance were shown to be effective and robust \cite{dgli,dgluo,dgwu,dgzhang}.

In additional to finite difference WENO and DG methods, high order HWENO methods \cite{hwenolim1,hwenolim2,qiushuhwenohj2005,lqfdhweno1} have recently gained many attention in solving hyperbolic conservation laws. Both the classical WENO and HWENO methods can achieve the high order accuracy and preserve the essentially non-oscillatory property.
The main difference lies in the fact that the HWENO scheme uses the Hermite reconstruction,  that involves both the unknown variable $\phi$ and its first order spatial derivative or first moment in the reconstruction. As a result, the reconstruction stencil becomes more compact, although more storage and some additional work are needed to evaluate the spatial derivatives.
The HWENO scheme was first proposed during the construction of a suitable limiter for the DG method \cite{hwenolim1,hwenolim2}, since it is more compact than the standard WENO scheme. In \cite{qiushuhwenohj2005}, the HWENO scheme was first used to solve the time-dependent HJ equation, and achieved very good numerical results. Compared with the WENO scheme, its boundary treatment is simpler and the numerical error is smaller with the same mesh, as shown in \cite{qiushuhwenohj2005}. The HWENO scheme was later extended to solve the hyperbolic conservation law in the finite difference framework \cite{lqfdhweno1}, where the same advantages can be observed.
Since then, a series of HWENO schemes \cite{taohwenosta,7thhweno,zhaohyhweno,zhengfdhwenohj} have been investigated to solve hyperbolic conservation laws and time-dependent HJ equations under the framework of finite difference or finite volume methods. Recently, a new HWENO scheme (denoted by HWENO-ZZQ) was developed by Zhu et. al. in \cite{newhwenozq} for time dependent HJ equations. Compared with the existing HWENO methods, the new HWENO-ZZQ reconstruction uses one big stencil and two small stencils. Also, during the reconstruction procedure, one only needs to apply the complicated HWENO reconstruction when updating the function values, and can use the simple high order linear reconstruction when updating the derivatives. As shown in \cite{newhwenozq}, the new HWENO-ZZQ scheme can obtain small errors with the same high order accuracy in the smooth areas, and maintain sharp transitions and non-oscillatory property near discontinuity.

In this work, we propose to combine the HWENO-ZZQ method with the fast sweeping idea, to provide an efficient solver for the static HJ equations. Two approaches to design finite difference HWENO FSM for the static HJ equations will be presented.
The first novel approach is unique to static HJ equation, and addresses the potential concern of HWENO method where additional equations governing the spatial derivatives of $\phi$ are introduced and increase computational costs. Here, we propose to use just one equation to update $\phi$, and then use the current $\phi$ and the old spatial derivatives of $\phi$ to update the spatial derivatives. This is very different from the traditional HWENO framework. By designing the algorithm in this way, there is no need to introduce additional equations for the spatial derivatives, which is commonly employed in the HWENO methods. This will lead to great savings in both computational time and storage, which improves the computational efficiency of the traditional HWENO scheme. For comparison, we also presented the standard HWENO fast sweeping method in this paper, where one equation is computed to update $\phi$, and the other one (or two) equation are computed to update the spatial derivative(s) of $\phi$. This can be viewed as a straightforward extension of HWENO-ZZQ method \cite{newhwenozq} in the framework of FSM.
Finally, we also introduced a hybrid strategy which leads to further saving in the computational resources.

The rest of the paper is organized as follows. In Section \ref{section4}, we present two novel HWENO fast sweeping methods for static HJ equations. The numerical tests are presented to demonstrate the effectiveness and efficiency of our schemes in Section \ref{section5}. A hybrid strategy, together with some numerical results, are presented in Section \ref{section6}. Conclusion remarks are given in Section \ref{section7}.



\section{HWENO FSMs for the static HJ equations}\label{section4}
\setcounter{equation}{0}
\setcounter{figure}{0}
\setcounter{table}{0}
In this section, we present two types of HWENO FSMs to efficiently solve the static HJ equations. The flowchart of these two algorithms and their numerical implementation will be provided. A quick review of finite difference WENO FSM will also be provided.

\subsection{Review of high order WENO FSM}\label{section2}
In this subsection, we briefly review the high order WENO FSM to solve the static HJ equations \cite{xiong,zhang2006}, where WENO reconstruction was used to approximate the first order derivatives appeared in the numerical Hamiltonian. For further details on this subject, we refer to \cite{xiong,zhang2006}.

For ease of presentation, we only consider the following two dimensional static HJ equation
\begin{equation}\label{Y1.1}
\begin{cases}
H(\phi_{x},\phi_{y})=f(x,y),    & (x,y)\in\Omega \setminus \Gamma,\\
\phi(x,y)=g(x,y),                         & (x,y)\in\Gamma\subset\Omega.
\end{cases}
\end{equation}
Suppose the computational domain $\Omega$ is discretized into the rectangular meshes $\Omega_{h}=\{(x_i, y_j), 1\leq i\leq
N_x, 1\leq j\leq N_y\}$, with $(x_i, y_j)$ being a grid point in $\Omega_{h}$.
We denote the numerical solution at the grid point $(x_i, y_j)$ by $\phi_{i,j}$.
$\Delta{x}$ and $\Delta{y}$ stand for the grid sizes in the $x$ and $y$ directions, respectively, and we assume $\Delta x=\Delta y=h$ for simplicity. The numerical approximation of \eqref{Y1.1} is given by
\begin{equation}\label{hamiltonian}
\begin{cases}
\widehat{H}(\phi_{x}^{-},\phi_{x}^{+},\phi_{y}^{-},\phi_{y}^{+})_{ij}=f_{ij},    & (x_i,y_j)\in\Omega_{h} \setminus \Gamma_{h},\\
\phi_{ij}=g_{ij},                         & (x_i,y_j)\in\Gamma_{h}\subset\Omega_{h} ,
\end{cases}
\end{equation}
where $\widehat{H}$ denotes a monotone numerical Hamiltonian which approximates the Hamiltonian $H$.
Such numerical Hamiltonian takes inputs $\phi^{\pm}_x$ and $\phi^{\pm}_y$ at the corresponding
grid point, which needs to be reconstructed from its neighboring point values using the high order WENO procedure.

Two types of numerical Hamiltonian are often considered in the literature.
For general static HJ equation, we adopt the Lax-Friedrichs (LF) numerical Hamiltonian \cite{Oshershu1991}:
  $$ \widehat{H}^{LF}_{i,j}=H\left(\frac{u_{i,j}^{-}+u_{i,j}^{+}}{2},\frac{v_{i,j}^{-}+v_{i,j}^{+}}{2}\right)-
\frac{1}{2}\alpha(u_{i,j}^{+}-u_{i,j}^{-})-
\frac{1}{2}\beta(v_{i,j}^{+}-v_{i,j}^{-}),$$
where
\begin{equation}\label{alphabeta}
\alpha=\max_{u,v}|H_{1}(u,v)|, \qquad \beta=\max_{u,v}|H_{2}(u,v)|.
\end{equation}
Here $H_{\ell}(u, v)\, (\ell=1,2)$ denotes the partial derivative of $H$ with respect to the $\ell$-th argument.
The updating procedure of the LF FSM for static HJ equations can be written as \cite{KaoLF,zhang2006}
\begin{equation}\label{lf}
  \begin{split}
\phi_{i,j}^{new} =&\left(\frac{h}{\alpha+\beta}\right)\bigg[f_{i,j}-H\left(\frac{(\phi_{x})^{+}_{i,j}+(\phi_{x})^{-}_{i,j}}{2},\frac{(\phi_{y})^{+}_{i,j}+(\phi_{y})^{-}_{i,j}}{2}\right) \\
 &\hskip2cm
 +\alpha\frac{(\phi_{x})^{+}_{i,j}-(\phi_{x})^{-}_{i,j}}{2}+\beta\frac{(\phi_{y})^{+}_{i,j}-(\phi_{y})^{-}_{i,j}}{2}\bigg]+\phi_{i,j}^{old}.
  \end{split}
\end{equation}
Here $\phi_{i,j}^{new}$ denotes the updated numerical approximations of $\phi$ at the grid point
$(x_i, y_j )$ and $\phi_{i,j}^{old}$ denotes the previous value of $\phi$ at the same grid point.

The other commonly used numerical Hamiltonian is the Godunov numerical Hamiltonian, often employed in the approximation of the
Eikonal equation \eqref{eq:eikonal}. Again, we consider the two dimensional version
\begin{equation}\label{41}
  \begin{cases}
 \displaystyle \sqrt{\phi_{x}^2+\phi_{y}^{2}}=f(x,y),&(x,y)\in\Omega,\\
 \displaystyle \phi(x,y)=g(x,y),&(x,y)\in\Gamma\subset\Omega,
  \end{cases}
\end{equation}
and utilize the following Godunov numerical Hamiltonian to approximate it on uniform meshes \cite{Zhaofsm1,zhang2006}
\begin{equation}\label{g1}
  \left[\left(\frac{{\phi}_{i,j}^{new}-\phi_{i,j}^{x min}}{h}\right)^{+}\right]^{2}+\left[\left(\frac{{\phi}_{i,j}^{new}-\phi_{i,j}^{y min}}{h}\right)^{+}\right]^{2}=f_{i,j}^{2},\quad x^{+}=\begin{cases}x,\quad x>0,\\
  0,\quad x<0,\end{cases}
\end{equation}
where
\begin{subequations}\label{phixymin}
\begin{align}
  &\phi_{i,j}^{xmin}=\min(\phi_{i,j}^{old}-h(\phi_{x})_{i,j}^{-},\phi_{i,j}^{old}+h(\phi_{x})_{i,j}^{+}),\\
  &\phi_{i,j}^{ymin}=\min(\phi_{i,j}^{old}-h(\phi_{y})_{i,j}^{-},\phi_{i,j}^{old}+h(\phi_{y})_{i,j}^{+}).
\end{align}
\end{subequations}
After obtaining $\phi^{xmin}_{i,j}$ and $\phi^{ymin}_{i,j}$ using the above formulas, ${\phi}_{i,j}^{new}$ can be computed as
\begin{equation}
\label{phinew}
{\phi}_{i,j}^{new} =
\begin{cases}
&\displaystyle \min (\phi^{xmin}_{i,j},\phi^{ymin}_{i,j}) + f_{i,j} h, \, \hskip25mm \text{ if } |\phi^{xmin}_{i,j}-\phi^{ymin}_{i,j}|\ge f_{i,j} h, \\
&\displaystyle  \frac12\left(\phi^{xmin}_{i,j}+\phi^{ymin}_{i,j}+(2f^2_{i,j}h^2-(\phi^{xmin}_{i,j}-\phi^{ymin}_{i,j})^2)^{1/2}\right), \,\hskip5mm \text{ otherwise}.
\end{cases}
\end{equation}
A systematic method for solving the Eikonal equations by first order FSM is developed in \cite{Zhaofsm1}. The essential idea of the FSM is to adopt nonlinear upwind difference and Gauss-Seidel (GS) iterations with alternating sweeping ordering. The FSM follows the causality along characteristics, namely, all characteristic lines are classified as finite groups according to their directions, and each GS iteration with a specific sweep order covers a set of characteristics lines. We can refer to \cite{Zhaofsm1,zhang2006} for more details including the flowchart of FSM.
\subsection{A novel HWENO FSM}\label{sec:approach1}

\indent Before introducing the new method, we start by briefly reviewing the traditional HWENO framework. We denote $u=\phi_{x}(x, y)$ and $v=\phi_{y}(x, y)$ as the first order partial derivatives of $\phi$ with respect to the variables $x$ and $y$, respectively. By taking spatial derivatives on both sides of \eqref{Y1.1}, we obtain the following system of equations:
\begin{equation}\label{Y1.2}
\begin{cases}
H(\phi_{x},\phi_{y})=f(x,y),    \\
H_{1}(u,v)u_{x}+H_{2}(u,v)u_{y}=f_{x},\\
H_{1}(u,v)v_{x}+H_{2}(u,v)v_{y}=f_{y},
\end{cases}
\end{equation}
where $H_{1}(u,v)=\frac{\partial H}{\partial u}$, $H_{2}(u,v)=\frac{\partial H}{\partial v}$, and $v_{x}=u_{y}$ is used in the derivation.

The first equation can be solved by the FSM \eqref{lf} or \eqref{phinew}, combined with the HWENO-ZZQ reconstruction to be discussed in Section \ref{section3}, which involves $\phi$, $u$ and $v$ simultaneously.
The auxiliary variables $u$ and $v$ are usually updated by solving these two equations arisen from the derivative of the HJ equation, which is a common approach in the traditional HWENO method for time-dependent problem \cite{qiushuhwenohj2005,shutimehj,zhengfdhwenohj}.

The targeting HJ equation in this paper is a steady state problem, and iterative method is used to update our approximation of $\phi$. During each iteration, we have already applied the HWENO-ZZQ reconstruction procedure to evaluate $\phi_{x}^{\pm}$ and $\phi_{y}^{\pm}$ when solving the first equation of \eqref{Y1.2}. Such information could be reused to generate our updated numerical approximation of $u$ and $v$, and there is no need to re-evaluate them from solving these two additional equations. In other words, we use $\phi^{new}$ (computed by \eqref{lf} or \eqref{phinew}) and $u^{old}$, $v^{old}$ to reconstruct $\phi_{x}^{\pm}$ and $\phi_{y}^{\pm}$ by HWENO-ZZQ reconstruction, and then define $u^{new}$ and $v^{new}$ from the following formulas:
\begin{equation}\label{app2}
  u_{i,j}^{new}=\begin{cases}
  (\phi_{x})_{i,j}^{-}, &if~ (\phi_{x})_{i,j}^{\pm}>0,\\
  (\phi_{x})_{i,j}^{+}, &if~ (\phi_{x})_{i,j}^{\pm}<0,\\
  u_{i,j}^{old},& otherwise,
  \end{cases}
  \qquad
       v_{i,j}^{new}=\begin{cases}
  (\phi_{y})_{i,j}^{-}, &if~ (\phi_{y})_{i,j}^{\pm}>0,\\
  (\phi_{y})_{i,j}^{+}, &if~ (\phi_{y})_{i,j}^{\pm}<0,\\
  v_{i,j}^{old},& otherwise.
  \end{cases},
\end{equation}
Here $u_{i,j}^{new}$ and $v_{i,j}^{new}$ denote the updated numerical approximations of $u$ and $v$ at the grid point
$(x_i, y_j )$, respectively. We use $u_{i,j}^{old}$
 and $v_{i,j}^{old}$ to denote the previous value of $u$ and $v$ at the same grid point. The HWENO-ZZQ reconstruction, to approximate the derivatives $\phi_{x}$, and $\phi_{y}$ at the grid point $(x_i, y_j )$ with high order accuracy, will be discuss in Section \ref{section3}.

The definition of $u^{new}$ and $v^{new}$ in \eqref{app2} comes from the fast sweeping method, that is, the solution of equation \eqref{Y1.1} is increasing along the characteristic lines \cite{qianzhangzhao,zhang2006}, and the information always comes from the upwind direction. Therefore we use this formulation to directly define the derivative values. Numerical tests in Section \ref{section5} also confirm that this simplified method is both robust and effective.

We would like to comment that this approach is different from the standard WENO method, and cannot be used in time-dependent problems. In the HWENO procedure to reconstruct the derivatives $\phi_{x}^{\pm}$ and $\phi_{y}^{\pm}$ (hence, $u^{new}$, $v^{new}$), the information of $\phi^{new}$, $u^{old}$ and $v^{old}$ are all used. Since this is an iterative method, $u^{old}$ and $v^{old}$ would also be ``good'' approximations of the exact derivatives. For time-dependent problems, $\phi^{new}$ would be the approximation at the next time step $t^n+\Delta t$, while $u^{old}$ and $v^{old}$ approximates the derivatives at the current time step $t^n$, hence such reconstruction cannot be applied.

 \subsection{The traditional HWENO FSM}\label{sec:approach2}
In this subsection, we present how to solve the static HJ equation when the traditional HWENO framework is used.
The first equation in \eqref{Y1.2} can be solved by the FSM \eqref{lf} or \eqref{phinew}, combined with the HWENO-ZZQ reconstruction. Here we describe how to solve the last two auxiliary equations to obtained the updated values of $u_{i,j}$ and $v_{i,j}$. These equations are approximated by the following scheme
\begin{equation}\label{fuzhu}
\begin{cases}
\widetilde{H}_{i,j}=(f_x)_{i,j},\\
\widetilde{\widetilde{H}}_{i,j}=(f_y)_{i,j},
\end{cases}
\end{equation}
where the $\widetilde{H}_{i,j}$ and $\widetilde{\widetilde{H}}_{i,j}$ are the numerical flux defined as
\begin{align}
      \widetilde{H}_{i,j}=&H_{1}\left(\frac{u_{i,j}^{-}+u_{i,j}^{+}}{2},\frac{v_{i,j}^{-}+v_{i,j}^{+}}{2}\right)\frac{u_{x~ij}^{+}+u_{x~ij}^{-}}{2}+
H_{2}\left(\frac{u_{i,j}^{-}+u_{i,j}^{+}}{2},\frac{v_{i,j}^{-}+v_{i,j}^{+}}{2}\right)\frac{u_{y~ij}^{+}+u_{y~ij}^{-}}{2} \notag \\
       &-
\frac{1}{2}\alpha(u_{x~ij}^{+}-u_{x~ij}^{-})- \frac{1}{2}\beta(u_{y~ij}^{+}-u_{y~ij}^{-}),	\label{huw} \\
      \widetilde{\widetilde{H}}_{i,j}=&H_{1}\left(\frac{u_{i,j}^{-}+u_{i,j}^{+}}{2},\frac{v_{i,j}^{-}+v_{i,j}^{+}}{2}\right)\frac{v_{x~ij}^{+}+v_{x~ij}^{-}}{2}+
H_{2}\left(\frac{u_{i,j}^{-}+u_{i,j}^{+}}{2},\frac{v_{i,j}^{-}+v_{i,j}^{+}}{2}\right)\frac{v_{y~ ij}^{+}+v_{y~ij}^{-}}{2}		\notag\\
       &-
\frac{1}{2}\alpha(v_{x~ij}^{+}-v_{x~ij}^{-})-\frac{1}{2}\beta(v_{y~ij}^{+}-v_{y~ij}^{-}),  \notag
\end{align}
with $\alpha$ and $\beta$ given in \eqref{alphabeta}.

The iteration scheme for updating $u$ will be discussed below and the procedure for $v$ is exactly the same. We start by discussing the simpler first order case, where $u^{\pm}$ are simply the backward and forward difference approximations. Hence, the first equation of \eqref{fuzhu} can be rewritten as
\begin{equation}\label{huw1}
  \begin{split}
      \widetilde{H}_{i,j}=&H_{1}\left(\frac{\phi_{i+1,j}-\phi_{i-1,j}}{2h},\frac{\phi_{i,j+1}-\phi_{i,j-1}}{2h}\right)\frac{u_{i+1,j}-u_{i-1,j}}{2h}\\
      &+H_{2}\left(\frac{\phi_{i+1,j}-\phi_{i-1,j}}{2h},\frac{\phi_{i,j+1}-\phi_{i,j-1}}{2h}\right)\frac{u_{i,j+1}-u_{i,j-1}}{2 h} \\
       &-
\frac{1}{2h}\alpha(u_{i+1,j}-2u_{i,j}+u_{i-1,j})-
\frac{1}{2h}\beta(u_{i,j+1}-2u_{i,j}+u_{i,j-1})=(f_{x})_{i,j}.
   \end{split}
\end{equation}
Therefore, one can solve for the first order approximation of $u_{i,j}$ with the following expression
\begin{equation}\label{huw2}
  \begin{split}
  u_{i,j}=&\frac{h}{\alpha+\beta}[(f_{x})_{i,j}-H_{1}\left(\frac{\phi_{i+1,j}-\phi_{i-1,j}}{2h},\frac{\phi_{i,j+1}-\phi_{i,j-1}}{2h}\right)\frac{u_{i+1,j}-u_{i-1,j}}{2h}\\
  & \qquad\quad
  -H_{2}\left(\frac{\phi_{i+1,j}-\phi_{i-1,j}}{2h},\frac{\phi_{i,j+1}-\phi_{i,j-1}}{2h}\right)\frac{u_{i,j+1}-u_{i,j-1}}{2 h} \\
  & \qquad\quad
+\frac{1}{2h}\alpha(u_{i+1,j}+u_{i-1,j})+
\frac{1}{2h}\beta(u_{i,j+1}+u_{i,j-1})].
   \end{split}
\end{equation}
To obtain the high order iterative scheme, we replace $u_{i+1,j}, u_{i-1,j},u_{i,j+1}$ and $u_{i,j-1}$ with $u_{i,j}+h(u_{x})_{i,j}^{+}, u_{i,j}-h(u_{x})_{i,j}^{-},u_{i,j}+h(u_{y})_{i,j}^{+}$ and $u_{i,j}-h(u_{y})_{i,j}^{-}$, respectively (see \cite{zhang2006}), where $(u_{x})_{i,j}^{\pm}$ and $(u_{y})_{i,j}^{\pm}$ are high order approximation of the partial derivatives of $u$. Similarly, one can apply this idea to $\phi_{i+1,j}, \phi_{i-1,j},\phi_{i,j+1}$ and $\phi_{i,j-1}$ as well, and the resulting high order schemes can be rewritten as
\begin{equation}\label{updateu}
  \begin{split}
  u^{new}_{i,j}=&\frac{1}{\frac{\alpha}{h}
  +\frac{\beta}{h}}\Big[(f_{x})_{i,j}-H_{1}\left(\frac{(\phi_{x}^{+})_{i,j}+(\phi_{x}^{-})_{i,j}}{2},\frac{(\phi_{y}^{+})_{i,j}+(\phi_{y}^{-})_{i,j}}{2}\right)
  \frac{(u_{x}^{+})_{i,j}+(u_{x}^{-})_{i,j}}{2}\\
  & \qquad\quad
  -H_{2}\left(\frac{(\phi_{x}^{+})_{i,j}+(\phi_{x}^{-})_{i,j}}{2},\frac{(\phi_{y}^{+})_{i,j}+(\phi_{y}^{-})_{i,j}}{2}\right)
  \frac{(u_{y}^{+})_{i,j}+(u_{y}^{-})_{i,j}}{2} \\
    & \qquad\quad
    +\frac{1}{2}\alpha((u_{x}^{+})_{i,j}-(u_{x}^{-})_{i,j})+
\frac{1}{2}\beta((u_{y}^{+})_{i,j}-(u_{y}^{-})_{i,j})\Big]+u_{i,j}^{old}.
   \end{split}
\end{equation}
The high order scheme to solve $v$ can be obtained in the similar way, and takes the form of
\begin{equation}\label{updatev}
  \begin{split}
  v^{new}_{i,j}=&\frac{1}{\frac{\alpha}{h}
  +\frac{\beta}{h}}\Big[(f_{y})_{i,j}-H_{1}\left(\frac{(\phi_{x}^{+})_{i,j}+(\phi_{x}^{-})_{i,j}}{2},\frac{(\phi_{y}^{+})_{i,j}+(\phi_{y}^{-})_{i,j}}{2}\right)
  \frac{(v_{x}^{+})_{i,j}+(v_{x}^{-})_{i,j}}{2}\\
  & \qquad\quad
  -H_{2}\left(\frac{(\phi_{x}^{+})_{i,j}+(\phi_{x}^{-})_{i,j}}{2},\frac{(\phi_{y}^{+})_{i,j}+(\phi_{y}^{-})_{i,j}}{2}\right)
  \frac{(v_{y}^{+})_{i,j}+(v_{y}^{-})_{i,j}}{2} \\
   & \qquad\quad
   +\frac{1}{2}\alpha((v_{x}^{+})_{i,j}-(v_{x}^{-})_{i,j})+
\frac{1}{2}\beta((v_{y}^{+})_{i,j}-(v_{y}^{-})_{i,j})\Big]+v_{i,j}^{old}.
   \end{split}
\end{equation}
The HWENO-ZZQ reconstruction will be used to approximate the derivatives $\phi_{x},\phi_{y},u_{x}, u_{y},v_{x}$ and $v_{y}$ at the grid point $(x_i, y_j )$ with high order accuracy.

\subsection{HWENO-ZZQ reconstruction}\label{section3}
In this subsection, the new finite difference HWENO-ZZQ reconstruction recently proposed in \cite{newhwenozq} will be briefly reviewed.
To save space, we only illustrate the reconstruction of $(\phi_{x})^{\pm}_{i,j}$ along $x$-direction here. The approximation of $(\phi_{y})^{\pm}_{i,j}$ along $y$-direction can be obtained similarly, and we refer to \cite{newhwenozq} for more details.
\begin{itemize}
  \item Reconstruction of $(\phi_{x})_{i,j}^{-}$ from upwind information: \\
Take a big stencil $S_{0}=\{x_{i-2}, x_{i-1}, x_{i}, x_{i+1}\}$ and two small stencils $S_{1}=\{x_{i-2},  x_{i-1}, x_{i}\}$, $S_{2}=\{x_{i-1}, x_{i}, ,x_{i+1}\}$, we compose a Hermite quintic polynomial $p^{-}_{1}(x)$, and two quadratic polynomials $p^{-}_{2}(x)$, $p^{-}_{3}(x)$ satisfying
\begin{equation*}
\begin{split}
 p^{-}_{1}(x_{k})&=\phi_{k,j},\quad k=i-2,\cdots,i+1,\,\,\mathrm{ and }\,\,~(p^{-}_{1})'|_{x_{k}}=u_{k,j},\quad k=i-1,i+1;\\
p^{-}_{2}(x_{k})&=\phi_{k,j},\quad k=i-2,i-1,i;\,\,\quad \mathrm{ and }\,\,~p^{-}_{3}(x_{k})=\phi_{k,j},\quad k=i-1,i,i+1;
\end{split}
\end{equation*}
The values of their first-order derivative at $x=x_i$ can be evaluated as
\begin{subequations}
\begin{align}
&(\phi_{x})_{i,j}^{-,1}=(p^{-}_{1})'|_{x_{i}}=\frac{\phi_{i-2,j}+18\phi_{i-1,j}-9\phi_{i,j}-10\phi_{i+1,j}+9hu_{i-1,j}+3hu_{i+1,j}}{-18h};	\label{hybridlinearf}\\
&(\phi_{x})_{i,j}^{-,2}=(p^{-}_{2})'|_{x_{i}}=\frac{\phi_{i-2,j}-4\phi_{i-1,j}+3\phi_{i,j}}{2h};\\	
 &(\phi_{x})_{i,j}^{-,3}=(p^{-}_{3})'|_{x_{i}}=\frac{-\phi_{i-1,j}+\phi_{i+1,j}}{2h}.
\end{align}
\end{subequations}
  \item Reconstruction of $(\phi_{x})_{i,j}^{+}$ from downwind information: \\
Take a big stencil $\widetilde{S}_{0}=\{x_{i-1}, x_{i}, x_{i+1},x_{i+2}\}$ and two small stencils $\widetilde{S}_{1}=\{x_{i-1}, x_{i},x_{i+1}\}$, $\widetilde{S}_{2}=\{x_{i}, x_{i+1}, ,x_{i+2}\}$, we compose a Hermite quintic polynomial $p^{+}_{1}(x)$, and two quadratic polynomials $p^{+}_{2}(x)$, $p^{+}_{3}(x)$ such that
\begin{equation*}
\begin{split}
 p^{+}_{1}(x_{k})&=\phi_{k,j},\quad k=i-1,\cdots,i+2, \,\,\mathrm{and}\,\,~(p^{+}_{1})'|_{x_{k}}=u_{k,j},\quad k=i-1,i+1;\\
p^{+}_{2}(x_{k})&=\phi_{k,j},\quad k=i-1,i,i+1;\,\,\quad \mathrm{ and }\,\,~p^{+}_{3}(x_{k})=\phi_{k,j},\quad k=i,i+1,i+2;
\end{split}
\end{equation*}
The values of their first-order derivative at $x=x_i$ can be evaluated as
\begin{subequations}
\begin{align}
&(\phi_{x})_{i,j}^{+,1}=(p^{+}_{1})'|_{x_{i}}=\frac{10\phi_{i-1,j}+9\phi_{i,j}-18\phi_{i+1,j}-\phi_{i+2,j}+3hu_{i-1,j}+9hu_{i+1,j}}{-18h};	\label{hybridlinearz}\\
&(\phi_{x})_{i,j}^{+,2}=(p^{+}_{2})'|_{x_{i}}=\frac{\phi_{i+1,j}-\phi_{i-1,j}}{2h};	\\
 &(\phi_{x})_{i,j}^{+,3}=(p^{+}_{3})'|_{x_{i}}=\frac{-3\phi_{i,j}+4\phi_{i+1,j}-\phi_{i+2,j}}{2h}.
\end{align}
\end{subequations}
\end{itemize}

In the nonlinear HWENO reconstructions, $(\phi_{x})_{i,j}^{\pm}$ are computed as a convex combination of these three corresponding values \cite{Levy1,wenozqhj,wenozqhy}
\begin{equation}\label{weno}
(\phi_{x})_{i,j}^{\pm}=\omega_{1}^{\pm}\left(\frac{1}{\gamma_{1}}(\phi_{x})_{i,j}^{\pm,1}-
\frac{\gamma_{2}}{\gamma_{1}}(\phi_{x})_{i,j}^{\pm,2}-\frac{\gamma_{3}}{\gamma_{1}}(\phi_{x})_{i,j}^{\pm,3}\right)
+\omega_{2}^{\pm}(\phi_{x})_{i,j}^{\pm,2}+\omega_{3}^{\pm}(\phi_{x})_{i,j}^{\pm,3},
\end{equation}
where the parameters $\omega_n~(n=1,2,3)$ and $\gamma_n~(n=1,2,3)$ are called the nonlinear weights and linear weights, respectively. The parameters $\gamma_n $ can be any positive constants that satisfy $\gamma_1+\gamma_2+\gamma_3=1$, and $\omega_{n}$ can be computed from
\begin{equation}\label{epsilon}
\omega_{n}^{\pm}=\frac{\overline{\omega}_{n}^{\pm}}{\sum_{l=1}^{3}\overline{\omega}_{l}^{\pm}},\qquad
\overline{\omega}_{n}=\gamma_{n}\left(1+\frac{\tau^{\pm}}{\epsilon+\beta_{n}^{\pm}}\right),\qquad n=1,2,3,
\end{equation}
where $\epsilon$ is a small positive number to avoid the denominator becoming $0$, and
$$\tau^{\pm}=\left(\frac{|\beta_{1}^{\pm}-\beta_{2}^{\pm}|+|\beta_{1}^{\pm}-\beta_{3}^{\pm}|}{2}\right)^{2},\mathrm{and}~\beta^{\pm}_{n}=\sum_{\alpha=2}^{r}\int^{x_{i+\frac{1}{2}}}_{x_{i-\frac{1}{2}}}h^{2\alpha-3}\left(\frac{d^{\alpha}p^{\pm}_{n}(x)}{dx^{\alpha}}\right)^2dx,~n=1,2,3,
$$
 where $\beta_{n}^{\pm}$ are the so-called smoothness indicators, which measure how smooth the first-order derivative functions of $p^{\pm}_{n}(x)$ are near the target point $x_{i}$, and $r=5$ for $n=1$, and $r=2$ for $n=2,3$, respectively.


The approximation of $(u^{\pm}_{x})_{i,j}$ is based on the high order linear reconstructions in x-direction, instead of the HWENO reconstructions. Given two big spatial stencils $Q=\{x_{i-2},x_{i-1},x_{i},x_{i+1}\}$ and $\widetilde{Q}=\{x_{i-1},x_{i},x_{i+1},x_{i+2}\}$, then we can construct two Hermite seventh order polynomials $q^{\pm}(x)$ such that
\begin{equation*}
\begin{split}
 q^{-}(x_{k})&=\phi_{k,j},\quad k=i-2,\cdots,i+1,\,\,\mathrm{and}\,\,~(q^{-})'|_{x_{k}}=u_{k,j},\quad k=i-1,i,i+1;\\
q^{+}(x_{k})&=\phi_{k,j},\quad k=i-1,\cdots,i+2, \,\,\mathrm{and}\,\,~(q^{+})'|_{x_{k}}=u_{k,j},\quad k=i-1,i,i+1.\\
\end{split}
\end{equation*}
One can evaluate their second-order derivatives at $x=x_i$ as:
\begin{equation*}
\begin{split}
(u^{-}_{x})_{i,j}=(q^{-})''|_{x_{i}}&=\frac{\phi_{i-2,j}+54\phi_{i-1,j}-81\phi_{i,j}+26\phi_{i+1,j}+18hu_{i-1,j}+18hu_{i,j}-6hu_{i+1,j}}{18h^{2}};\\
(u^{+}_{x})_{i,j}=(q^{+})''|_{x_{i}}&=\frac{26\phi_{i-1,j}-81\phi_{i,j}+54\phi_{i+1,j}+\phi_{i+2,j}+6hu_{i-1,j}-18hu_{i,j}-18hu_{i+1,j}}{18h^{2}}.
\end{split}
\end{equation*}
 The approximation of $(v_{y})^{\pm}_{i,j}$ along y-direction can be obtained in a similar way and is skipped here. The mixed derivative $u_{y}$ and $v_{x}$ can be evaluated using central difference easily, since they play smaller role on the spurious oscillations according to \cite{qiushuhwenohj2005}. Therefore, one can use the fourth order central approximations in the $x$ and
$y$ directions, and obtain
\begin{equation*}
  \begin{split}
     (u_{y})_{i,j}\approx\frac{-u_{i,j+2}+8u_{i,j+1}-8u_{i,j-1}+u_{i,j-2}}{12h},(v_{x})_{i,j}\approx\frac{-v_{i+2,j}+8v_{i+1,j}-8v_{i-1,j}+v_{i-2,j}}{12h}.
   \end{split}
\end{equation*}
We then set $u_{y}^{+}=u_{y}^{-}=u_{y}$, and $v_{x}^{+}=v_{x}^{-}=v_{x}$ at the point $(x_i,y_j)$. This finishes the description of the HWENO-ZZQ reconstruction.
\subsection{The flowchart of both approaches} \label{sec4.3}

 We have discussed two approaches to solve the static HJ equations. For simplicity, let us denote the novel HWENO FSM in Section \ref{sec:approach1}
 by Approach 1, and the traditional HWENO FSM in Section \ref{sec:approach2} by Approach 2.
Next we will summarize the detailed procedure of these two approaches, and provide a flowchart for them. We start by labelling the points $\{(x_i, y_j)\}$ into several categories as in \cite{me}:

\noindent\emph{Category I}: For points on the boundary $\Gamma$, the exact values are assigned for these points.

\noindent\emph{Category II}: For ghost points (exterior of the boundary), we use the high order extrapolation to compute their numerical
solution $\phi_{i,j}$.

\noindent\emph{Category III}: For points near the $\Gamma$ (whose distances to $\Gamma$ are less than or equal to $2h$). The numerical boundary treatment from \cite{xiong,tanshuinverse} is used (i.e., Richardson extrapolation for a single point or a set of isolated points, while Lax-Wendroff type procedure for continuous $\Gamma$).


\noindent\emph{Category IV}: All remaining points. which will be updated by FSM.

Note that only \emph{Category IV} points need to be updated by following sweepings. We now summarize our flowchart for two approaches as follows:

\noindent \textbf{Step 1}. \emph{Initialization}:
The numerical solution from the first order fast sweeping method \cite{Zhaofsm1} is taken as the initial guess of $\phi$. The forward or backward difference of this $\phi$ is used as the initial guess of $u$ and $v$.

\noindent\textbf{Step 2}. \emph{Gauss-Seidel iteration}. We solve the discretized nonlinear system \eqref{Y1.2} by GS iterations with four alternating direction sweepings:
$$(1)\quad i=1:N_x, ~j=1:N_y; \qquad (2)\quad i=N_x:1,~ j=1:N_y;$$
$$(3)\quad i=N_x:1,~ j=N_y:1; \qquad (4)\quad i=1:N_x,~ j=N_y:1.$$
In this step, we first compute the $\phi_{i,j}^{new}$ by \eqref{lf} or \eqref{phinew}, then
\begin{itemize}
  \item for Approach 1, $\phi_{i,j}^{new}$ $u_{i,j}^{old}$ and $v_{i,j}^{old}$ are used to reconstruct $u^{new}_{i,j}$ and $v^{new}_{i,j}$ following \eqref{app2}.
  \item for Approach 2, $u_{i,j}^{new}$ and $v_{i,j}^{new}$ are evaluated by \eqref{updateu}-\eqref{updatev} in each sweeping direction.
  Note that we need to apply linear reconstruction to obtain $u_{x}, u_{y},v_{x}$ and $v_{y}$ first, as explained in Section \ref{section3}.
\end{itemize}
The values at ghost points will be updated by high order extrapolations in both approaches.

\noindent\textbf{Step 3}. \emph{Convergence}: In general, the iteration will stop if, for two consecutive iteration steps,
$$\delta=||\phi^{new}-\phi^{old}||_{L_{1}}<10^{-14}.$$

\subsection{Comments and remarks}
At the end of this section, we would like to present some comments and remarks about the proposed algorithms.
Approach 1 needs only one equation to update $\phi$, while its derivative approximations are obtained by applying the HWENO-ZZQ reconstruction on the updated $\phi$ and old $u$, $v$. As a comparison, Approach 2 adopted the traditional HWENO idea, namely, one equation to update original variable $\phi$, and two auxiliary equations to update the derivative values. The major difference between Approach 1 and Approach 2 is that the latter one need extra
work and storage for the auxiliary variables. As a result, Approach 1 could greatly save computational cost.
\begin{remark} For the Eikonal equation, the numerical tests \cite{xiong,me} show that if the Godunov flux \eqref{g1} is used to solve the first equation, both approaches will not converge to machine epsilon, i.e. $\delta$ will not decrease to $10^{-14}$, especially for the examples with singularities. Therefore, following the idea in \cite{relaxlf}, we propose to update the solution by
$$\phi^{new}=\omega\phi^{new}+(1-\omega)\phi^{old}, ~0<\omega<1.$$
This fix is shown to yield good convergence, although it may slightly increase the number of iterations. Numerically, one observes that $\omega=0.7$ or $0.8$ is the optimal choice.
When the LF numerical Hamiltonian is considered, it is not necessary to take $\omega<1$.
\end{remark}

\begin{remark}
In Approach 2, we use \eqref{updateu} and \eqref{updatev} to update the auxiliary variables $u$ and $v$, respectively. We could increase the values of $\alpha$ and $\beta$ appropriately to improve the convergence speed. Numerically, we observe that doubling their values can reduce the iteration numbers.
\end{remark}

\section{Numerical examples}\label{section5}
\setcounter{equation}{0}
\setcounter{figure}{0}
\setcounter{table}{0}

In this section, we will present extensive numerical examples by testing the proposed fifth order finite difference HWENO FSM for the Eikonal equations and general static HJ equations in two dimensions. We will compare the numerical results of these two approaches with results of WENO-JP FSM \cite{WENOJP,xiong}, and list their errors, convergence rates and the numbers of iterations. In all the numerical examples, $\epsilon$ in \eqref{epsilon} is taken as $10^{-6}$ unless otherwise specified. We use ``iter'' to indicates the number of iterations (noting that one iteration includes four alternating sweepings) in all the tables. The total number of grid points is assumed to be $N_x=N_y=N$. We take $\omega=0.7$ for Approach 1 and $\omega=0.8$ for Approach 2 in Example 1-6. While for Example 7 P-wave, we take $\omega=1.2$ for both Approach 1 and 2, and $\omega=0.9$ for both Approach 1 and 2 in Example 7 SV-wave. All the computations are implemented by using MATLAB 2019b on ThinkPad computer with 1.70 GHz Intel Core i5 processor and 4GB RAM.

\bigskip
\noindent \textbf{Example 1}. We solve the Eikonal equation with
$$f(x,y)=\frac{\pi}{2}\sqrt{\sin^{2}\Big(\pi+\frac{\pi}{2}x\Big)+\sin^{2}\Big(\pi+\frac{\pi}{2}y\Big)},$$
on the computational domain $[-1,1]^2$, with the inflow boundary $\Gamma=(0,0)$. The exact solution is given by
\begin{equation*}
\phi(x,y)=\cos\Big(\pi+\frac{\pi}{2}x\Big)+\cos\Big(\pi+\frac{\pi}{2}y\Big).
\end{equation*}

\begin{figure}[!h]
\begin{center}
	\includegraphics[width=7.6cm]{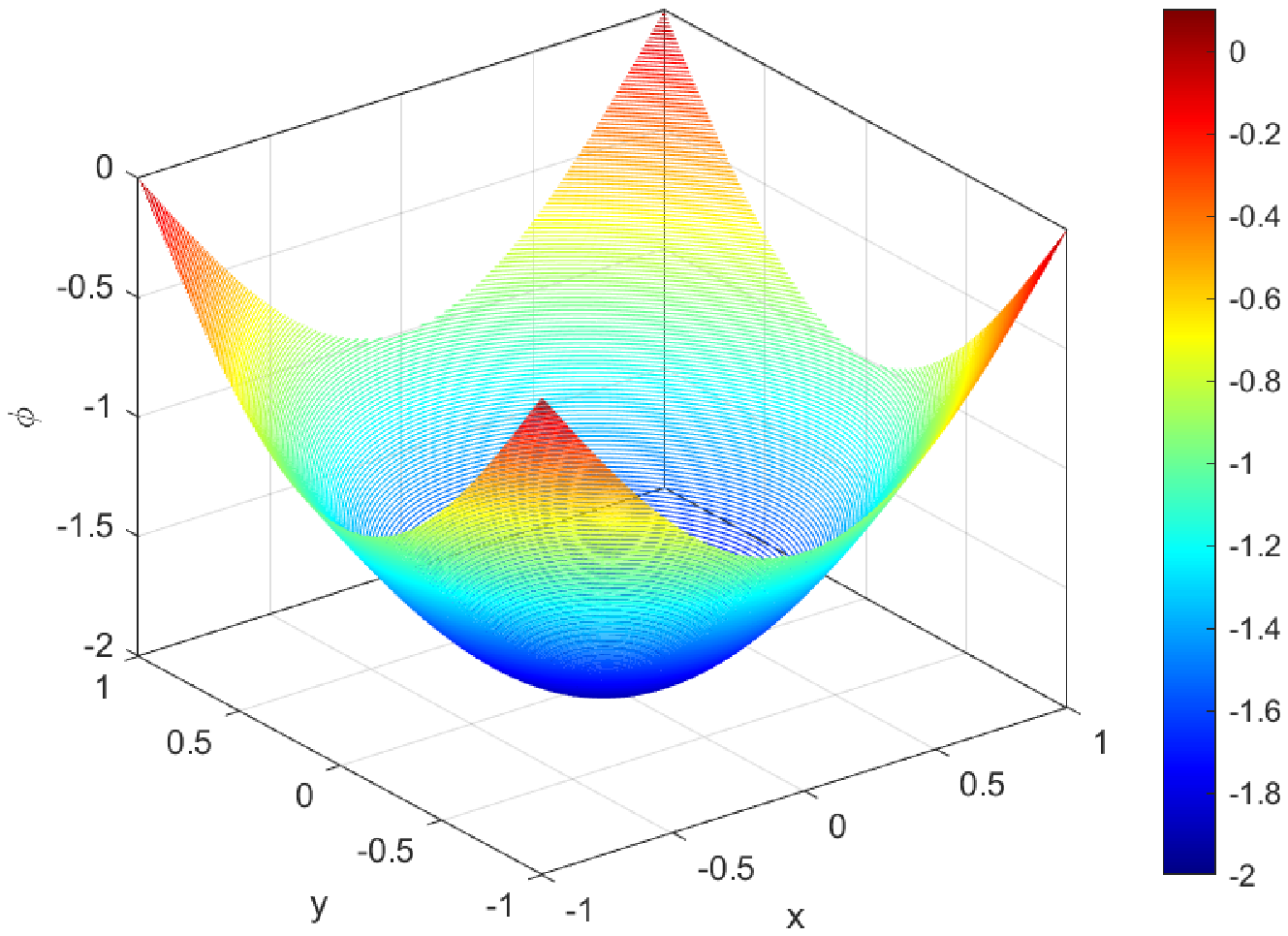}
	\includegraphics[width=7.6cm]{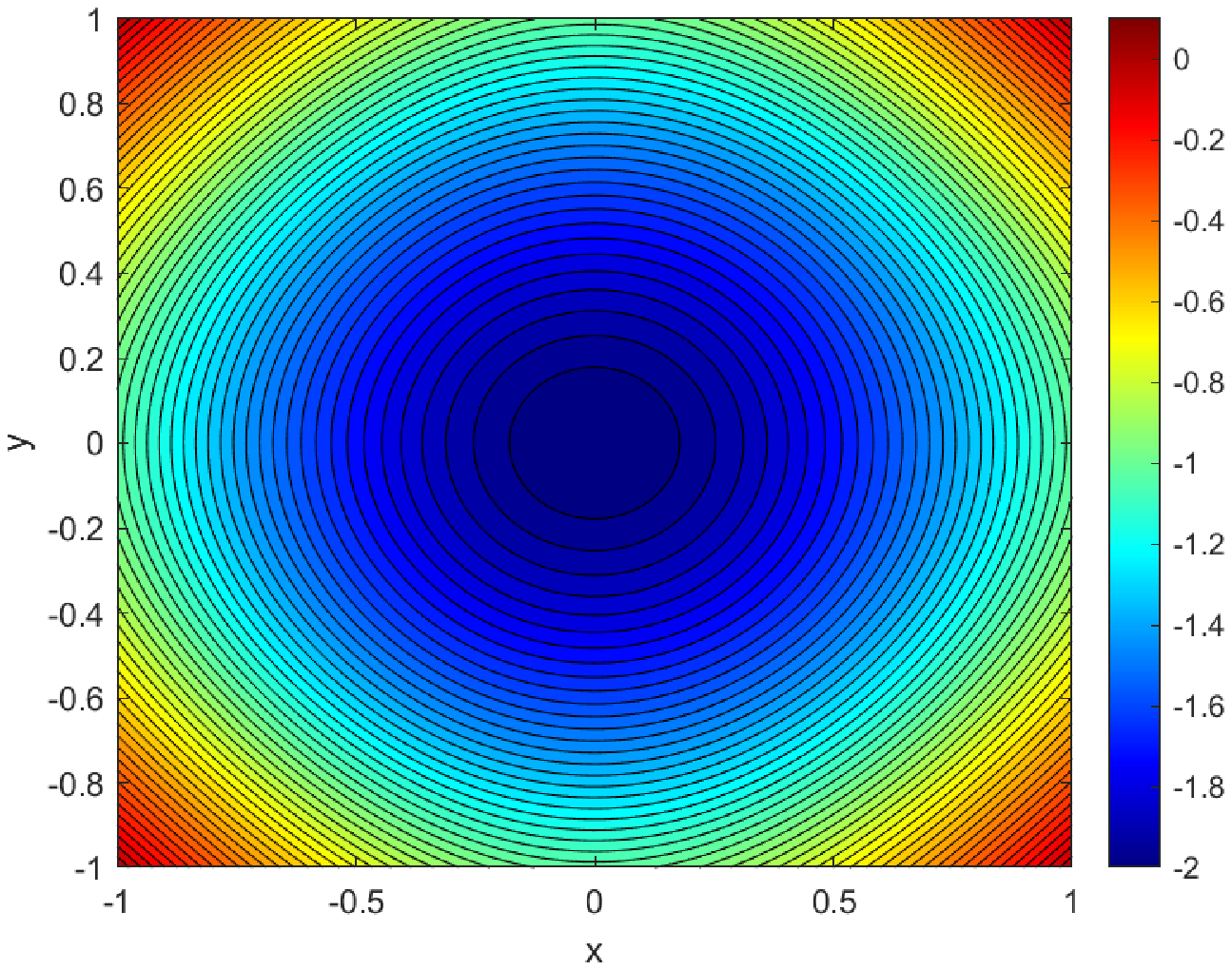}
\end{center}
\caption{Example 1. The numerical solution by Approach 1 on mesh $N=160$. Left: the 3D plot of numerical solution $\phi$; Right: the contour plot for $\phi$ .}\label{e1solutioncontour}
\end{figure}

\begin{figure}[!h]
\begin{center}
	\includegraphics[width=7.6cm]{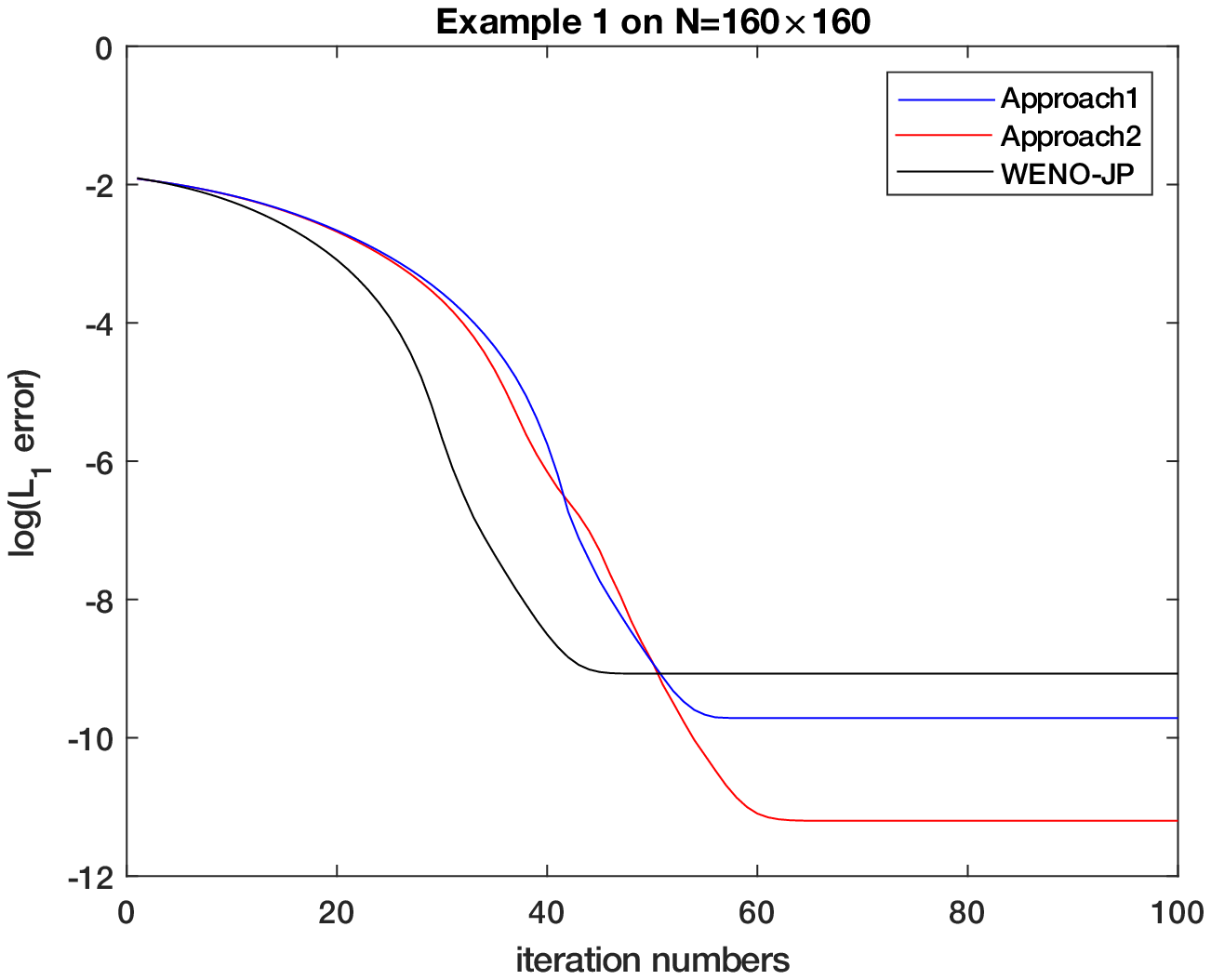}
	\includegraphics[width=7.6cm]{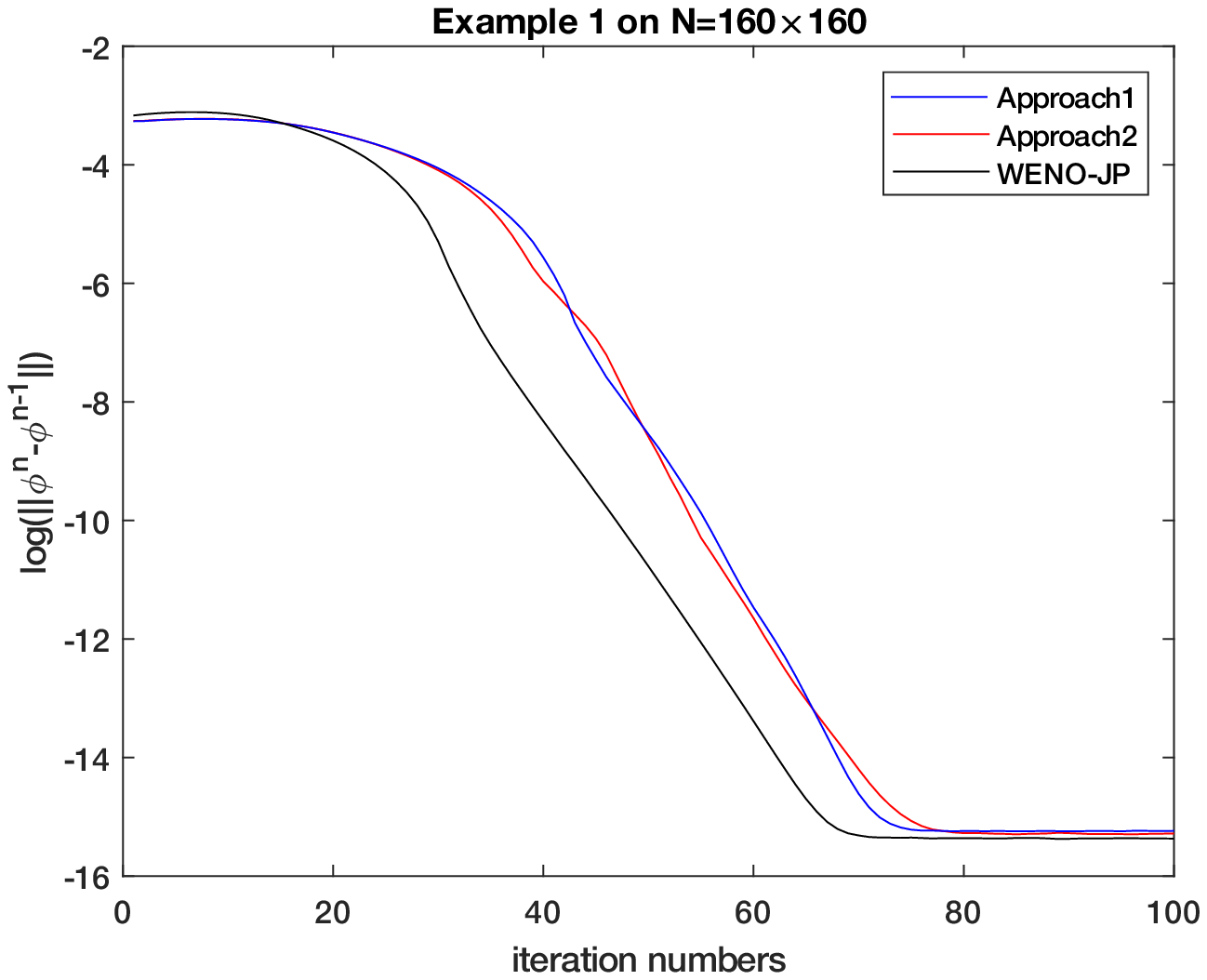}
\end{center}
\caption{Three methods solving Example 1. Left: number of iterations vs $log(L_{1}~error$); Right: number of iterations vs $log(||\phi^{n}-\phi^{n-1}||_{L_{1}}$).}\label{fig11}
\end{figure}

The Godunov numerical Hamiltonian \eqref{g1} is used. The picture of numerical solution by Approach 1 on mesh $N=160$ are presented in Figure \ref{e1solutioncontour}. The numerical errors and orders of convergence for these different schemes are provided in Table \ref{tab1} for comparison.
We can see that the errors of Approach 1 and 2 are smaller than that of the WENO-JP scheme on the same mesh size, although they have slightly more iteration numbers when the mesh is refined. We fix the mesh to be $N=160$, and run all three methods with $100$ iterations. The history of their numerical errors is shown in Figure \ref{fig11}, from which we can clearly observe that Approach 2 yields the smallest error (around $10^{-11}$), while the error of WENO-JP scheme stays at the level of $10^{-9}$. The error of Approach 1 stays between them. Although the error of WENO-JP scheme reaches $10^{-9}$ with the least number of iterations, its error will not decay further and stay at that level. We also plot the history of the quantity $\delta=||\phi^{new}-\phi^{old}||_{L_{1}}$ for comparison, where we can observe that all three methods converge to $\delta = 10^{-15}$ after the 80th iteration.

\begin{table}[h!]
\caption{Example 1. Comparison of the three methods: The errors of the numerical solution, the accuracy obtained and the number of iterations for convergence}\label{tab1}
	\begin{center}
		\begin{tabular}{|c|c|c|c|c|c|c|c|}
			\hline
Approach 1 & $L_{1}$~error&order&  $L_{\infty}$~error&order&iter\\ \hline
40 &3.21e-06 &- &2.61e-05 &- &40 \\ \hline
80 &3.80e-08 &6.40 &7.13e-07 &5.19 &46 \\ \hline
160 &1.92e-10 &7.62 &7.68e-09 &6.53 &69 \\ \hline
320 &2.18e-13 &9.78 &7.28e-12 &10.00 &111\\ \hline
Approach 2 & $L_{1}$~error&order&  $L_{\infty}$~error&order&iter\\ \hline
40 &2.38e-07 &- &4.09e-06 &- &41 \\ \hline
80 &8.21e-10 &8.18 &3.10e-08 &7.03 &50 \\ \hline
160 &6.28e-12 &7.02 &2.02e-11 &11.05 &70 \\ \hline
320 &2.75e-13 &4.51 &8.73e-12 &1.21 &111 \\ \hline
WENO-JP FSM & $L_{1}$~error&order&  $L_{\infty}$~error&order&iter\\ \hline
40 &1.54e-05 &- &9.83e-05 &- &47\\ \hline
80 &1.13e-07 &7.09 &1.60e-06 &5.94 &51\\ \hline
160 &8.45e-10 &7.06 &2.24e-08 &6.15 &63\\ \hline
320 &3.01e-12 &8.13 &3.29e-11 &9.40 &80\\ \hline
		\end{tabular}
	\end{center}
\end{table}

\bigskip
\noindent\textbf{Example 2}. We solve the Eikonal equation with $f(x,y)=1$. The computational domain is set as $[-1,1]^2$, and the inflow boundary $\Gamma$ is the circle with center at $(0,0)$ and radius $0.5$, that is
$$\Gamma=\left\{(x,y)|x^2+y^2=\frac{1}{4}\right\}.$$
The boundary condition is given as $\phi(x,y)=0$ on $\Gamma$. The exact solution is a distance function to the circle $\Gamma$, and it has a singularity at the center of the circle (due to the intersection of characteristic lines).

\begin{figure}[!h]
\begin{center}
	\includegraphics[width=7.6cm]{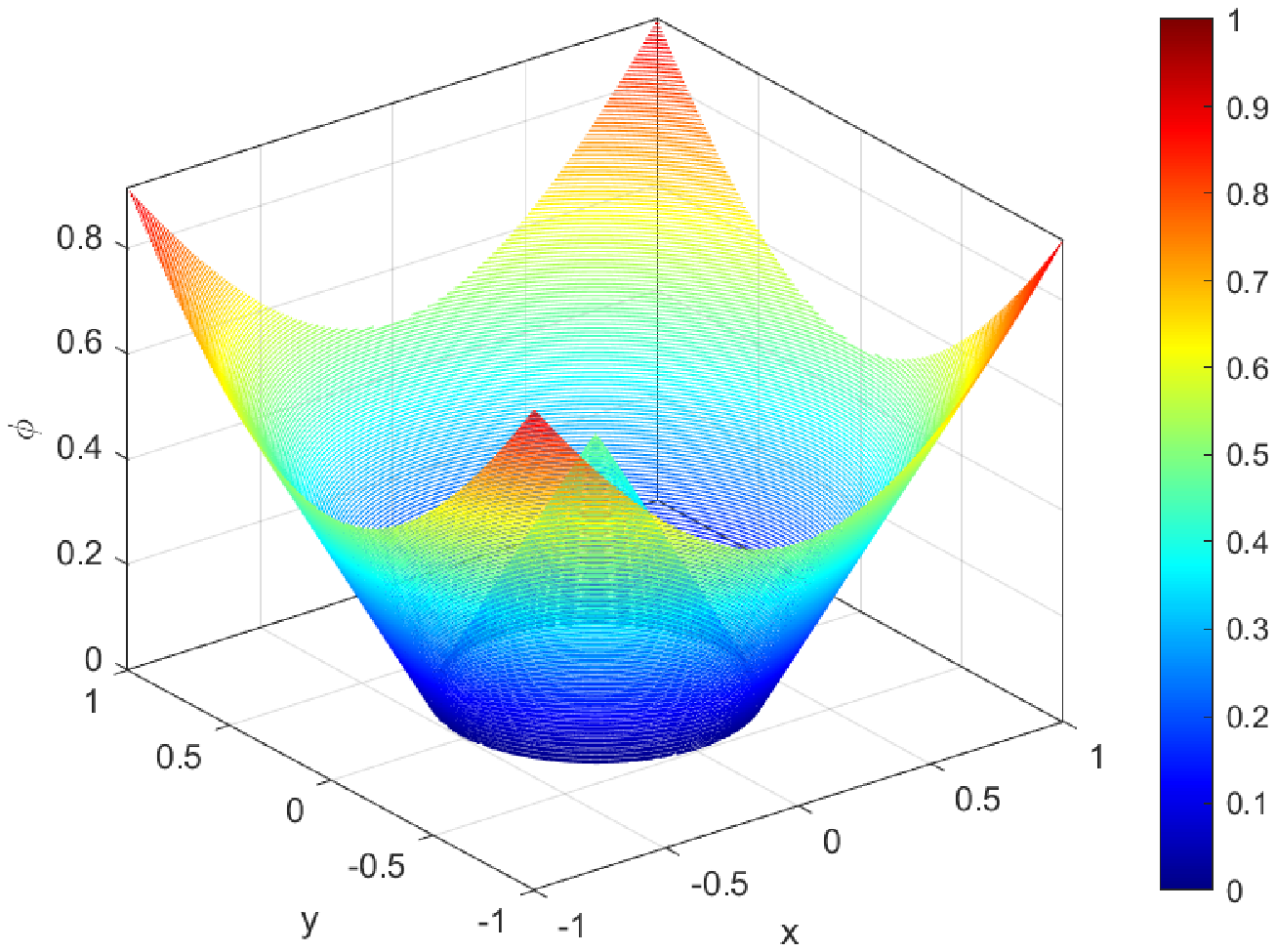}
	\includegraphics[width=7.6cm]{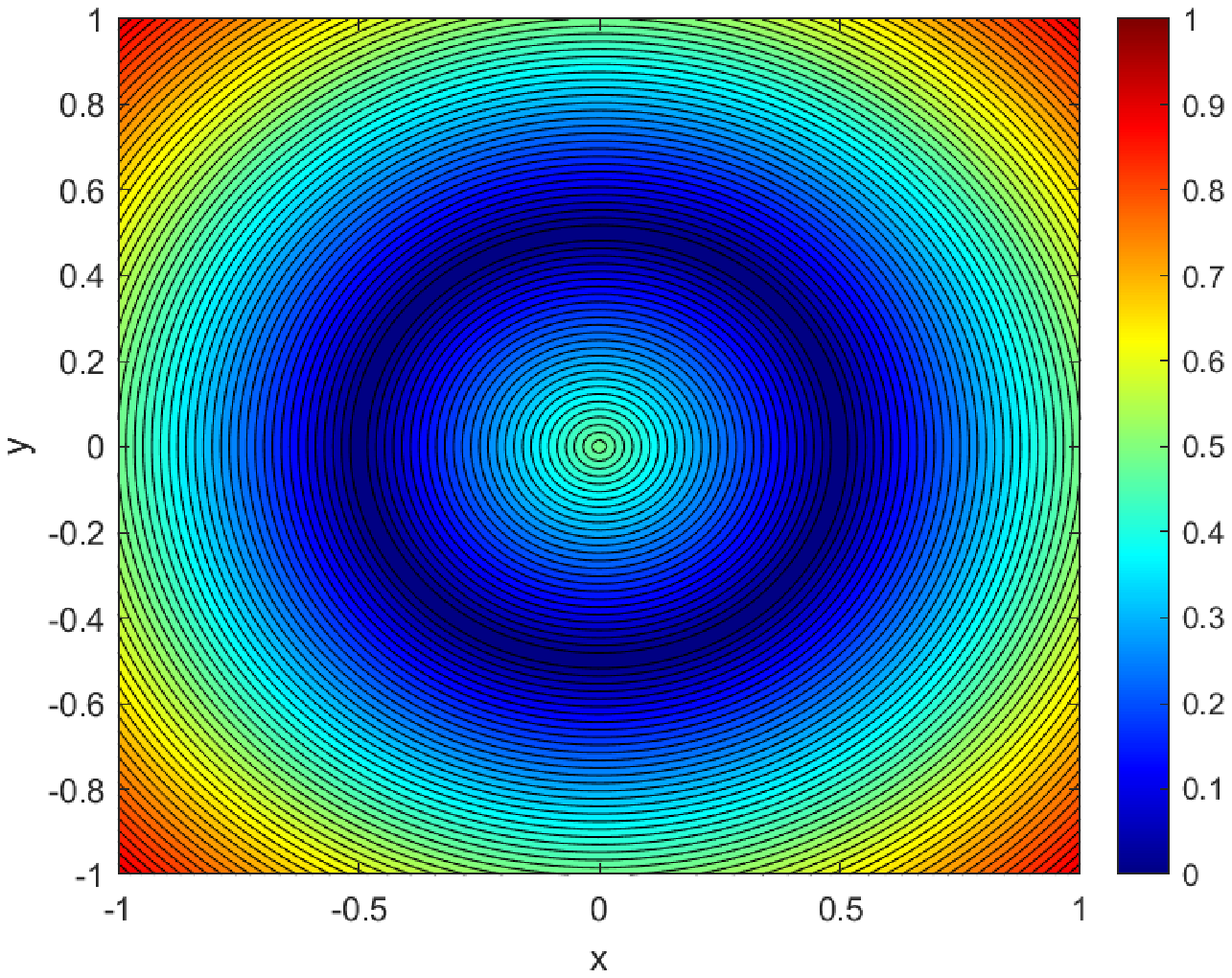}
\end{center}
\caption{Example 2. The numerical solution by Approach 1 on mesh $N=160$. Left: the 3D plot of numerical solution $\phi$; Right: the contour plot for $\phi$ .}\label{e2solutioncontour}
\end{figure}

\begin{figure}[!h]
\begin{center}
	\includegraphics[width=7.6cm]{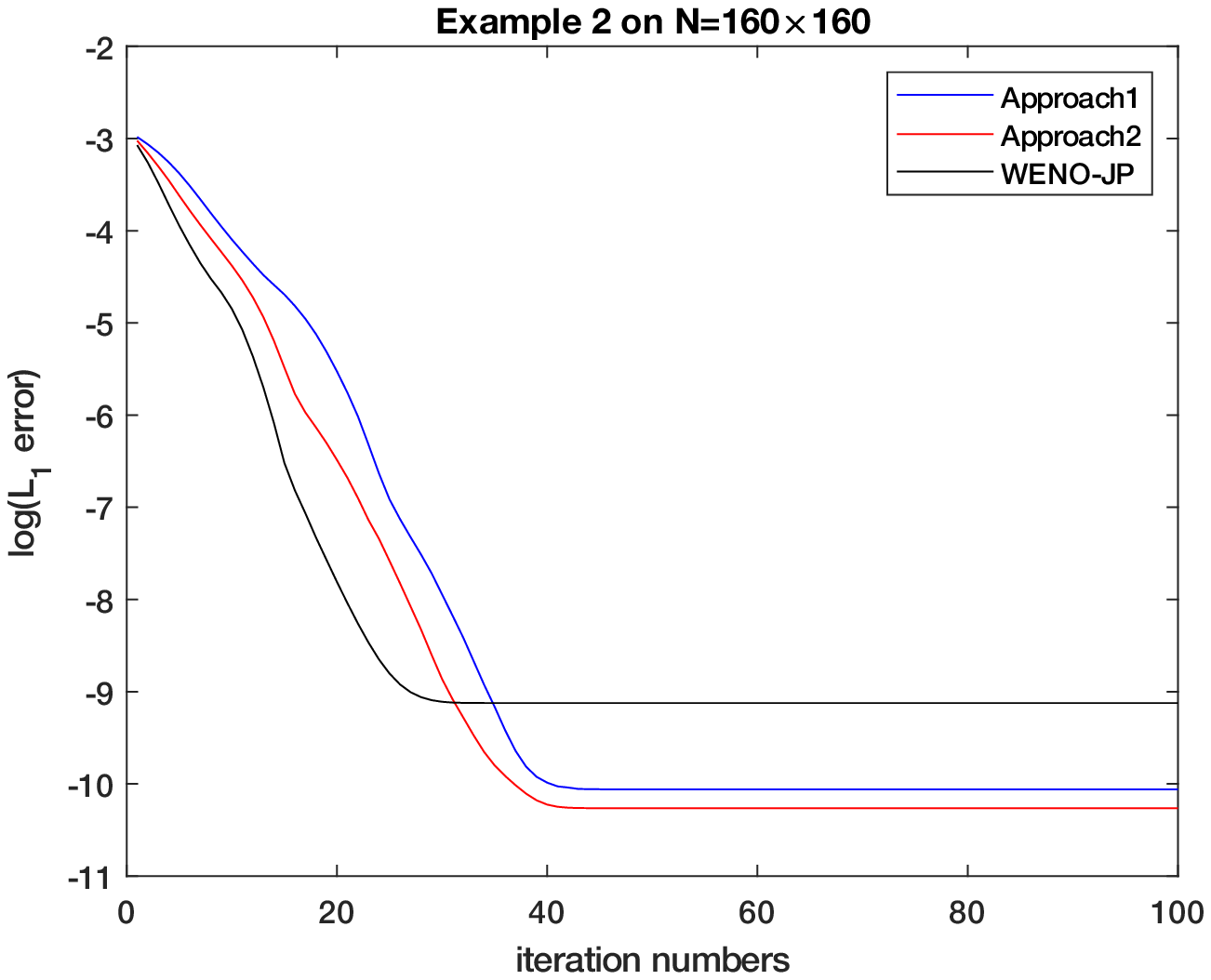}
	\includegraphics[width=7.6cm]{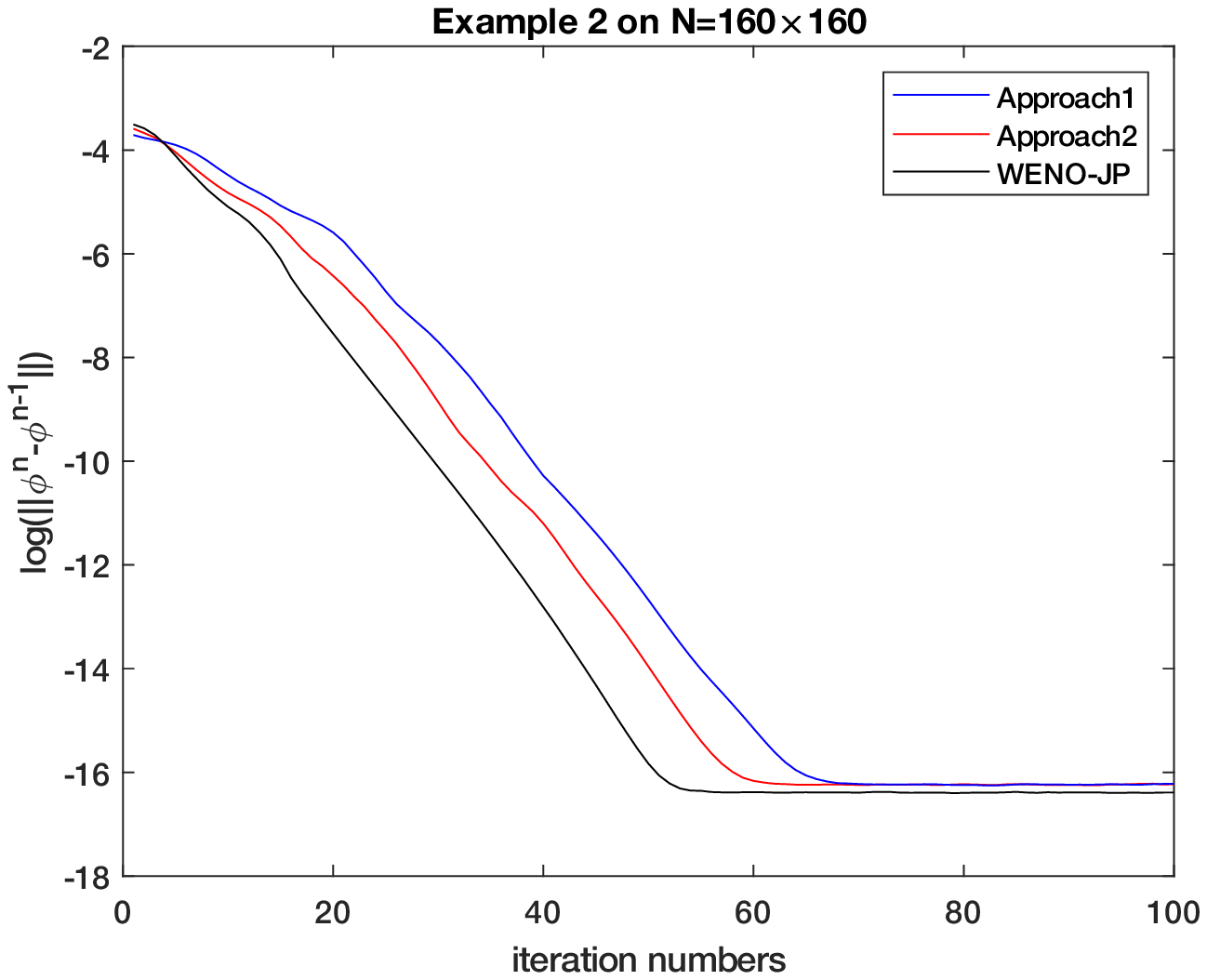}
\end{center}
\caption{Three methods solving Example 2. Left: number of iterations vs $log(L_{1}~error$); Right: number of iterations vs $log(||\phi^{n}-\phi^{n-1}||_{L_{1}}$).}\label{fig12}
\end{figure}

The Godunov numerical Hamiltonian is used, and the numerical errors are measured in the box $[-0.9,0.9]^2$ and outside the box $[-0.15,0.15]^2$, which aims to remove the influence of singularity and outflow boundary treatment. The picture of numerical solution by Approach 1 are presented in Figure \ref{e2solutioncontour}. The numerical errors and orders of convergence are listed in Table \ref{tab2}. Again, we can observe that the errors of the Approach 1 and 2 are smaller than that of the WENO-JP schemes. Moreover, all these three methods achieve the designed fifth order accuracy on this test example. The number of iterations of WENO-JP scheme is slightly smaller than the other two methods. We also fix the mesh to be $N=160$, and run all three methods with $100$ iterations. The history of their numerical errors is shown in Figure \ref{fig12}, from which we can clearly observe that Approach 2 yields the smallest error, WENO-JP method produces the largest error, and the error of Approach 1 stays between them. We plot the history of the quantity $\delta$ for comparison, where we can observe that all three methods converge to $\delta = 10^{-16}$ after the 70th iteration. 

\begin{table}[h!]
\caption{Example 2. Comparison of the three methods: The errors of the numerical solution, the accuracy obtained and the number of iterations for convergence}\label{tab2}
	\begin{center}
		\begin{tabular}{|c|c|c|c|c|c|c|c|}
			\hline
Approach 1 & $L_{1}$~error&order&  $L_{\infty}$~error&order&iter\\ \hline
40 &6.06e-07 &- &2.34e-05 &- &37\\ \hline
80 &1.17e-08 &5.69 &1.74e-06 &3.74 &44\\ \hline
160 &8.72e-11 &7.07 &2.20e-08 &6.30 &55\\ \hline
320 &1.78e-12 &5.60 &1.30e-10 &7.40 &81\\ \hline
Approach 2 & $L_{1}$~error&order&  $L_{\infty}$~error&order&iter\\ \hline
40 &1.92e-07 &- &1.09e-05 &- &33 \\ \hline
80 &1.66e-09 &6.85 &1.87e-07 &5.85 &38 \\ \hline
160 &5.45e-11 &4.92 &3.67e-09 &5.67 &51 \\ \hline
320 &2.01e-12 &4.75 &1.23e-10 &4.90 &73 \\ \hline
WENO-JP FSM & $L_{1}$~error&order&  $L_{\infty}$~error&order&iter\\ \hline
40 &5.18e-07 &- &3.66e-05 &- &35\\ \hline
80 &3.84e-08 &3.75 &3.40e-06 &3.42 &35\\ \hline
160 &7.52e-10 &5.67 &2.72e-07 &3.64 &44\\ \hline
320 &1.65e-11 &5.50 &1.09e-09 &7.95 &59\\ \hline
		\end{tabular}
	\end{center}
\end{table}

\bigskip
\noindent \textbf{Example 3}. We solve the Eikonal equation with $f(x,y)=1$. The computational domain is set as $[-3,3]^2$, and the inflow boundary $\Gamma$ consists of two circles of equal radius $0.5$ with the centers located at $(-1,0)$ and $(\sqrt{1.5},0)$, respectively, that is
$$\Gamma=\left\{(x,y)|(x+1)^2+y^2=\frac{1}{4}\quad or\quad (x-\sqrt{1.5})^2+y^2=\frac{1}{4}\right\}.$$
The exact solution is a distance function to the inflow boundary $\Gamma$, containing the singularities at the center of each circle and the line $x=0.5(\sqrt{1.5}-1)$ that is of equal distance to two circle centers.

Again, the Godunov numerical Hamiltonian is used. We measure the numerical errors within the box of $[-2.85,2.85]^2$, which also excludes the boxes $[-1.15,-0.85]\times[-0.15,0.15]$, $[\sqrt{1.5}-0.15,\sqrt{1.5}+0.15]\times[-0.15,0.15]$
and $[\sqrt{0.375}-0.65,\sqrt{0.375}-0.35]\times[-2.85,2.85]$. These excluded boxes contain two centers of $\Gamma$ and the singular line.

\begin{figure}[!h]
\begin{center}
	\includegraphics[width=7.6cm]{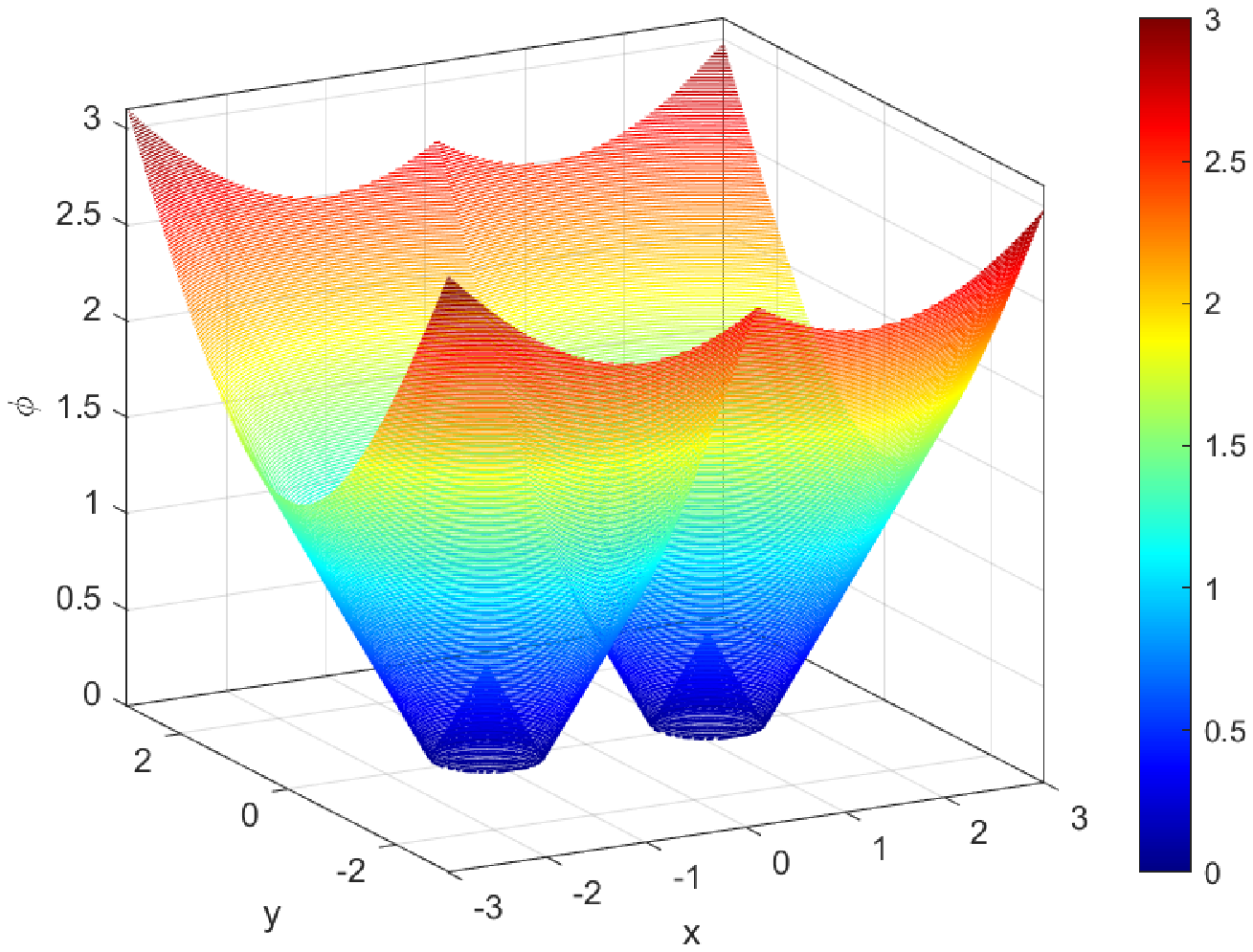}
	\includegraphics[width=7.6cm]{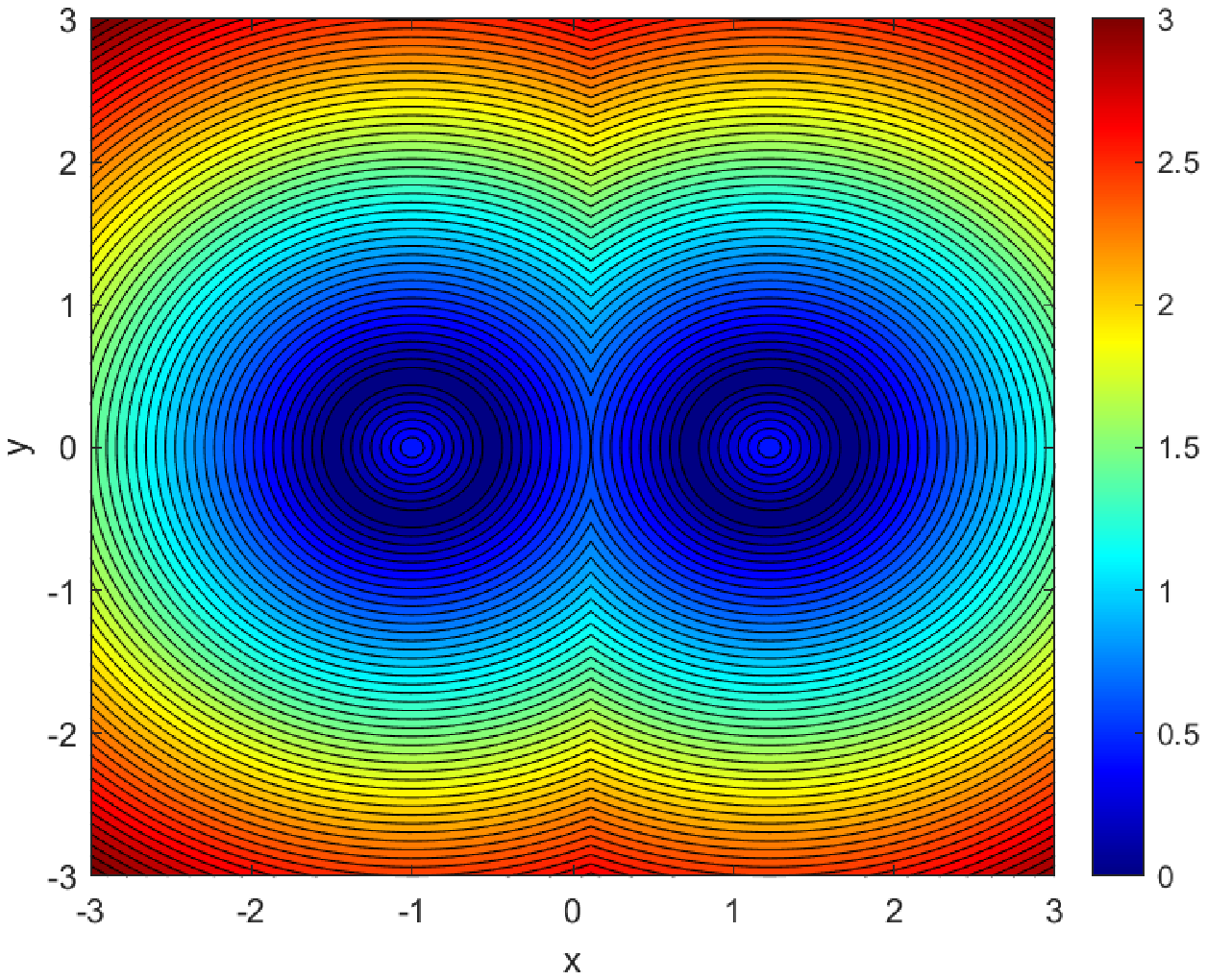}
\end{center}
\caption{Example 3. The numerical solution by Approach 1 on mesh $N=160$. Left: the 3D plot of numerical solution $\phi$; Right: the contour plot for $\phi$ .}\label{e3solutioncontour}
\end{figure}

The Figure \ref{e3solutioncontour} shows that the numerical solution by Approach 1 on mesh $N=160$. The numerical errors and orders of convergence are shown in Table \ref{tabl3}. All three methods achieved the designed high order accuracy. Although WENO-JP method needs the least number of iterations, its numerical errors are also the largest among all three methods. On the other hand, Approach 1 and Approach 2 yield similar numerical errors and the number of iterations. 

\begin{table}[h!]
\caption{Example 3. Comparison of the three methods: The errors of the numerical solution, the accuracy obtained and the number of iterations for convergence}\label{tabl3}
	\begin{center}
		\begin{tabular}{|c|c|c|c|c|c|c|c|}
			\hline
Approach 1 & $L_{1}$~error&order&  $L_{\infty}$~error&order&iter\\ \hline
80 &7.13e-07 &- &1.04e-04 &- &58\\ \hline
160 &1.09e-07 &2.70 &9.45e-06 &3.47 &55\\ \hline
320 &2.97e-10 &8.52 &3.86e-07 &4.61 &85\\ \hline
640 &6.82e-12 &5.44 &3.22e-9 &6.90 &139\\ \hline
Approach 2 & $L_{1}$~error&order&  $L_{\infty}$~error&order&iter\\ \hline
80 &1.34e-06 &- &8.57e-05 &- &38\\ \hline
160 &2.04e-07 &2.72 &1.86e-05 &2.20 &55\\ \hline
320 &7.53e-10 &8.08 &4.45e-07 &5.38 &93\\ \hline
640 &6.26e-12 &6.91 &7.86e-10 &9.14 &150\\ \hline
WENO-JP FSM & $L_{1}$~error&order&  $L_{\infty}$~error&order&iter\\ \hline
80 &1.12e-06 &- &1.22e-04 &- &55\\ \hline
160 &1.19e-07 &3.22 &6.73e-06 &4.18 &45\\ \hline
320 &2.19e-09 &5.76 &4.30e-07 &4.02 &52\\ \hline
640 &5.90e-11 &5.21 &8.87e-09 &5.59 &78\\ \hline
		\end{tabular}
	\end{center}
\end{table}

\bigskip
\noindent \textbf{Example 4}. We solve the Eikonal equation with $f(x,y)=1$. The computational domain is set as $[-1,1]^2$, and the inflow boundary is given by $\Gamma=(0,0)$. The exact solution for this problem is a distance function to $\Gamma$,
and it contains a singularity at $\Gamma$.

The Godunov numerical Hamiltonian is used. Due to the singularity, we follow the setup in \cite{xiong}, and pre-assign the exact solution in a small box with length $0.3$ around the source point. Numerical errors and orders are listed in Table \ref{tabl4}. Again, we can observe that fifth order accuracy can be obtained for all schemes, and Approach 1 and 2 yield smaller numerical errors than the WENO-JP scheme. 

\begin{table}[h!]
\caption{Example 4. Comparison of the three methods: The errors of the numerical solution, the accuracy obtained and the number of iterations for convergence}\label{tabl4}
	\begin{center}
		\begin{tabular}{|c|c|c|c|c|c|c|c|}
			\hline
Approach 1 & $L_{1}$~error&order&  $L_{\infty}$~error&order&iter\\ \hline
40 &3.11e-07 &- &4.60e-06 &- &41 \\ \hline
80 &6.95e-09 &5.48 &1.62e-07 &4.82 &47 \\ \hline
160 &1.19e-10 &5.85 &1.78e-09 &6.50 &62\\ \hline
320 &3.26e-12 &5.20 &8.67e-12 &7.68 &95 \\ \hline
Approach 2 & $L_{1}$~error&order&  $L_{\infty}$~error&order&iter\\ \hline
40 &8.42e-08 &- &7.14e-07 &- &32 \\ \hline
80 &2.88e-09 &4.87 &1.14e-08 &5.96 &40 \\ \hline
160 &1.04e-10 &4.78 &2.55e-10 &5.47 &53\\ \hline
320 &3.57e-12 &4.87 &9.84e-12 &4.70 &77 \\ \hline
WENO-JP FSM & $L_{1}$~error&order&  $L_{\infty}$~error&order&iter\\ \hline
40 &6.17e-06 &- &9.31e-05 &- &35\\ \hline
80 &4.85e-07 &5.03 &7.34e-06 &3.66 &40\\ \hline
160 &8.92e-09 &5.76 &2.05e-07 &5.16 &50\\ \hline
320 &2.16e-10 &5.36 &2.68e-09 &6.25 &66\\ \hline
		\end{tabular}
	\end{center}
\end{table}

\bigskip
\noindent\textbf{Example 5}. We solve the Eikonal equation with $f(x,y)=1$. The computational domain is set as $[-1,1]^2$, and the inflow boundary
$\Gamma$ is a sector of three quarters of the circle centered at $(0,0)$ with radius $0.5$, closed with the x-axis and y-axis in the first quadrant, which can be described as
$$\Gamma=\left\{(x,y):\sqrt{x^2+y^2=0.5},~\mathrm{if}~x<0,y<0\}\cup\{(x,0):0\leq x\leq0.5\}\cup\{(0,y):0\leq y\leq0.5\right\}.$$
The exact solution is still the distance function to $\Gamma$. Singularities appear at the two corners in $\Gamma$, which give rise to both shock and rarefaction wave in the solution.

\begin{figure}[!h]
\begin{center}
	\includegraphics[width=7.6cm]{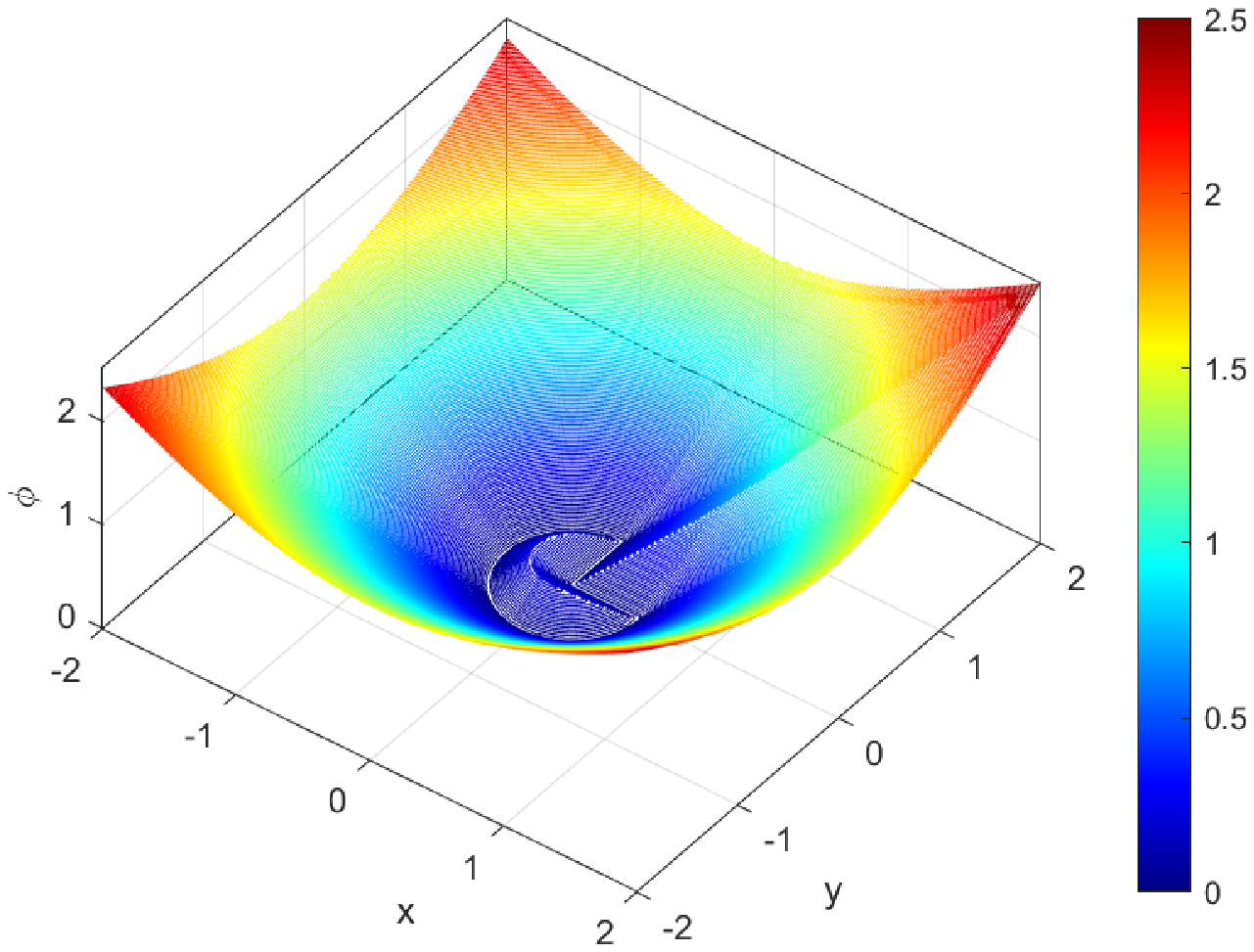}
	\includegraphics[width=7.6cm]{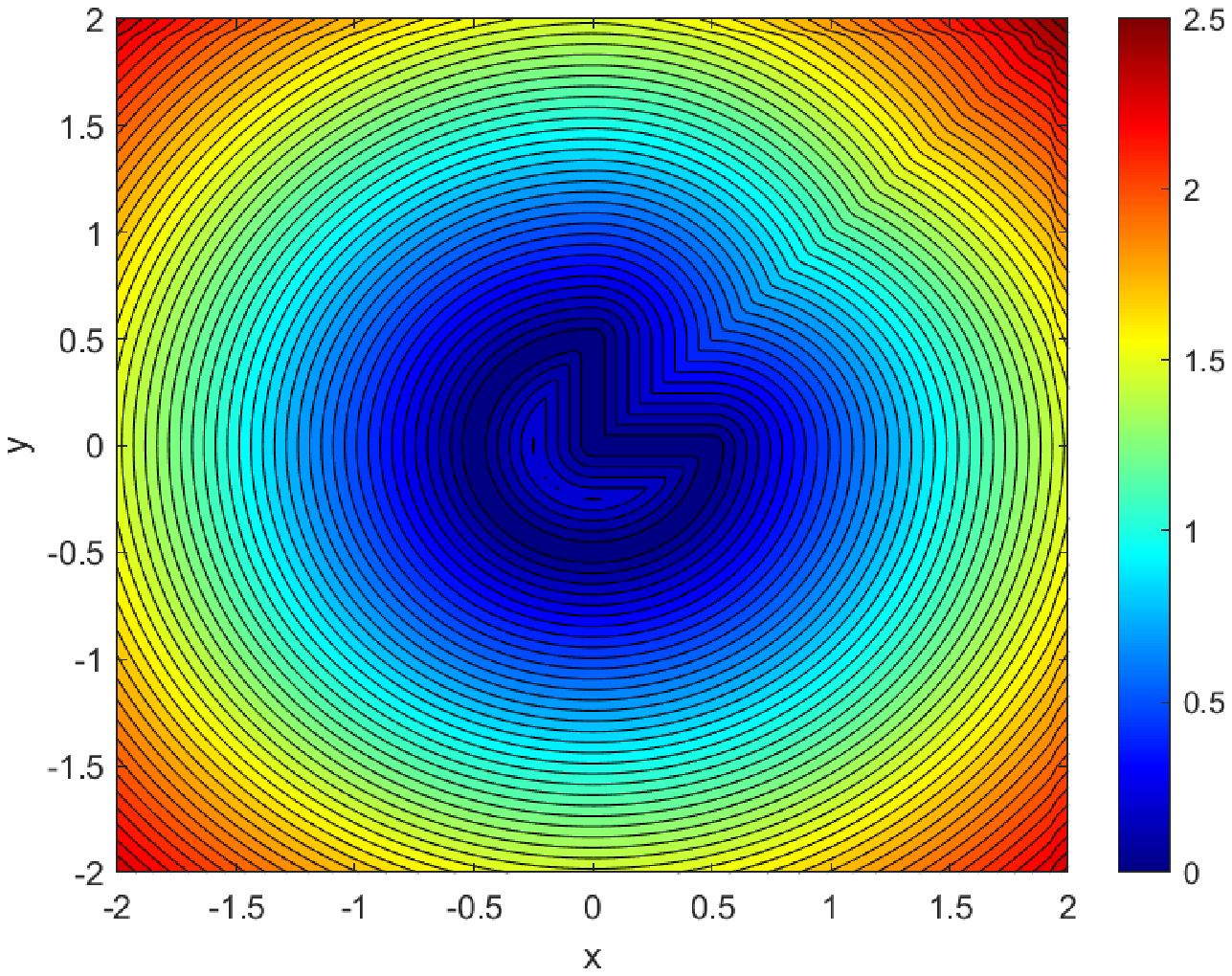}
\end{center}
\caption{Example 5. The numerical solution by Approach 1 on mesh $N=160$. Left: the 3D plot of numerical solution $\phi$; Right: the contour plot for $\phi$ .}\label{e5solutioncontour}
\end{figure}

\begin{figure}[!h]
\begin{center}
	\includegraphics[width=7.6cm]{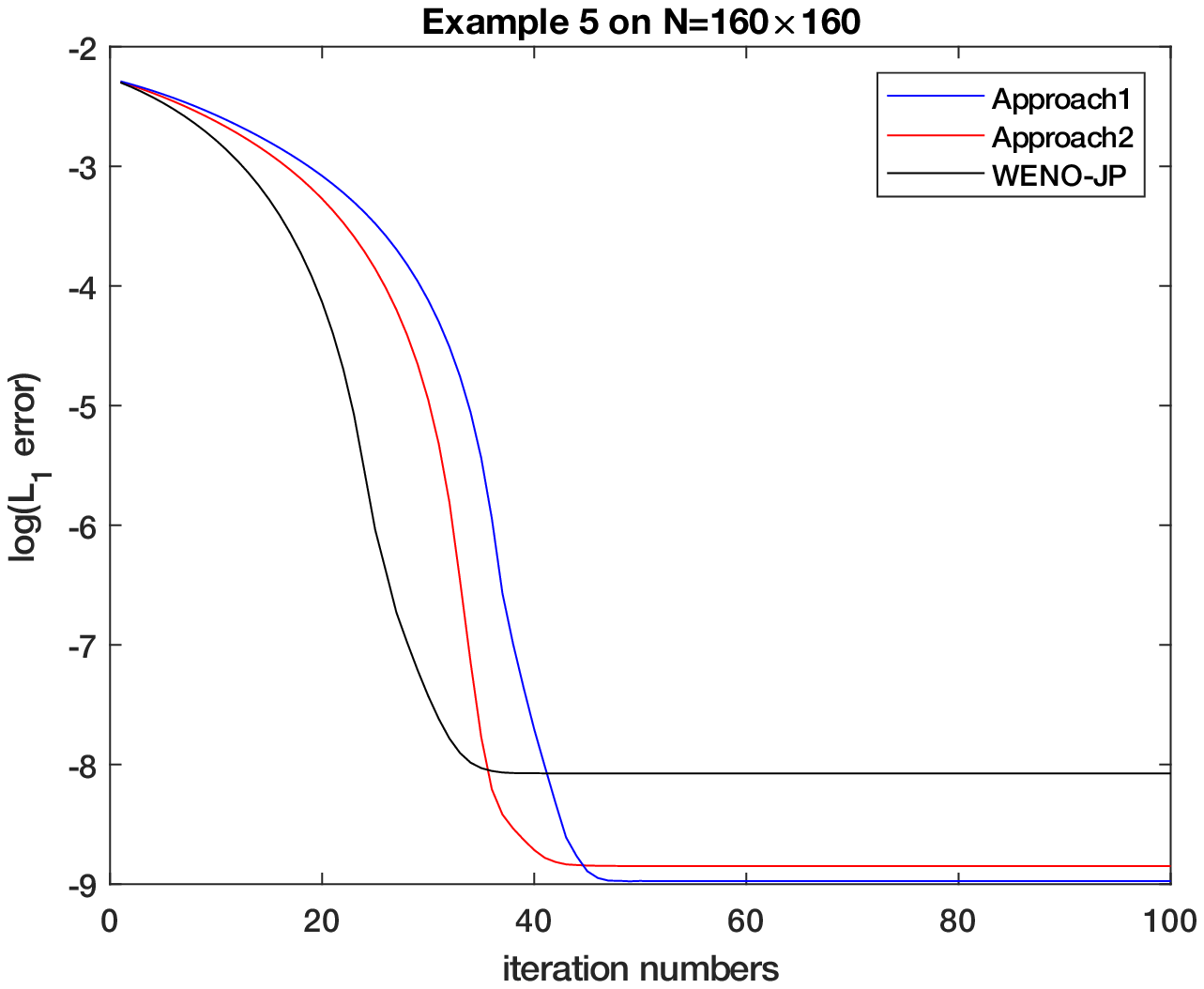}
	\includegraphics[width=7.6cm]{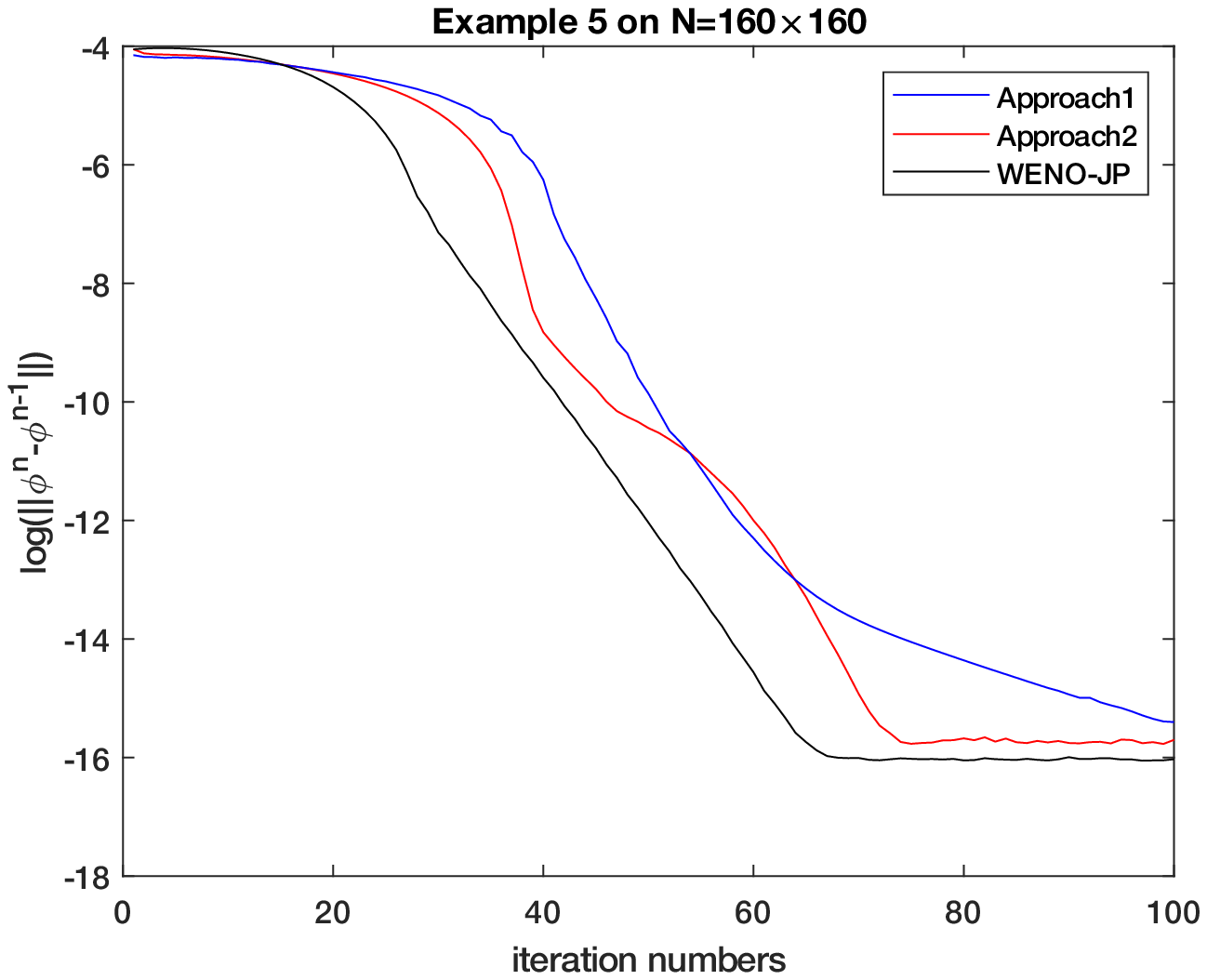}
\end{center}
\caption{Three methods solving Example 5. Left: number of iterations vs $log(L_{1}~error$); Right: number of iterations vs $log(||\phi^{n}-\phi^{n-1}||_{L_{1}}$).}\label{fig13}
\end{figure}

The Godunov numerical Hamiltonian is used. We measure the errors in smooth regions inside the box of $[-1.9,1.9]^2$ with $x\leq 0$ or $y\leq 0$, and outside the box $[-0.5,0.5]^2$. The picture of numerical solution are presented in Figure \ref{e5solutioncontour}. Numerical errors and orders of convergence are listed in Table \ref{tabl5}. Again, the fifth order accuracy can be obtained in the smooth regions for all three methods, and Approach 1 and Approach 2 yield smaller numerical errors on the same mesh size. We also fix the mesh to be $N=160$, and run all three methods with $100$ iterations. The history of their numerical errors is shown in Figure \ref{fig13}, from which we can clearly observe that Approach 1 yields the smallest error this time. WENO-JP method produces the largest error, and the error of Approach 2 stays between them. We plot the history of the quantity $\delta$ for comparison, where we can observe that all three methods converge to $\delta = 10^{-16}$.  

\begin{table}[h!]
\caption{Example 5. Comparison of the three methods: The errors of the numerical solution, the accuracy obtained and the number of iterations for convergence}\label{tabl5}
	\begin{center}
		\begin{tabular}{|c|c|c|c|c|c|c|c|}
			\hline
Approach 1 & $L_{1}$~error&order&  $L_{\infty}$~error&order&iter\\ \hline
40 &1.12e-06 &- &1.61e-05 &- &71 \\ \hline
80 &4.28e-08 &4.71 &1.08e-06 &3.89 &59 \\ \hline
160 &1.05e-09 &5.33 &3.93e-08 &4.78 &75 \\ \hline
320 &2.11e-11 &5.64 &4.44e-10 &6.46 &104 \\ \hline
Approach 2 & $L_{1}$~error&order&  $L_{\infty}$~error&order&iter\\ \hline
40 &8.84e-07 &- &1.52e-05 &- &33 \\ \hline
80 &3.04e-08 &4.86 &9.49e-07 &4.00 &41 \\ \hline
160 &1.07e-09 &4.82 &4.83e-08 &4.29 &65 \\ \hline
320 &2.66e-11 &5.33 &1.21e-09 &5.31 &108 \\ \hline
WENO-JP FSM & $L_{1}$~error&order&  $L_{\infty}$~error&order&iter\\ \hline
40 &3.46e-06 &- &2.68e-05 &- &35\\ \hline
80 &2.41e-07 &3.84 &1.10e-06 &4.60 &44\\ \hline
160 &8.43e-09 &4.84 &4.83e-08 &4.51 &58\\ \hline
320 &1.67e-10 &5.65 &9.94e-10 &5.60 &74\\ \hline
		\end{tabular}
	\end{center}
\end{table}

\bigskip
 \noindent\textbf{Example 6}. We solve the Eikonal equation with
 $$f(x,y)=2\pi\sqrt{[\cos(2\pi x)\sin(2\pi y)]^2+[\sin(2\pi x)\cos(2\pi y)]^2}.$$
The computational domain is set as $\Omega=[0,1]^2$, and the inflow boundary condition is given by $\Gamma=\{(\frac{1}{4},\frac{1}{4}),(\frac{3}{4},\frac{3}{4}),(\frac{1}{4},\frac{3}{4}),(\frac{3}{4},\frac{1}{4}),(\frac{1}{2},\frac{1}{2})\},$ consisting of five isolated points.
$\phi(x,y)=0$ is prescribed at the boundary of the unit square. The exact solution of this problem is the shape function \cite{zhang2006}. Two cases are considered here, based on different boundary conditions.
\indent \emph{Case a:}
$$g\left(\frac{1}{4},\frac{1}{4}\right)=g\left(\frac{3}{4},\frac{3}{4}\right)=1,\quad g\left(\frac{1}{4},\frac{3}{4}\right)=g\left(\frac{3}{4},\frac{1}{4}\right)=-1,\quad
g\left(\frac{1}{2},\frac{1}{2}\right)=0,$$
with the exact solution being
$$\phi(x,y)=\sin(2\pi x)\sin(2\pi y);$$
\indent\emph{Case b:}
 $$g\left(\frac{1}{4},\frac{1}{4}\right)=g\left(\frac{3}{4},\frac{3}{4}\right)=g\left(\frac{1}{4},\frac{3}{4}\right)=g\left(\frac{3}{4},\frac{1}{4}\right)=1,\quad
g\left(\frac{1}{2},\frac{1}{2}\right)=2,$$
with the exact solution being
\begin{equation*}
  \phi(x,y)=\begin{cases}
\max(|\sin(2\pi x)\sin(2\pi y)|,1+\cos(2\pi x)\cos(2\pi y)),
&\mathrm{if}~|x+y-1|<\frac{1}{2}~\mathrm{and}~|x-y|<\frac{1}{2},\\
|\sin(2\pi x)\sin(2\pi y)|,&\mathrm{otherwise},
\end{cases}
\end{equation*}
which is not smooth.

Due to the singularity of these point sources, the exact solutions are placed in a small box with a length $2h$ around these isolated points in both test cases. The Godunov numerical Hamiltonian is used in this test.

\begin{figure}[!h]
\begin{center}
	\includegraphics[width=7.6cm]{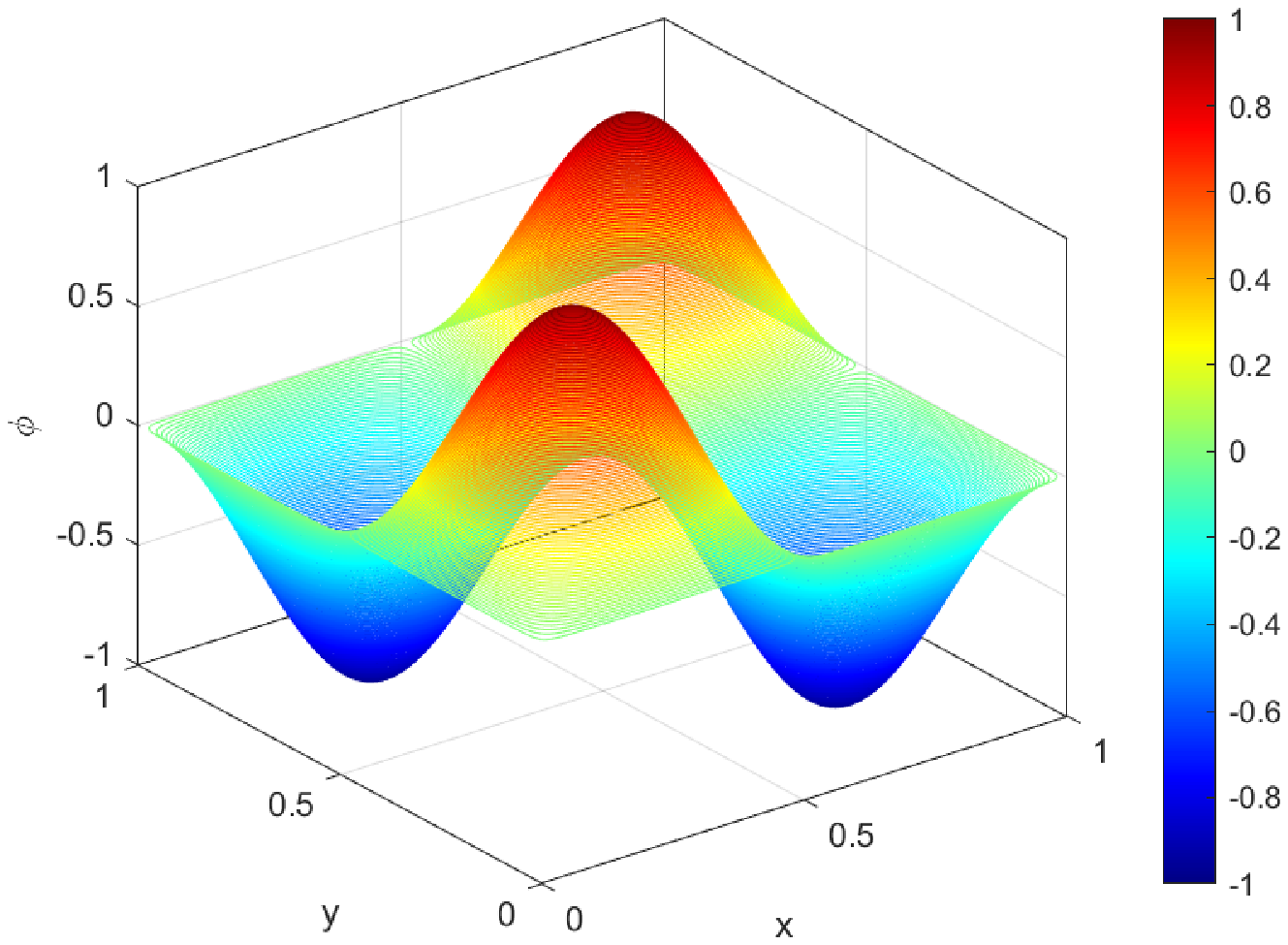}
	\includegraphics[width=7.6cm]{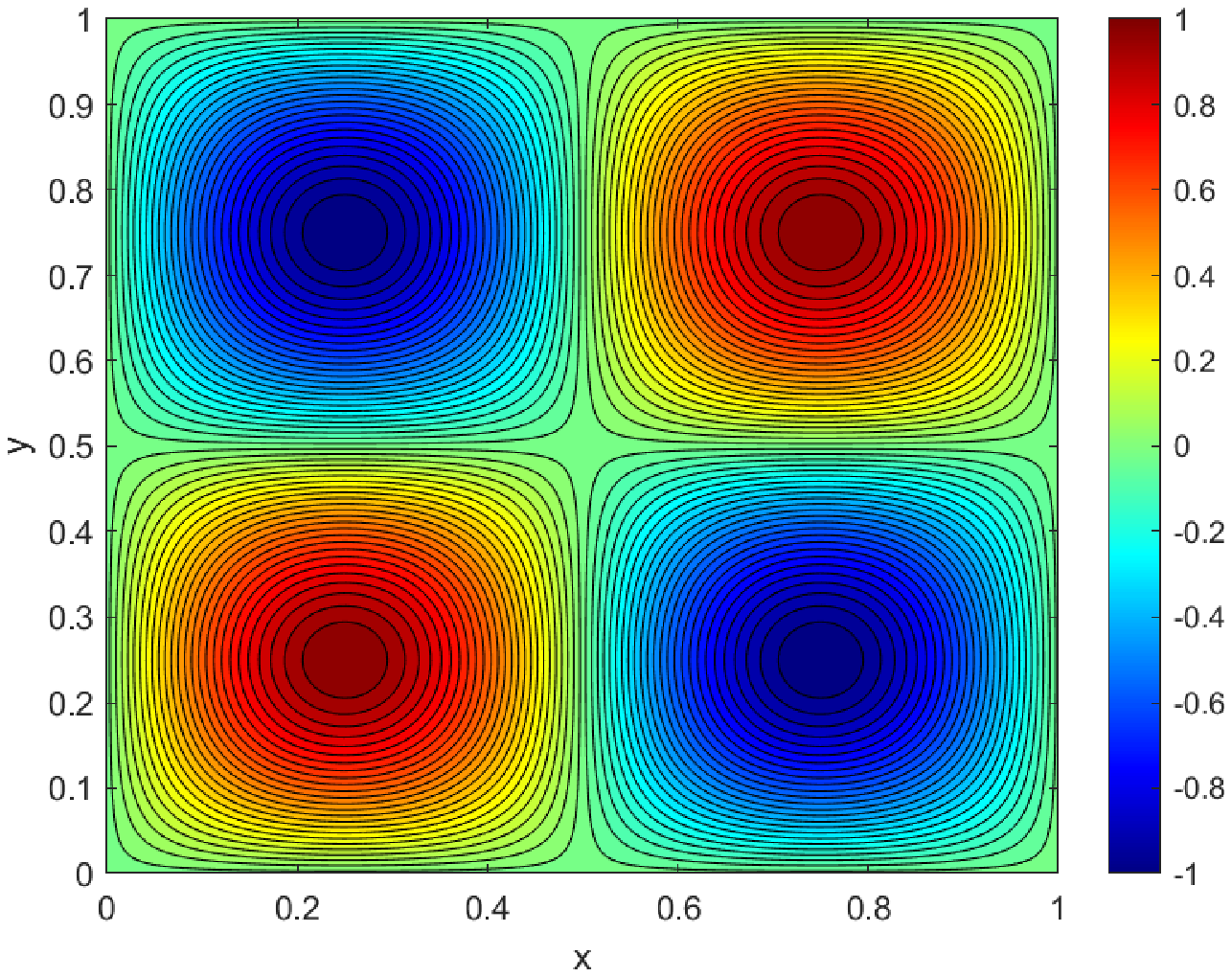}
\end{center}
\caption{Example 6 (case a). The numerical solution by Approach 1 on mesh $N=160$. Left: the 3D plot of numerical solution $\phi$; Right: the contour plot for $\phi$ .}\label{e61solution}
\end{figure}

\begin{figure}[H]
\begin{center}
	\includegraphics[width=7.6cm]{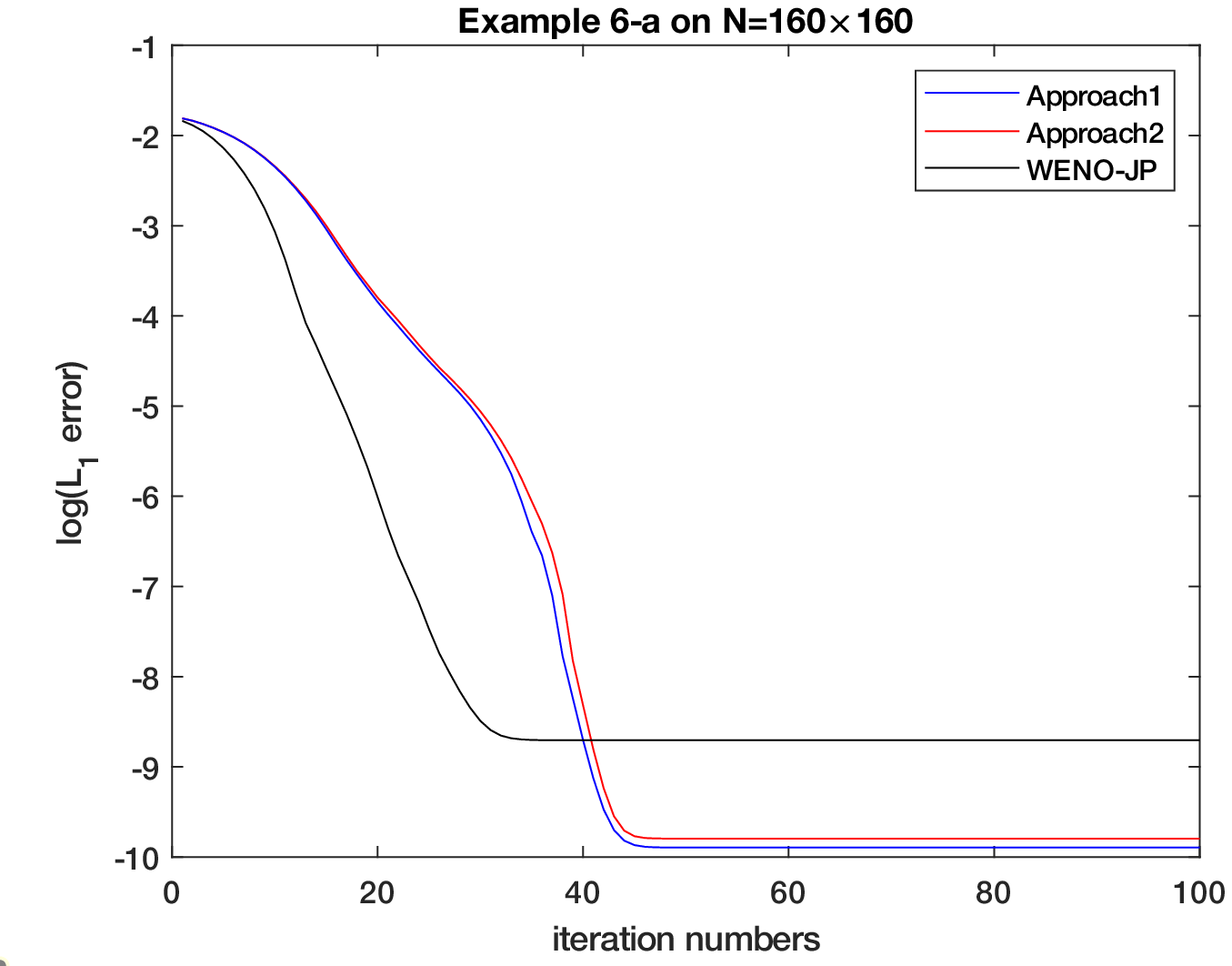}
	\includegraphics[width=7.6cm]{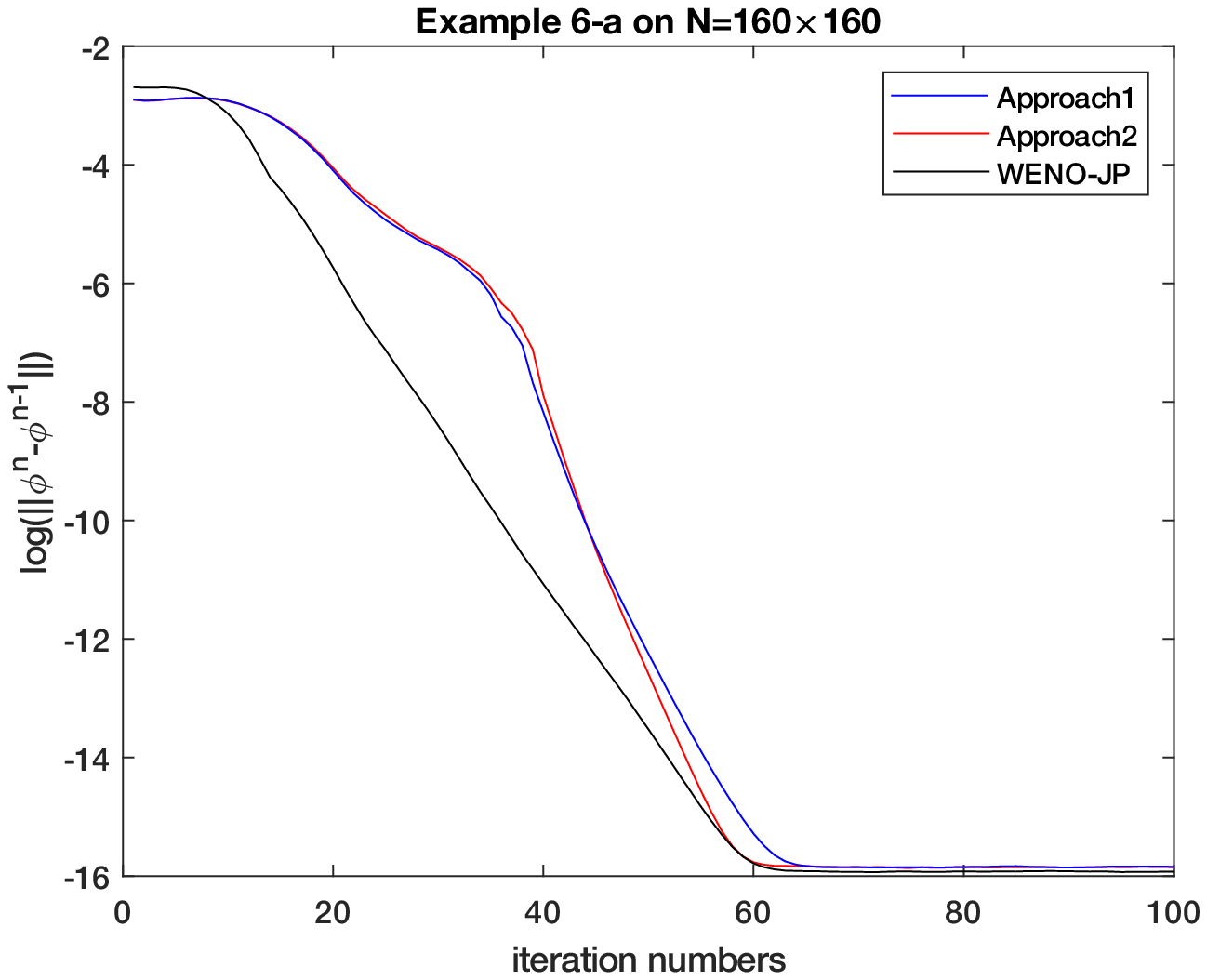}
\end{center}
\caption{Three methods solving Example 6-a. Left: number of iterations vs $log(L_{1}~error$); Right: number of iterations vs $log(||\phi^{n}-\phi^{n-1}||_{L_{1}}$).}\label{fig14}
\end{figure}

In \emph{case a}, The Figure \ref{e61solution} shows that numerical solution in case a. The numerical errors and orders of convergence of three methods are listed in Table \ref{tabl61}. We can see that the fifth order accuracy can be obtained, and the errors of the Approach 1 and 2 are smaller than that of the WENO-JP schemes. For this example, we would comment that, the number of iterations depends on the parameter $\epsilon$ in \eqref{epsilon} in order to achieve the desired fifth order.  We also fix the mesh to be $N=160$, and run all three methods with $100$ iterations. The history of their numerical errors is shown in Figure \ref{fig14}, from which we can clearly observe that Approach 1 and Approach 2 yield the smallest error, and WENO-JP method produces the largest error. We plot the history of the quantity $\delta$ for comparison, where we can observe that all three methods converge to $\delta = 10^{-16}$ after the 65th iteration. 

For \emph{case b}, we modified the convergence criteria to $\delta<10^{-12}$ since the solution is not smooth. The numerical errors and orders of convergence are listed in Table \ref{tabl62}. Due to the non-smoothness of the exact solution, we can only achieve second order accuracy. 
\begin{remark}
Similar to the observations in \cite{xiong,me}, we would like to comment that if we take a fixed $\epsilon$, e.g., $\epsilon=10^{-3}$ or $\epsilon=10^{-6}$, the three methods may either lose order or even blow up during the mesh refinement.
\end{remark}
\begin{table}[h!]
\caption{Example 6 case a. Comparison of the three methods: The errors of the numerical solution, the accuracy obtained and the number of iterations for convergence}\label{tabl61}
	\begin{center}
		\begin{tabular}{|c|c|c|c|c|c|c|c|}
			\hline
Approach 1 & $L_{1}$~error&order&  $L_{\infty}$~error&order&iter&$\epsilon$\\ \hline
40 &1.91e-07 &- &2.05e-06 &- &46&$10^{-2}$\\ \hline
80 &3.96e-09 &5.59 &1.59e-08 &7.01 &45&$10^{-3}$\\ \hline
160 &1.27e-10 &4.96 &5.33e-10 &4.90 &56&$10^{-4}$\\ \hline
320 &4.19e-12 &4.92 &1.73e-11 &4.94 &89&$10^{-5}$\\ \hline
Approach 2 & $L_{1}$~error&order&  $L_{\infty}$~error&order&iter&$\epsilon$\\\hline
40 &2.58e-07 &- &3.25e-06 &- &41&$10^{-2}$\\ \hline
80 &4.97e-09 &5.69 &2.46e-08 &7.04 &45&$10^{-3}$\\ \hline
160 &1.60e-10 &4.95 &5.23e-10 &5.56 &54&$10^{-4}$\\ \hline
320 &5.38e-12 &4.89 &1.51e-11 &5.11 &91&$10^{-5}$\\ \hline
WENO-JP FSM & $L_{1}$~error&order&  $L_{\infty}$~error&order&iter&$\epsilon$\\ \hline
40 &4.61e-07 &- &2.10e-06 &- &36&$10^{-2}$\\ \hline
80 &4.83e-08 &3.25 &1.93e-07 &3.44 &37&$10^{-3}$\\ \hline
160 &1.97e-09 &4.61 &7.18e-09 &4.75 &52&$10^{-4}$\\ \hline
320 &4.12e-11 &5.58 &1.50e-10 &5.58 &71&$10^{-4}$\\ \hline
		\end{tabular}
	\end{center}
\end{table}

\begin{table}[h!]
\caption{Example 6 case b. Comparison of the three methods: The errors of the numerical solution, the accuracy obtained and the number of iterations for convergence}\label{tabl62}
	\begin{center}
		\begin{tabular}{|c|c|c|c|c|c|c|c|}
			\hline
Approach 1 & $L_{1}$~error&order&  $L_{\infty}$~error&order&iter&$\epsilon$\\\hline
40 &4.71e-04 &- &4.05e-03 &- &34&$10^{-2}$\\ \hline
80 &1.43e-04 &1.71 &1.23e-03 &1.71 &42&$10^{-3}$\\ \hline
160 &3.16e-05 &2.18 &3.02e-04 &2.02 &56&$10^{-4}$\\ \hline
320 &6.75e-06 &2.22 &7.30e-05 &2.04 &92&$10^{-5}$\\ \hline
Approach 2 & $L_{1}$~error&order&  $L_{\infty}$~error&order&iter&$\epsilon$\\\hline
40 &6.48e-04 &- &4.43e-03 &- &38&$10^{-2}$\\ \hline
80 &1.63e-04 &1.98 &1.31e-03 &1.75 &43&$10^{-3}$\\ \hline
160 &3.42e-05 &2.25 &3.22e-04 &2.02 &60&$10^{-4}$\\ \hline
320 &8.22e-06 &2.06 &6.15e-05 &2.38 &100&$10^{-5}$\\ \hline
WENO-JP FSM & $L_{1}$~error&order&  $L_{\infty}$~error&order&iter&$\epsilon$\\ \hline
40 &1.75e-04 &- &2.11e-03 &- &34&$10^{-2}$\\ \hline
80 &1.14e-04 &0.61 &1.12e-03 &0.91 &38&$10^{-3}$\\ \hline
160 &3.78e-05 &1.60 &3.18e-04 &1.81 &49&$10^{-4}$\\ \hline
320 &1.05e-05 &1.83 &8.79e-05 &1.85 &68&$10^{-5}$\\ \hline

		\end{tabular}
	\end{center}
\end{table}
\bigskip
\noindent \textbf{Example 7} The travel-time problem in elastic wave propagation is considered in this example. The
quasi-P and the quasi-SV slowness surfaces are defined as follows \cite{qianchengosher}
\begin{equation*}
  c_{1}\phi_{x}^{4}+c_{2}\phi_{x}^{2}\phi_{y}^{2}+c_{3}\phi_{y}^{4}+c_{4}\phi_{x}^{2}+c_{5}\phi_{y}^{2}+1=0,
\end{equation*}
where
\begin{equation*}
\begin{split}
  c_{1}=a_{11}a_{44}, \quad c_{2}=a_{11}a_{33}+a_{44}^{2}-(a_{13}+a_{44})^{2},&  \\
    c_{3}=a_{33}a_{44},\quad c_{4}=-(a_{11}+a_{44}),\quad c_{5}=-(a_{33}+a_{44}),&
\end{split}
\end{equation*}
in which $a_{i,j}$ are given elastic parameters. The quasi-P wave Eikonal equation is
\begin{equation*}
  \sqrt{-\frac{1}{2}(c_{4}\phi_{x}^{2}+c_{5}\phi_{y}^{2})+\sqrt{\frac{1}{4}(c_{4}\phi_{x}^{2}+c_{5}\phi_{y}^{2})^2-(c_{1}\phi_{x}^{4}+c_{2}\phi_{x}^{2}\phi_{y}^{2}+c_{3}\phi_{y}^{4})}}=1,
\end{equation*}
which is a convex HJ equation, and the elastic parameters are taken to be
\begin{equation*}
  a_{11}=15.0638,~ a_{33}=10.8373,~a_{13}=1.6381,~a_{44}=3.1258.
\end{equation*}
The corresponding quasi-SV wave Eikonal equation is given by
 \begin{equation*}
  \sqrt{-\frac{1}{2}(c_{4}\phi_{x}^{2}+c_{5}\phi_{y}^{2})-\sqrt{\frac{1}{4}(c_{4}\phi_{x}^{2}+c_{5}\phi_{y}^{2})^2-(c_{1}\phi_{x}^{4}+c_{2}\phi_{x}^{2}\phi_{y}^{2}+c_{3}\phi_{y}^{4})}}=1,
\end{equation*}
which is a nonconvex HJ equation, and the elastic parameters are taken to be
\begin{equation*}
  a_{11}=15.90,~ a_{33}=6.21,~a_{13}=4.82,~a_{44}=4.00.
\end{equation*}
The computational domain is set as $\Omega=[-1,1]^2$, and the inflow boundary is given by $\Gamma={(0,0)}$. Exact values are assigned in a small box with length $0.3$ around the source point. Because these Hamiltonians are pretty complicated, we use the Lax-Friedrich numerical Hamiltonian for both equations. In addition, since we only know the numerical solution of $\phi$, the ``exact solution" of $u$ and $v$ on \emph{Category I} will be obtained by fifth order WENO-JP reconstruction.

For quasi-P wave equation, the numerical errors and orders of convergence are presented in Table \ref{tabl9p} for three methods. We can observe that all three methods achieve the fifth order accuracy as we expected. Moreover, the numbers of iterations required by three methods are basically the same, while Approach 1 and Approach 2 produce smaller numerical errors. The CPU time are plotted in Figure \ref{fig1h8}, which shows that WENO-JP costs the least CPU time, but it has the largest error. And, Approach 1 and Approach 2 have the same error, but Approach 1 costs less CPU time than Approach 2.

\begin{figure}[!h]
\begin{center}
	\includegraphics[width=7.6cm]{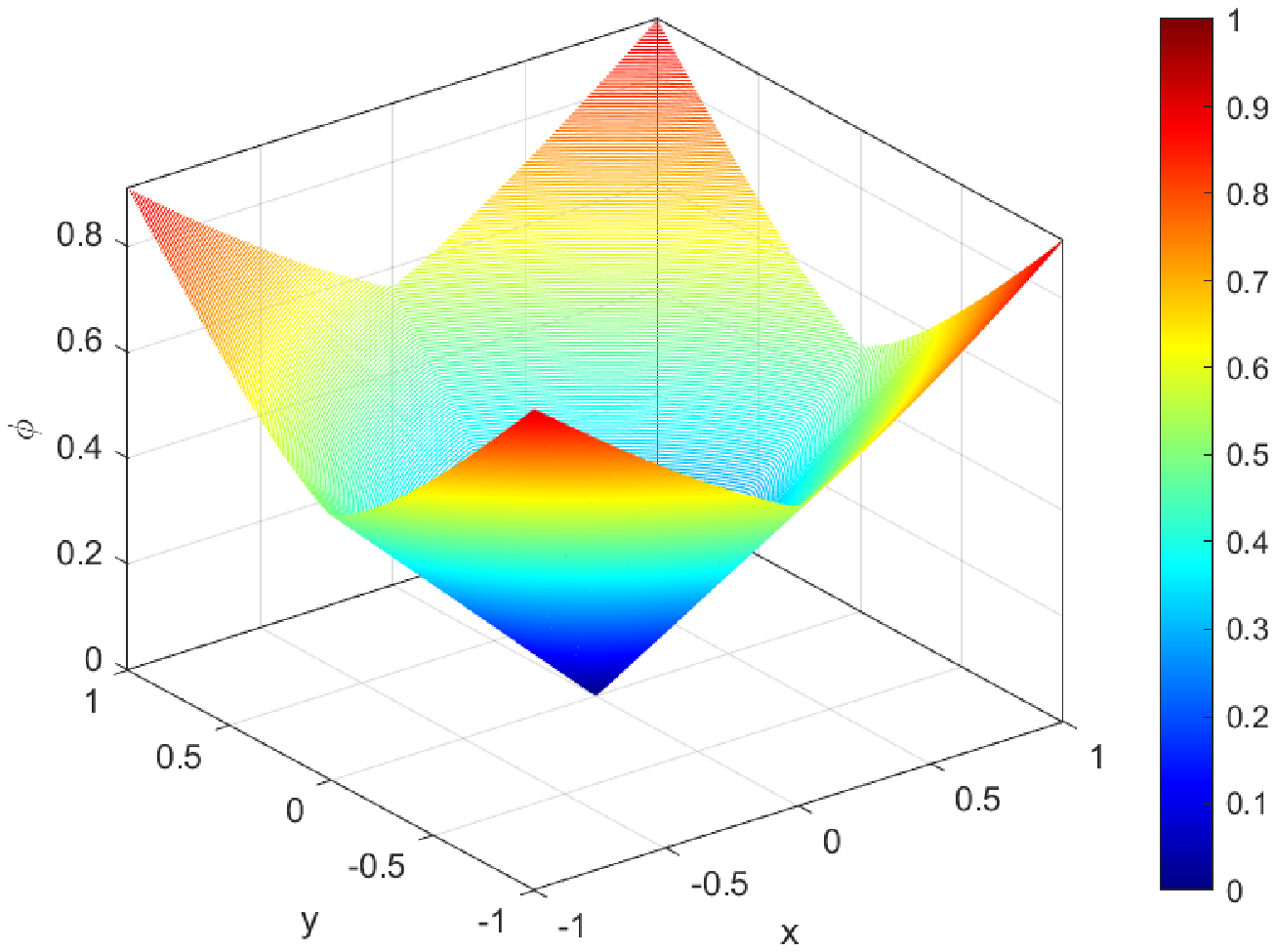}
	\includegraphics[width=7.6cm]{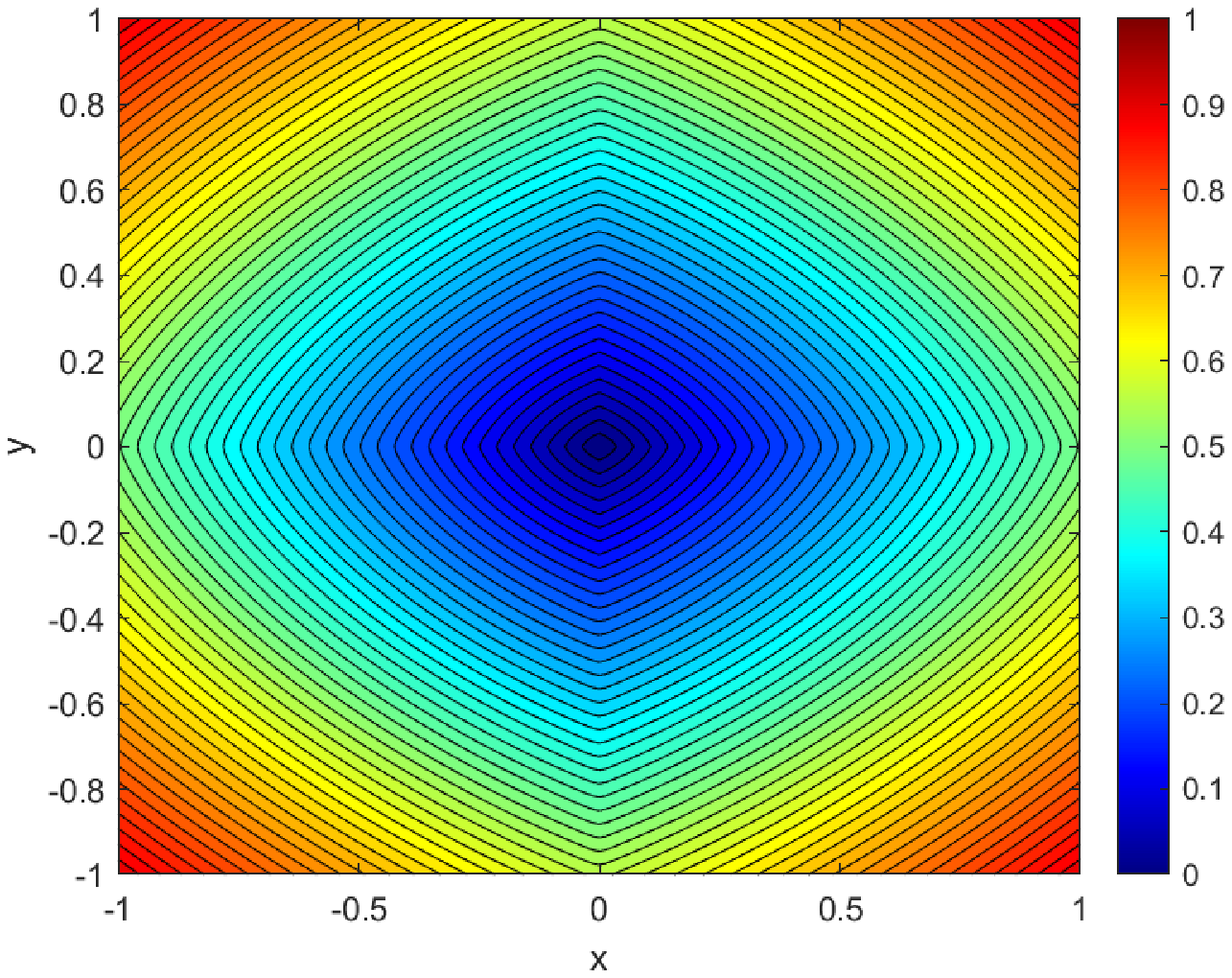}
\end{center}
\caption{Example 7 SV-wave. The numerical solution by Approach 1 on mesh $N=160$. Left: the 3D plot of numerical solution $\phi$; Right: the contour plot for $\phi$ .}\label{e72solution}
\end{figure}

For the quasi-SV wave equation, we set the $\delta<10^{-9}$ for Approach 2 on mesh $N=80,160,320$, otherwise the Approach 2 will not convergence.  The errors are measured in the region away from the singular lines of $x=0$ and $y=0$. The Figure \ref{e72solution} shows that picture of numerical solution for SV-wave. The numerical errors and orders of convergence are shown in Table \ref{tabl9sv}, from which we observe that the designed fifth order accuracy is again achieved. Approach 1 and Approach 2 yield smaller numerical errors for this test case, and in the same time, Approach 2 requires slightly less numbers of iterations when compared with WENO-JP method. 

\begin{table}[!h]
\caption{Example 7 P-wave. Comparison of the three methods: The errors of the numerical solution, the accuracy obtained and the number of iterations for convergence}\label{tabl9p}
	\begin{center}
		\begin{tabular}{|c|c|c|c|c|c|c|c|}
			\hline
Approach 1 & $L_{1}$~error&order&  $L_{\infty}$~error&order&iter\\ \hline
40 &4.78e-06 &- &3.95e-05 &- &41\\ \hline
80 &2.07e-07 &4.52 &2.25e-06 &4.13 &44\\ \hline
160 &7.05e-09 &4.87 &8.09e-08 &4.80 &56\\ \hline
320 &2.29e-10 &4.93 &2.61e-09 &4.95 &77\\ \hline
Approach 2 & $L_{1}$~error&order&  $L_{\infty}$~error&order&iter\\ \hline
40 &5.77e-06 &- &4.35e-05 &- &34\\ \hline
80 &2.27e-07 &4.66 &2.30e-06 &4.23 &42\\ \hline
160 &7.42e-09 &4.93 &7.44e-08 &4.95 &54\\ \hline
320 &2.40e-10 &4.94 &2.24e-09 &5.05 &78\\ \hline
WENO-JP FSM & $L_{1}$~error&order&  $L_{\infty}$~error&order&iter\\ \hline
40 &2.56e-05 &- &2.67e-04 &- &37\\ \hline
80 &1.64e-06 &3.96 &2.05e-05 &3.70 &44\\ \hline
160 &5.46e-08 &4.90 &6.85e-07 &4.90 &56\\ \hline
320 &1.51e-09 &5.17 &1.44e-08 &5.57 &83\\ \hline
		\end{tabular}
	\end{center}
\end{table}

\begin{table}[!h]
\caption{Example 7 SV-wave. Comparison of the three methods: The errors of the numerical solution, the accuracy obtained and the number of iterations for convergence}\label{tabl9sv}
	\begin{center}
		\begin{tabular}{|c|c|c|c|c|c|c|c|}
			\hline
Approach 1 & $L_{1}$~error&order&  $L_{\infty}$~error&order&iter \\\hline
80 &9.63e-07 &- &1.60e-05 &- &71\\ \hline
160 &1.95e-08 &5.61 &8.97e-07 &4.15 &92\\ \hline
320 &6.37e-11 &8.26 &9.12e-09 &6.61 &112\\ \hline
640 &5.60e-13 &6.82 &1.65e-11 &9.10 &190\\ \hline
Approach 2 & $L_{1}$~error&order&  $L_{\infty}$~error&order&iter\\\hline
80 &1.76e-06 &- &2.39e-05 &- &43\\ \hline
160 &5.03e-09 &7.64 &3.91e-07 &5.93 &62\\ \hline
320 &2.85e-11 &7.46&1.03e-09 &8.56 &95\\ \hline
640 &7.88e-13 &5.17 &5.95e-11 &4.11 &170\\ \hline
WENO-JP FSM & $L_{1}$~error&order&  $L_{\infty}$~error&order&iter\\ \hline
80 &1.28e-06 &- &2.02e-05 &- &50\\ \hline
160 &1.99e-08 &6.01 &8.26e-07 &4.61 &70\\ \hline
320 &1.80e-10 &6.79 &1.15e-08 &6.16 &108\\ \hline
640 &4.11e-12 &5.45 &1.16e-10 &6.62 &181\\ \hline
		\end{tabular}
	\end{center}
\end{table}

\section{Hybrid strategy}\label{section6}
\setcounter{equation}{0}
\setcounter{figure}{0}
\setcounter{table}{0}

In the HWENO reconstruction procedure, the evaluation of the smoothness indicators occupies most of the extra computational costs when compared with the simple linear reconstruction. In this section, we explore a hybrid method which combines both linear and HWENO method, similar to the hybrid WENO fast sweeping method studied in \cite{me}. The main idea is to adopt the fifth order linear reconstruction on the big stencil $S_0$ or $\tilde{S}_0$ if the numerical solution is monotone, that is $u_{i,j}$ do not change sign in either $S_0$ or $\tilde{S}_0$. Otherwise, the HWENO reconstruction is used.

In Section \ref{sec4.3}, the points $\{(x_i, y_j)\}$ were classified into two categories. Here, we further separate the \emph{Category IV} into the following two subcategories, which will be handled slightly differently in the hybrid method. \\
\noindent\emph{Category IV.1}: For points whose distances to \emph{Category III} are less than or equal to $2h$ (excluding \emph{Category I}). These points will be updated by FSM. \\
\noindent\emph{Category IV.2}: All remaining points, which will also be updated by FSM. \\
Below, we present the flowchart for \textbf{Approach 1}, coupled with the hybrid strategy.

\noindent \textbf{Step 1}. \emph{Initialization}:
The numerical solution from the first order fast sweeping method \cite{Zhaofsm1} is taken as the initial guess of $\phi$. The forward or backward difference of this $\phi$ is used as the initial guess of $u$ and $v$.

\noindent \textbf{Step 2}. \emph{Gauss-Seidel iteration}. We Solve the discretized nonlinear system by GS iterations with four alternating direction sweepings:
$$(1)\quad i=1:N_x, ~j=1:N_y; \qquad (2)\quad i=N_x:1,~ j=1:N_y;$$
$$(3)\quad i=N_x:1,~ j=N_y:1; \qquad (4)\quad i=1:N_x,~ j=N_y:1.$$
During each sweeping, the updating strategy for the points in \emph{Category IV} is outlined below.
For the points in \emph{Category IV.1}:  the HWENO reconstruction (\ref{weno}) is applied to evaluate $(\phi_{x})_{i,j}^{\pm}$, and the similar procedure for $(\phi_{y})_{i,j}^{\pm}$ is used. For the points in \emph{Category IV.2}, the hybrid strategy is applied here, and we define
\begin{equation}\label{signf}
(\phi_{x})_{i,j}^{\pm}=
\begin{cases}
\eqref{hybridlinearf} ~\mathrm{or}~ \eqref{hybridlinearz}, &\mathrm{if~} \{u_{i,j}\}\mathrm{~have~the~same~sign~on~}S_{0}~or~ \widetilde{S_{0}},\\
\eqref{weno}; &      \mathrm{otherwise}.
\end{cases}
\end{equation}
Similarly, one can evaluate $(\phi_{y})_{i,j}^{\pm}$ along the $y$-direction. The rest of the algorithm is the same as before. $\phi_{i,j}^{new}$ is updated using \eqref{lf} or \eqref{phinew}, and $u_{i,j}^{new}$ and $v_{i,j}^{new}$  are evaluated by \eqref{app2} in each sweeping direction. High order extrapolations are used at the ghost points.

\noindent \textbf{Step 3}. \emph{Convergence}:  In general, the iteration will stop if, for two consecutive iteration steps,
$$\delta=||\phi^{new}-\phi^{old}||_{L_{1}}<10^{-14}.$$
The procedure of \textbf{Approach 2} with hybrid strategy is similar, and is omitted here.

\begin{figure}[htbp!]
\begin{center}
	\includegraphics[width=5.4cm]{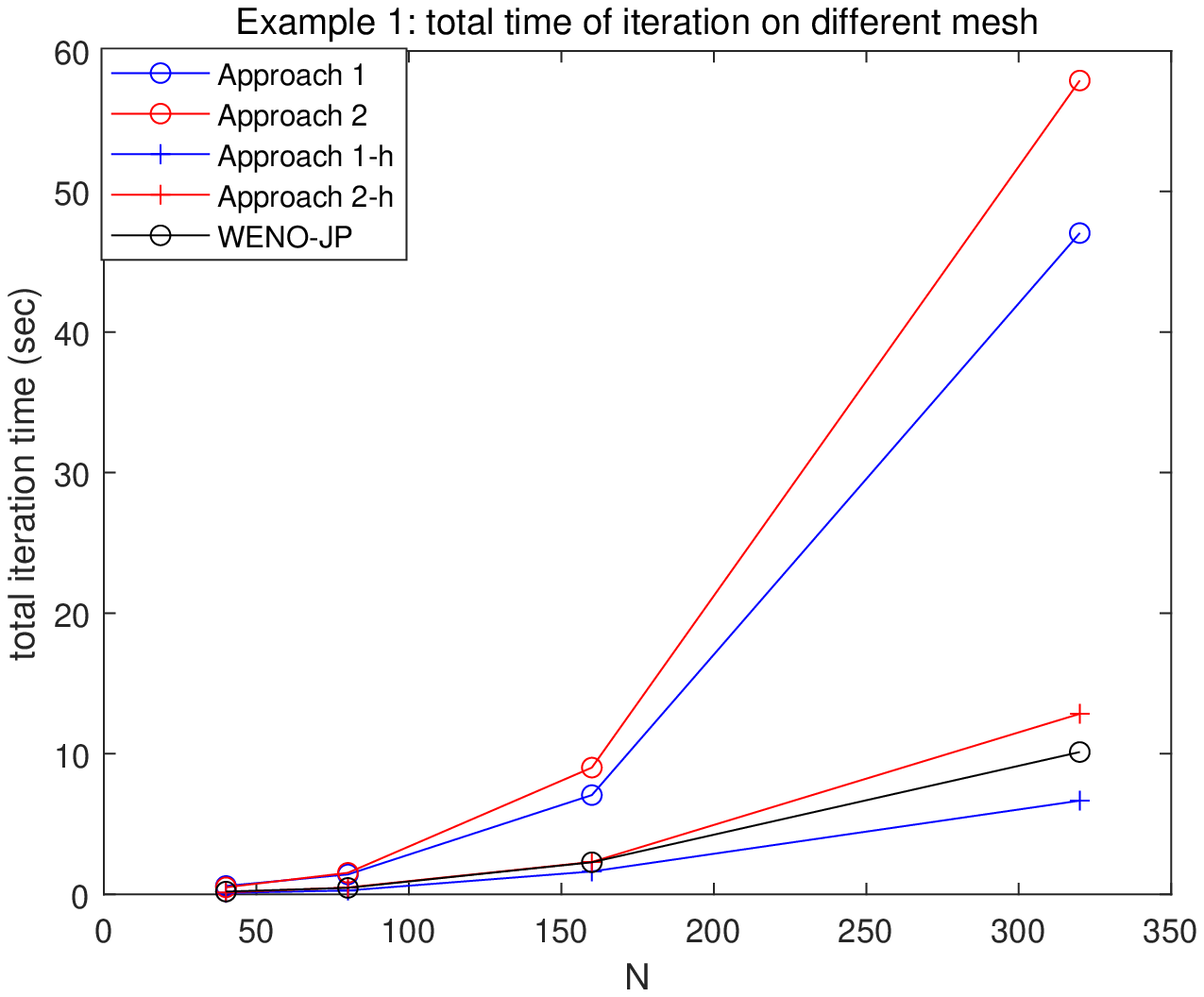}
	\includegraphics[width=5.4cm]{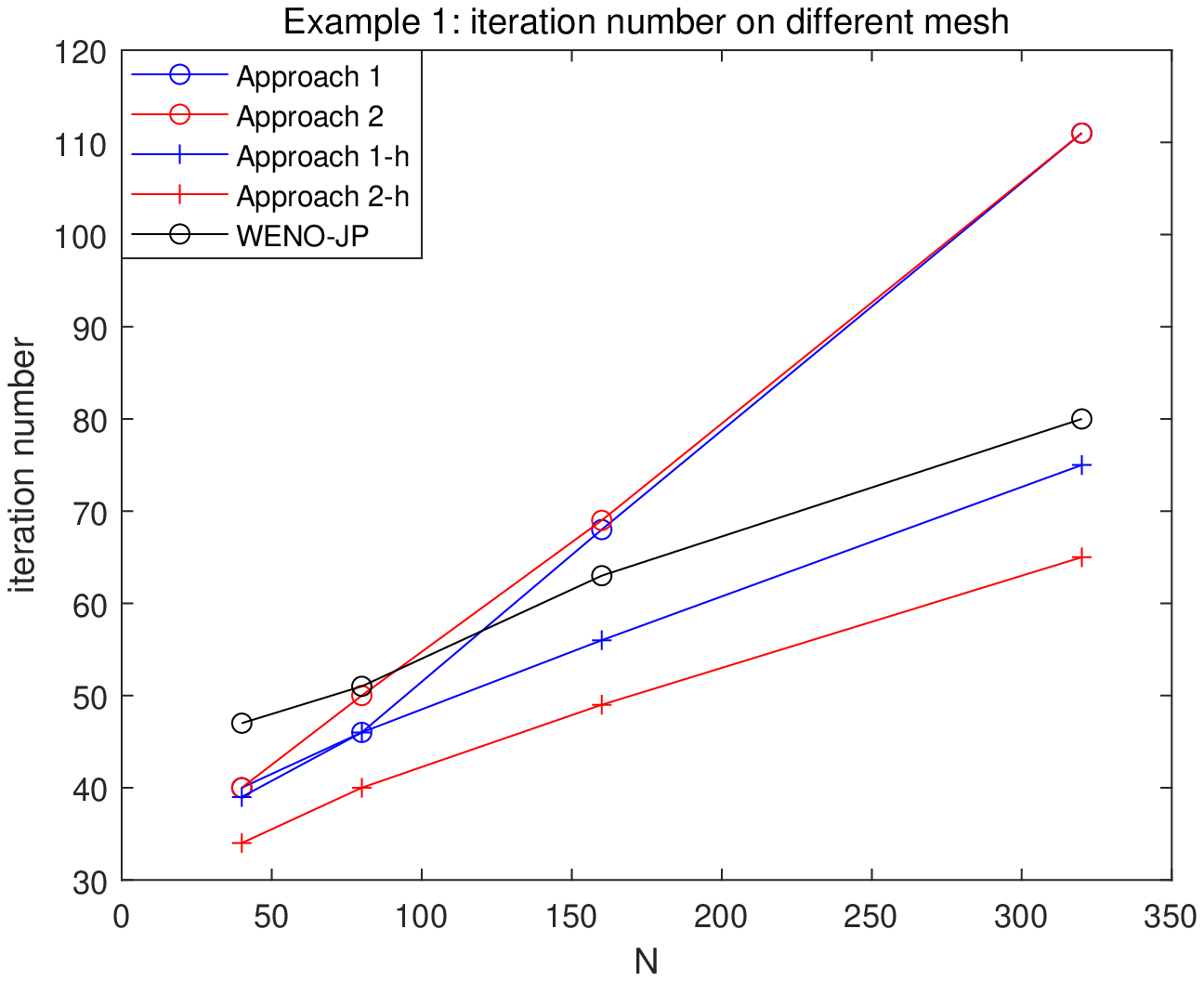}
	\includegraphics[width=5.4cm]{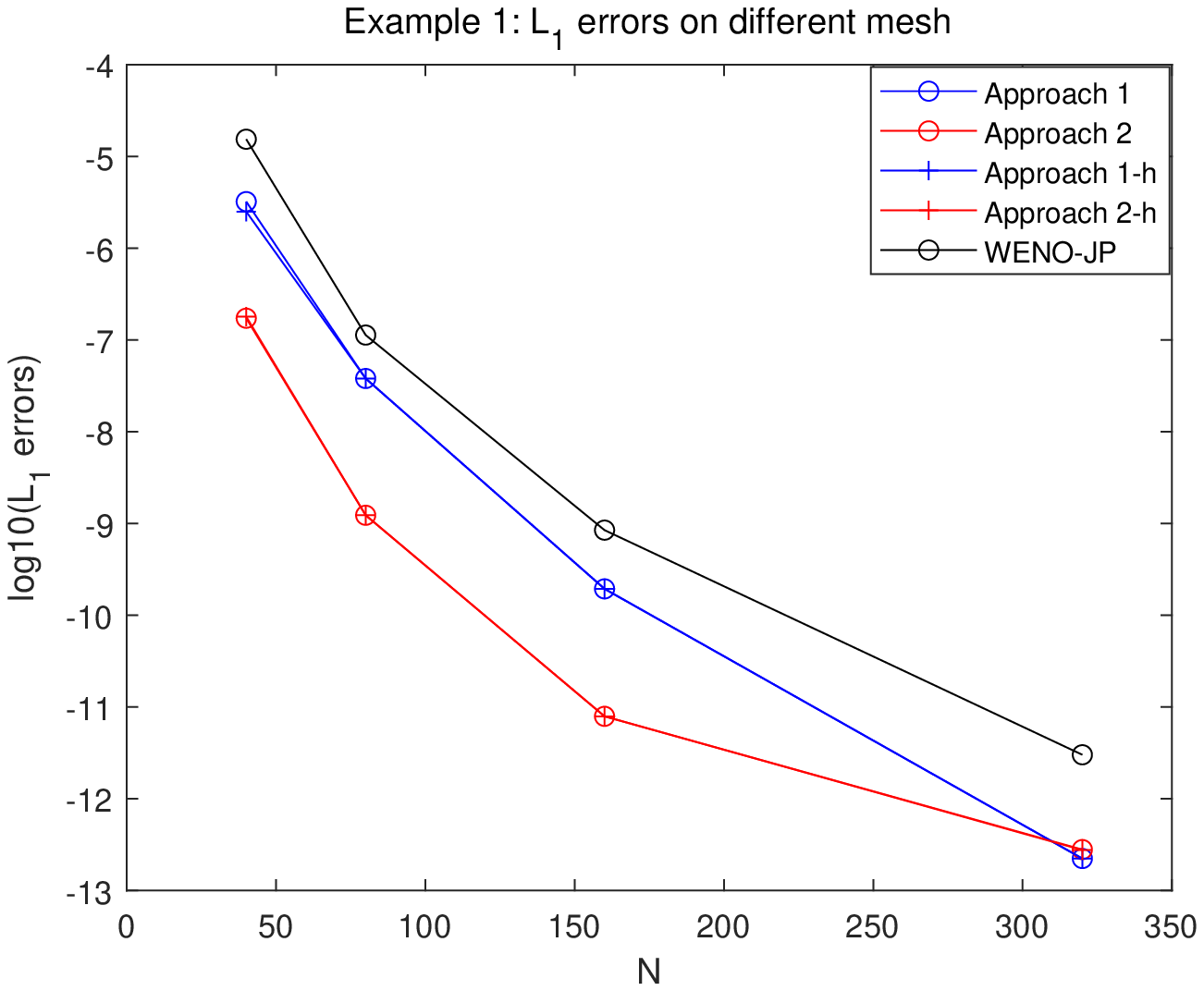}
\end{center}
\caption{Example 1 with Approach 1, 2, the hybrid Approach 1, 2 and WENO-JP. Left: mesh number $N$ vs CPU time; Middle: mesh number $N$ vs number of iterations; Right: mesh number N vs $L_{1}$ error.
}\label{fig1h}
\end{figure}

We expect these hybrid HWENO FSM to be more efficient. Next, some numerical results of the hybrid HWENO FSM will be presented.
All the seven examples in Section \ref{section5} have been tested using the hybrid algorithms, and the comparison of their numerical performance with those of two HWENO methods in Section \ref{section4} are presented in Table \ref{tab1h}. In Figures \ref{fig1h}-\ref{fig1h9}, the comparison of their CPU time, number of iterations and $L_1$ numerical errors on various mesh size is provided. In the Table, we denote the $L_{1}$ error and $L_{\infty}$ error in Approach 1 with hybrid strategy by ``A1-$L_1$" and ``A1-$L_{\infty}$", respectively. Similar notations are adopted for Approach 2 with hybrid strategy.
In these figures, we denote Approach 1 with hybrid strategy by ``Approach 1-h", and Approach 2 with hybrid strategy by ``Approach 2-h". The ``blue -o-" lines indicates Approach 1, the ``blue -+-" lines indicates hybrid Approach 1, the ``red -o-" lines indicates Approach 2, the ``red -+-" lines indicates hybrid Approach 2, and the ``black -o-" lines indicates WENO-JP.
The table and all of these figures demonstrate that the hybrid schemes cost much less CPU time, and converge with a smaller number of iterations when the mesh is refined. In addition, for most of these examples, Approach 1 enjoys more savings in computational time than Approach 2. The computational time of the hybrid HWENO FSM is comparable to that of the WENO-JP scheme, yet the numerical error of WENO-JP scheme is the largest among all these five methods, when the same mesh size is considered.

\begin{table}[htbp!]
\caption{All examples. Comparison of the two hybrid methods: The errors of the numerical solution, the accuracy obtained and the number of iterations for convergence}\label{tab1h}
	\begin{center}
		\begin{tabular}{cc|ccccc|ccccc}
			\hline
Test&N & $A1-L_{1}$&order&  $A1-L_{\infty}$&order&iter & $A2-L_{1}$&order&  $A2-L_{\infty}$&order&iter\\\hline
1&40 &2.51e-06 &- &1.94e-05 &- &39&1.81e-07 &- &1.35e-06 &- &33     \\
&80 &3.80e-08 &6.04 &7.13e-07 &4.76 &46&1.33e-09 &7.08 &5.25e-08 &4.68 &41     \\
&160 &1.94e-10 &7.60 &7.69e-09 &6.53 &58&9.82e-12 &7.08 &1.62e-11 &11.65 &51      \\
&320 &2.49e-13 &9.60 &8.37e-12 &9.84 &76&3.04e-13 &5.01 &2.17e-12 &2.90 &65       \\ \hline
2&40 &6.04e-07 &- &2.34e-05 &- &37&1.92e-07 &- &1.09e-05 &- &25              \\
&80 &1.17e-08 &5.68 &1.74e-06 &3.74 &44&1.65e-09 &6.85 &1.87e-07 &5.85 &31          \\
&160 &8.74e-11 &7.06 &2.21e-08 &6.30 &57&5.45e-11 &4.92 &3.67e-09 &5.67 &42          \\
&320 &1.78e-12 &5.61 &1.30e-10 &7.40 &83&2.01e-12 &4.75 &1.23e-10 &4.89 &61            \\ \hline
3&80 &7.14e-07 &- &1.04e-04 &- &59&1.57e-06 &- &8.57e-05 &- &35                 \\
&160 &1.94e-07 &1.87 &1.81e-05 &2.53 &55&3.19e-07 &2.30 &4.26e-05 &1.00 &47            \\
&320 &3.14e-10 &9.27 &3.86e-07 &5.55 &67&9.77e-10 &8.35 &6.55e-07 &6.02 &69            \\
&640 &6.41e-12 &5.61 &1.25e-09 &8.26&121&6.44e-12 &7.24 &7.85e-10 &9.70 &117            \\ \hline
4&40 &3.10e-07 &- &4.60e-06 &- &42&8.42e-08 &- &7.14e-07 &- &30                \\
&80 &6.95e-09 &5.48 &1.62e-07 &4.82 &49&2.88e-09 &4.87 &1.14e-08 &5.96 &36            \\
&160 &1.19e-10 &5.85 &1.78e-09 &6.51 &64&1.04e-10 &4.78 &2.55e-10 &5.47 &48               \\
&320 &3.26e-12 &5.20 &8.67e-12 &7.68 &88&3.57e-12 &4.87 &1.80e-11 &3.82 &73               \\ \hline
5&40 &1.12e-06 &- &1.61e-05 &- &84&8.84e-07 &- &1.52e-05 &- &32                     \\
&80 &4.28e-08 &4.71 &1.08e-06 &3.89 &60&4.11e-08 &4.42 &1.11e-06 &3.77 &42                \\
&160 &1.05e-09 &5.33 &3.93e-08 &4.78 &74&1.41e-09 &4.85 &5.95e-08 &4.22 &66                     \\
&320 &2.11e-11 &5.64 &4.44e-10 &6.46 &92&3.17e-11 &5.47 &1.55e-09 &5.26 &105                  \\ \hline
6a&40 &1.71e-07 &- &2.05e-06 &- &47&2.26e-07 &- &3.25e-06 &- &44             \\
&80 &3.45e-09 &5.63 &1.57e-08 &7.02 &46&4.28e-09 &5.72 &2.46e-08 &7.04 &46            \\
&160 &1.20e-10 &4.84 &5.25e-10 &4.90 &52&1.52e-10 &4.81 &5.23e-10 &5.56 &51                    \\
&320 &4.11e-12 &4.87 &1.72e-11 &4.92 &74&5.29e-12 &4.84 &1.39e-11 &5.23 &80                     \\ \hline
6b&40 &3.35e-04 &- &2.79e-03 &- &55&5.18e-04 &- &4.05e-03 &- &30                      \\
&80 &9.15e-05 &1.87 &6.44e-04 &2.11 &35&8.45e-05 &2.61 &4.86e-04 &3.05 &34                   \\
&160 &2.39e-05 &1.93 &1.67e-04 &1.94 &42&1.55e-05 &2.44 &1.28e-04 &1.91 &47                  \\
&320 &4.28e-06 &2.48 &4.41e-05 &1.92 &71&4.56e-06 &1.76 &3.91e-05 &1.71 &80                \\ \hline
7p&40 &4.78e-06 &- &3.95e-05 &- &46&5.77e-06 &- &4.35e-05 &- &33                   \\
&80 &2.07e-07 &4.52 &2.25e-06 &4.13 &48&2.27e-07 &4.66 &2.30e-06 &4.23 &41                 \\
&160 &7.05e-09 &4.87 &8.09e-08 &4.80 &56&7.42e-09 &4.93 &7.44e-08 &4.95 &51                 \\
&320 &2.29e-10 &4.93 &2.61e-09 &4.95 &82&2.40e-10 &4.94 &2.24e-09 &5.05 &80                  \\ \hline
7sv&80 &8.61e-07 &- &2.19e-05 &- &68&8.58e-07 &- &2.09e-05 &- &35                         \\
&160 &2.25e-08 &5.25 &1.44e-06 &3.92 &78&6.08e-09 &7.14 &5.04e-07 &5.37 &54                      \\
&320 &7.85e-11 &8.16 &7.90e-09 &7.51 &113&2.90e-11 &7.71 &1.45e-09 &8.44 &126                      \\
&640 &5.68e-13 &7.11 &1.84e-11 &8.74 &189&7.85e-13 &5.20 &2.57e-11 &5.82 &162                          \\ \hline
		\end{tabular}
	\end{center}
\end{table}
\begin{figure}[htbp!]
\begin{center}
	\includegraphics[width=7.6cm]{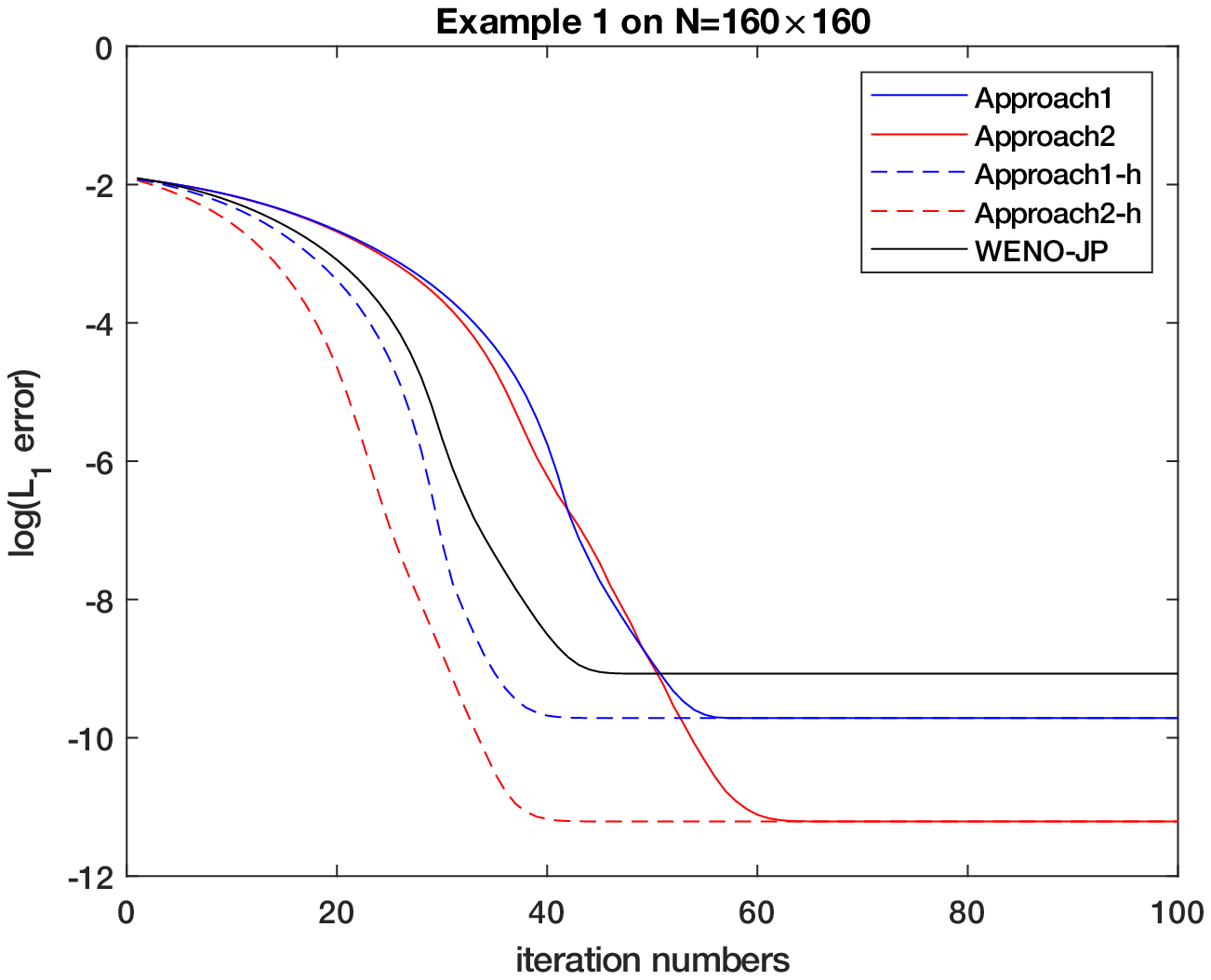}
	\includegraphics[width=7.6cm]{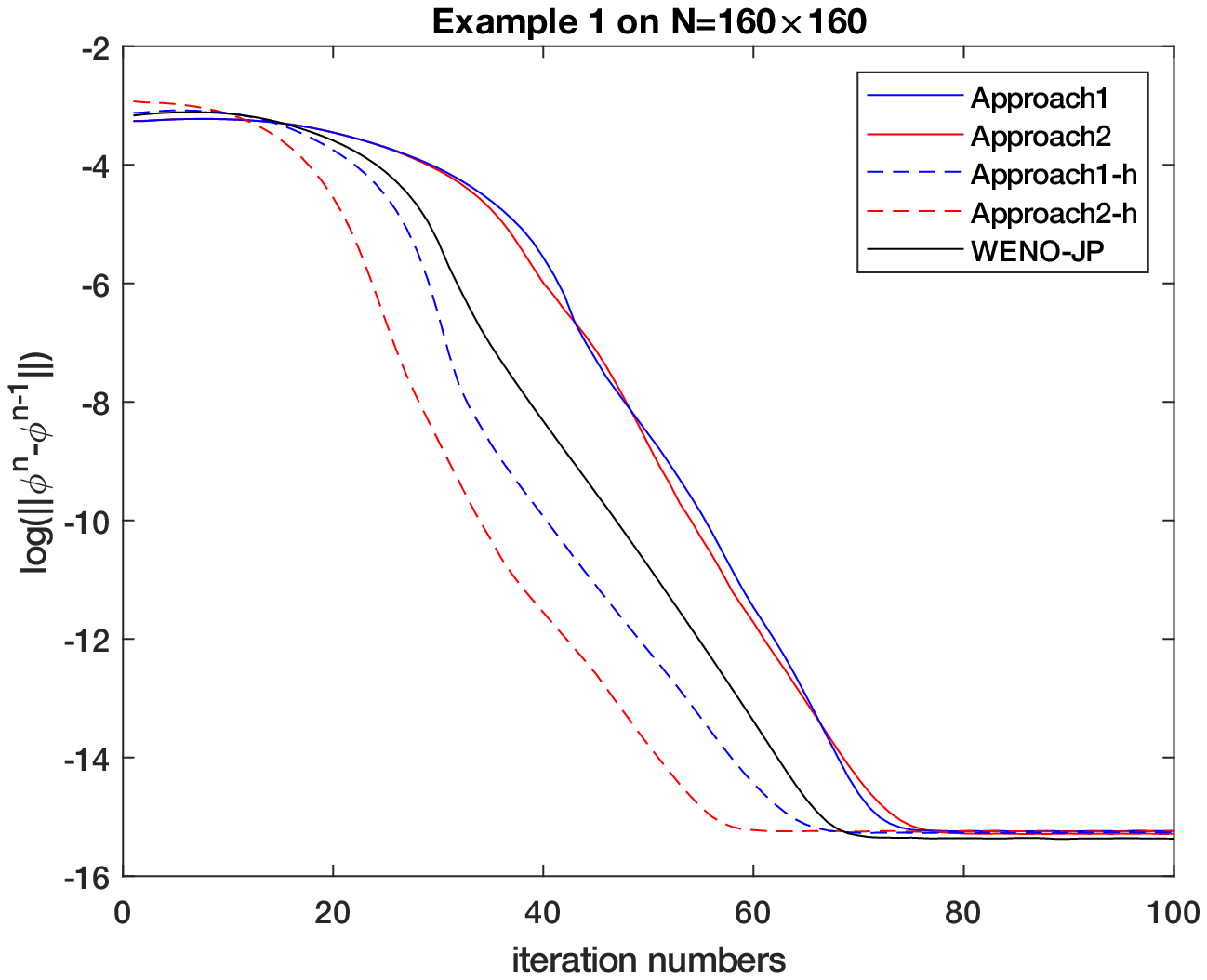}
    \end{center}
\caption{Fifth methods solving Example 1. Left: number of iterations vs $log(L_{1}~error$); Right: number of iterations vs $log(||\phi^{n}-\phi^{n-1}||_{L_{1}}$).}
\end{figure}
\begin{figure}[htbp!]
\begin{center}
	\includegraphics[width=5.4cm]{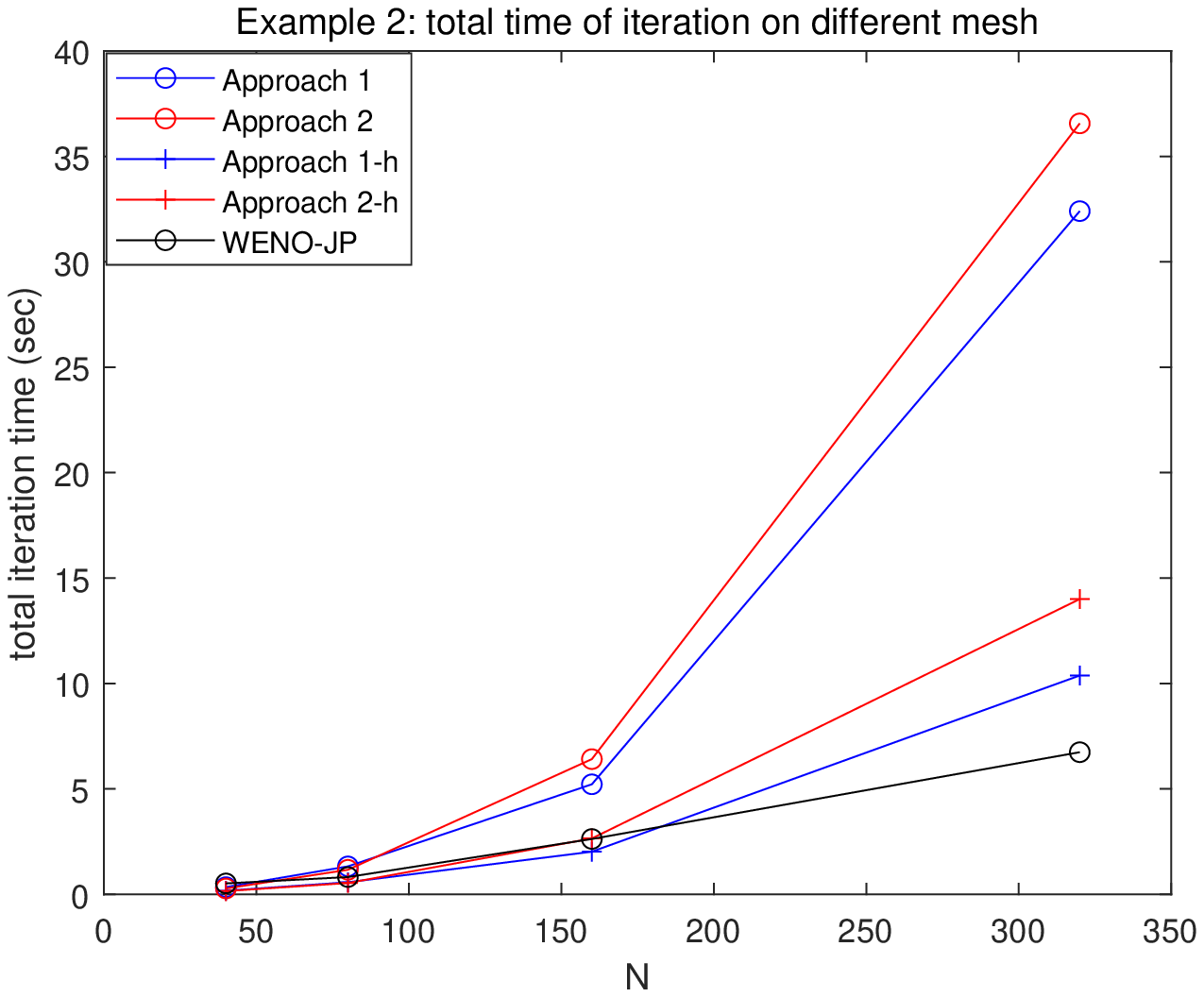}
	\includegraphics[width=5.4cm]{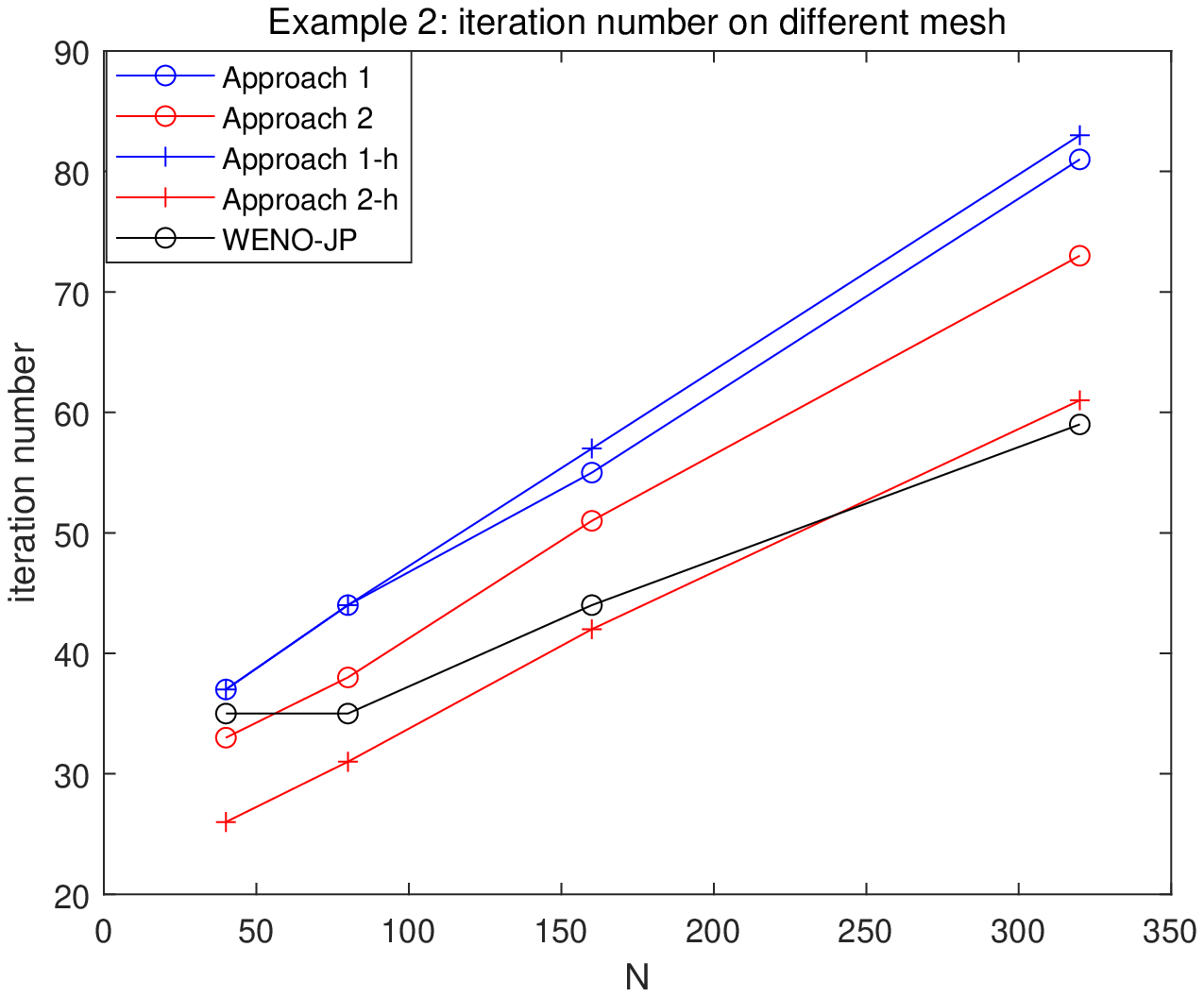}
	\includegraphics[width=5.4cm]{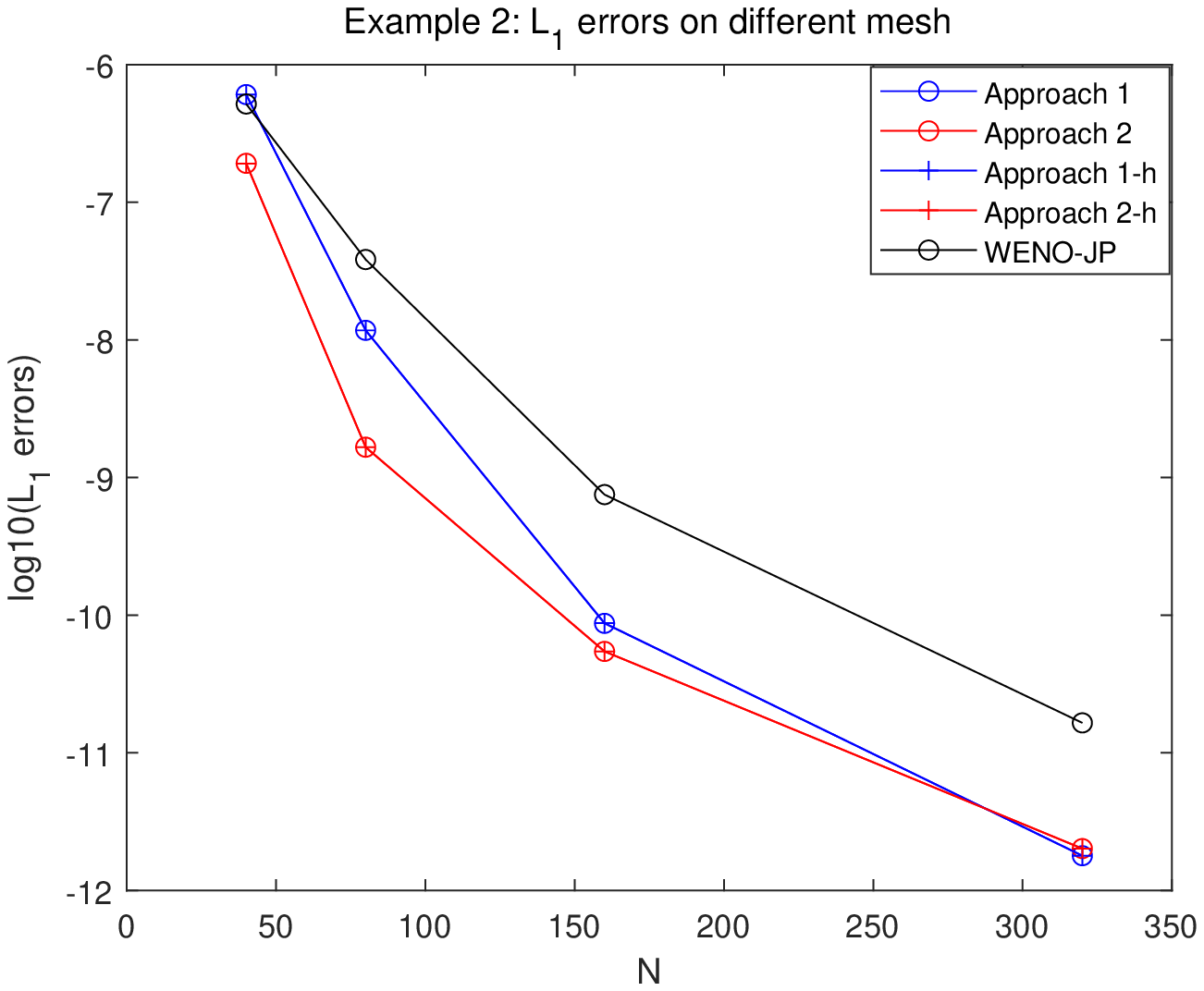}
\end{center}
\caption{Example 2 with Approach 1, 2, the hybrid Approach 1, 2 and WENO-JP. Left: mesh number $N$ vs CPU time; Middle: mesh number $N$ vs number of iterations; Right: mesh number N vs $L_{1}$ error.}\label{fig1h2}
\end{figure}
\begin{figure}[htbp!]
\begin{center}
	\includegraphics[width=7.6cm]{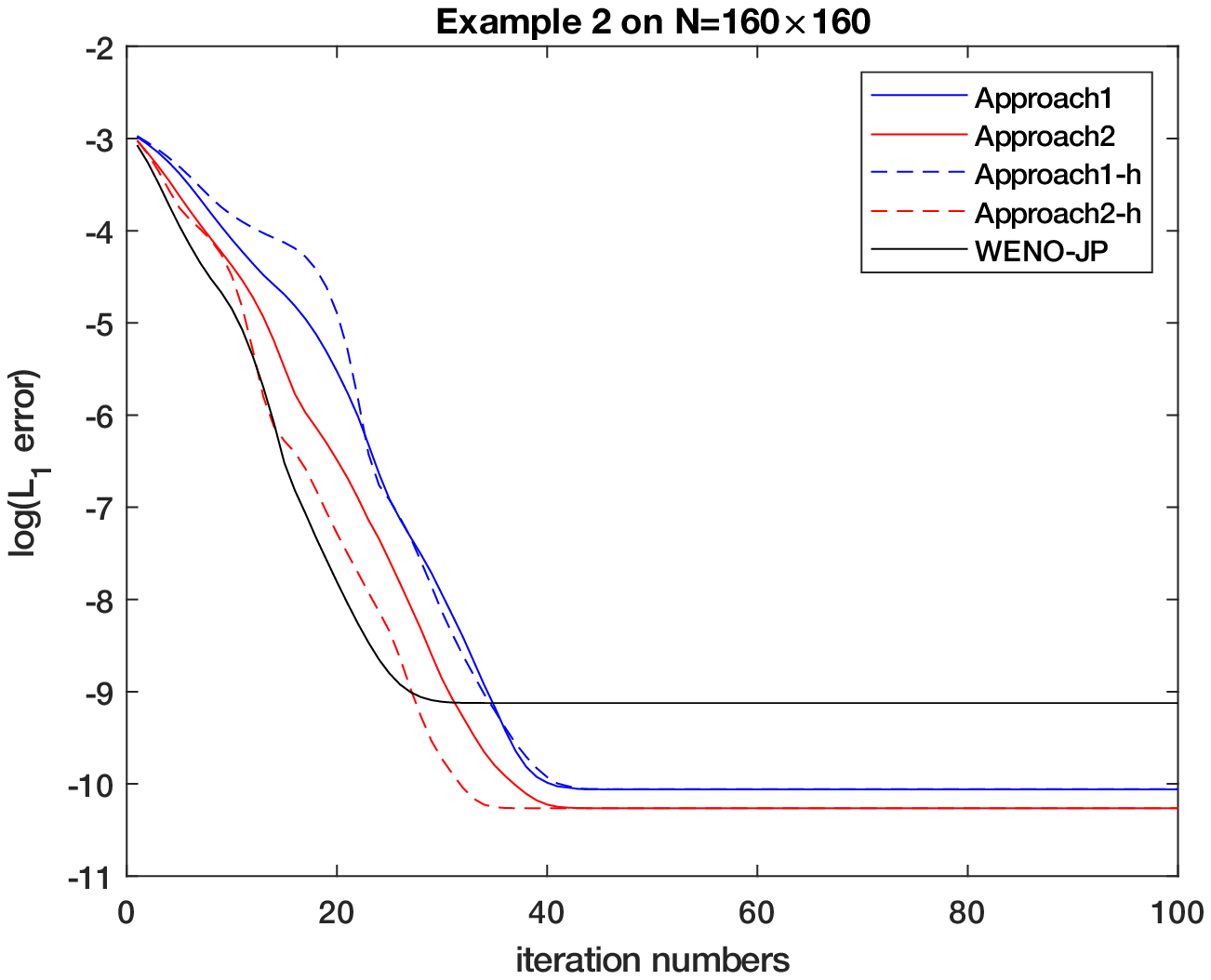}
	\includegraphics[width=7.6cm]{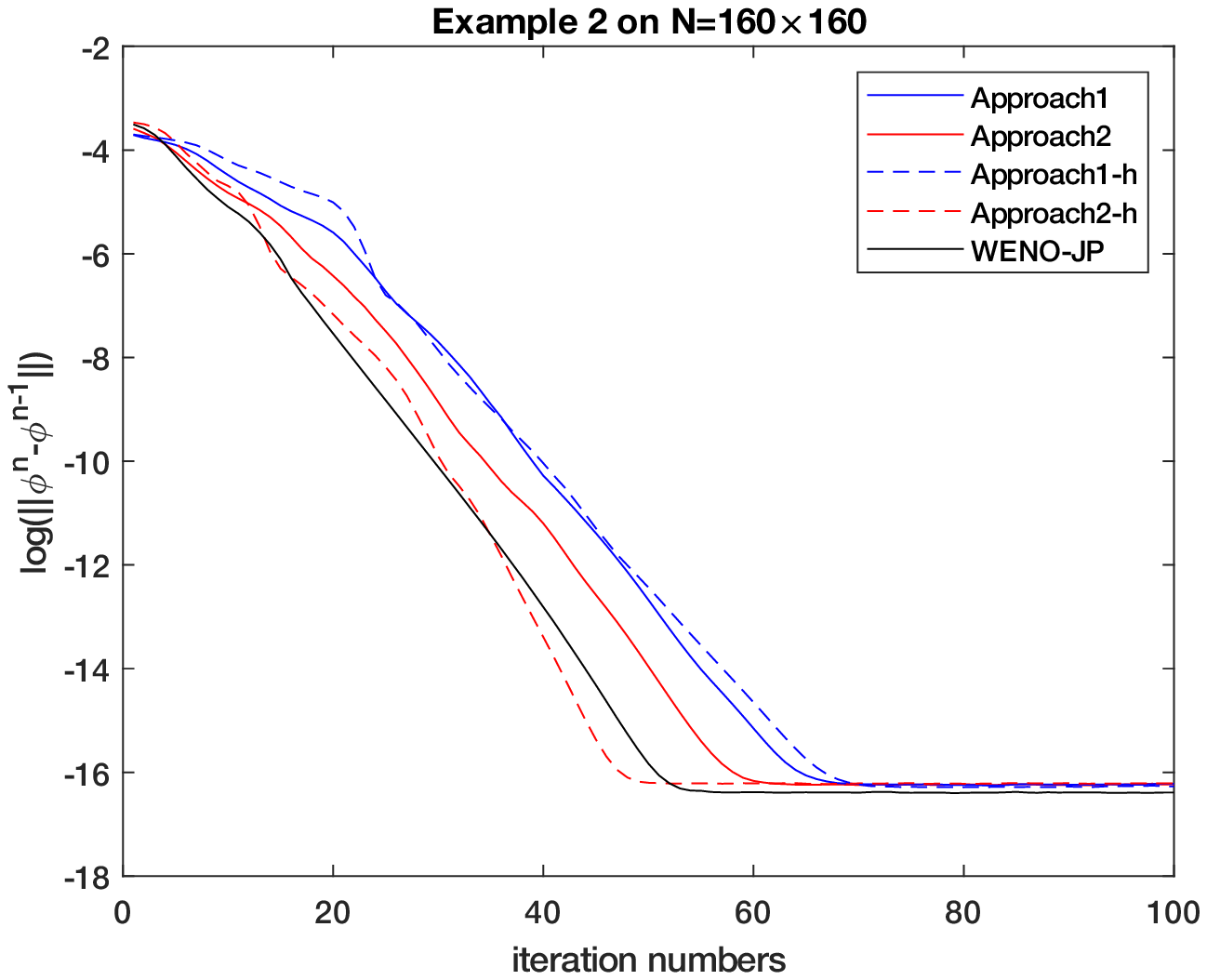}
    \end{center}
\caption{Three methods solving Example 2. Left: number of iterations vs $log(L_{1}~error$); Right: number of iterations vs $log(||\phi^{n}-\phi^{n-1}||_{L_{1}}$).}
\end{figure}

\begin{figure}[htbp!]
\begin{center}
	\includegraphics[width=5.4cm]{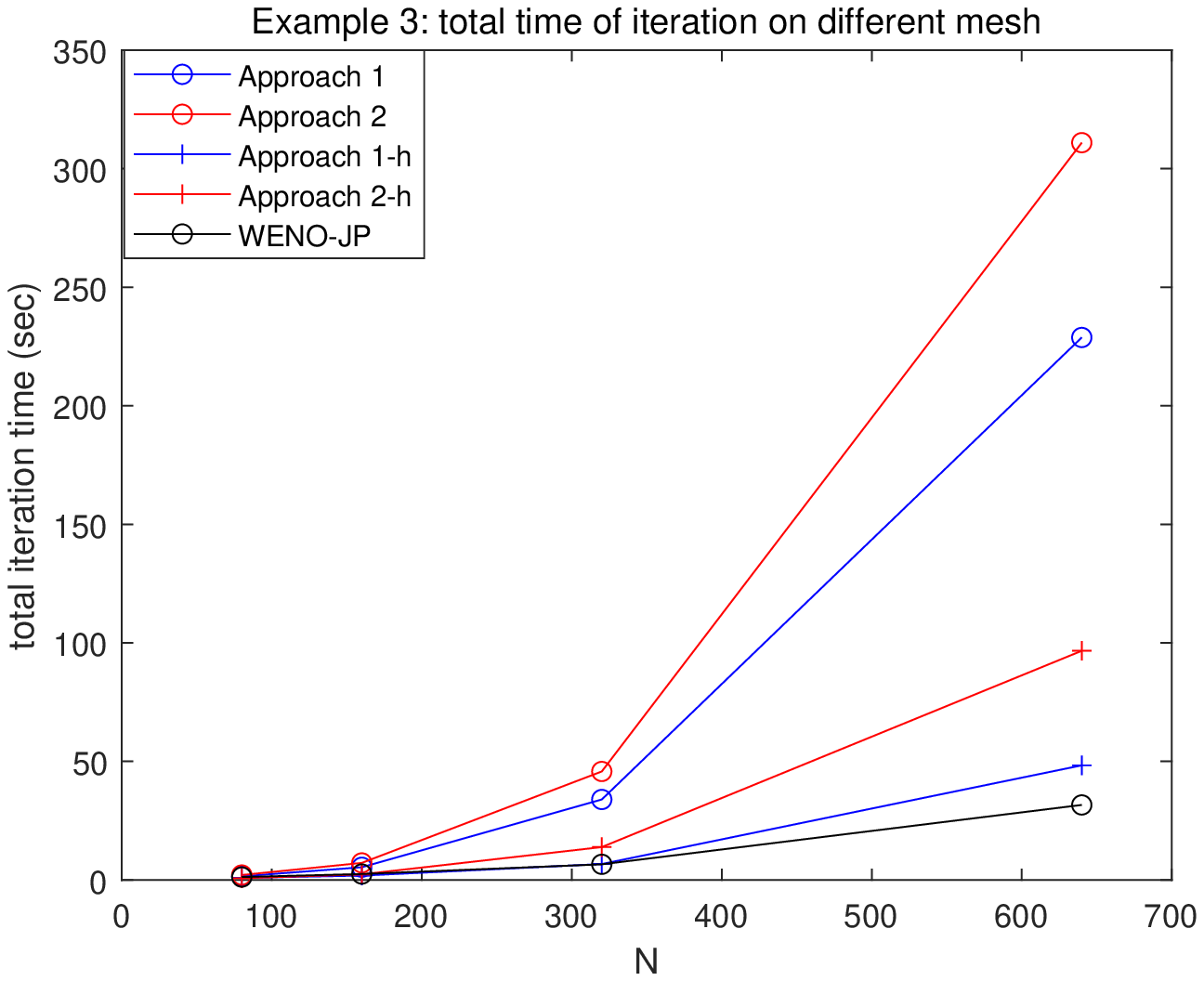}
	\includegraphics[width=5.4cm]{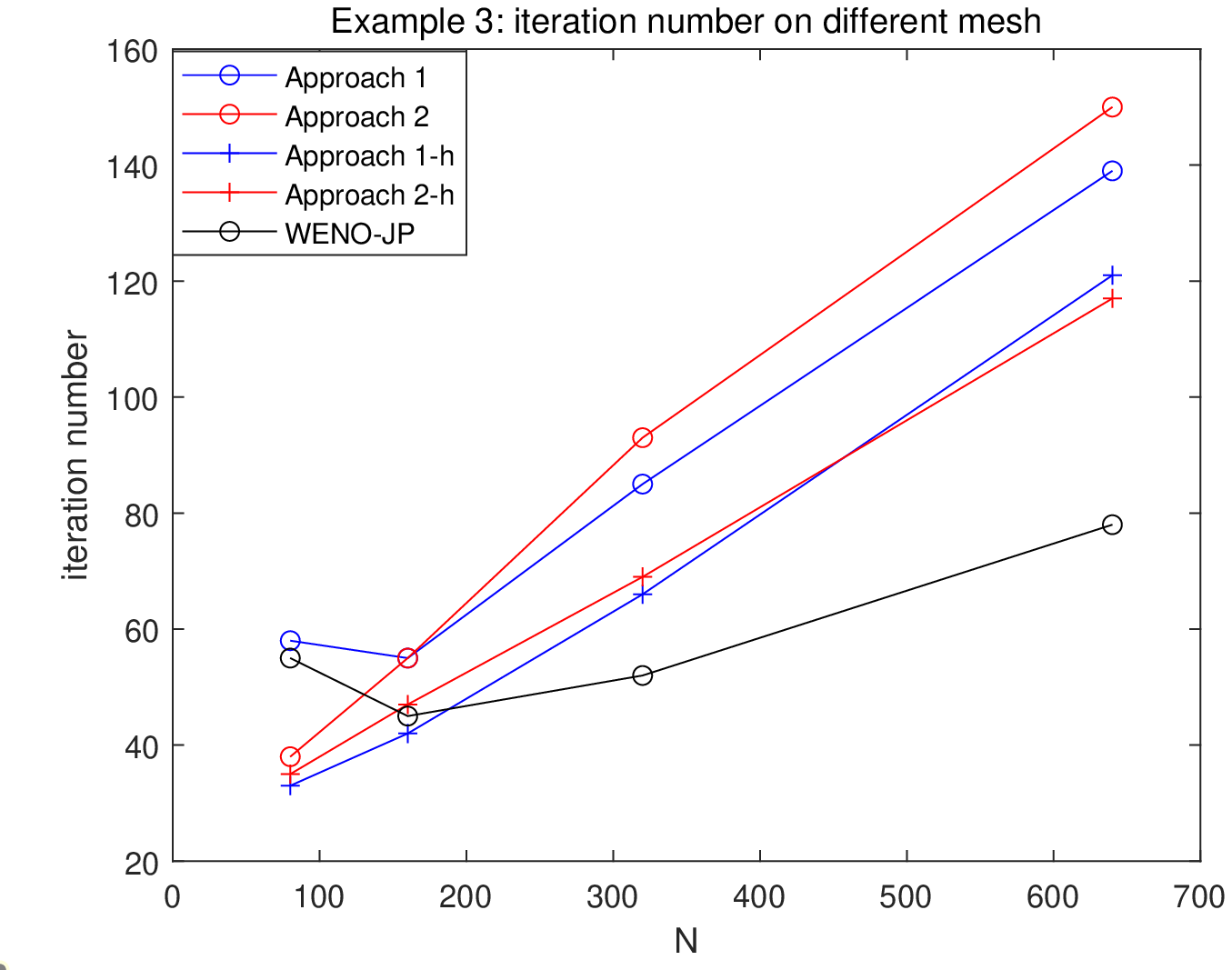}
	\includegraphics[width=5.4cm]{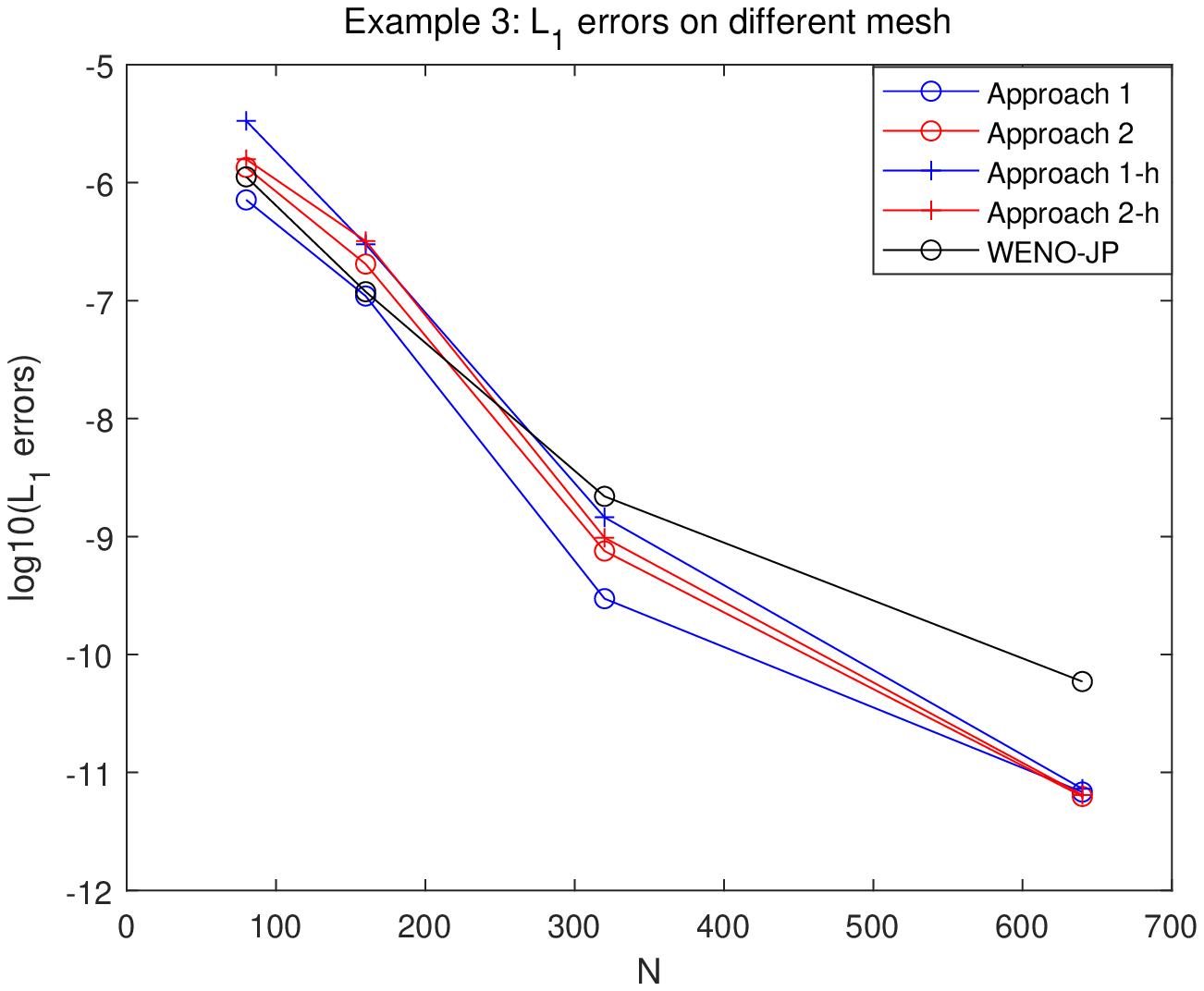}
\end{center}
\caption{Example 3 with Approach 1, 2, the hybrid Approach 1, 2 and WENO-JP. Left: mesh number $N$ vs CPU time; Middle: mesh number $N$ vs number of iterations; Right: mesh number N vs $L_{1}$ error.}\label{fig1h3}
\end{figure}
\begin{figure}[htbp!]
\begin{center}
	\includegraphics[width=5.4cm]{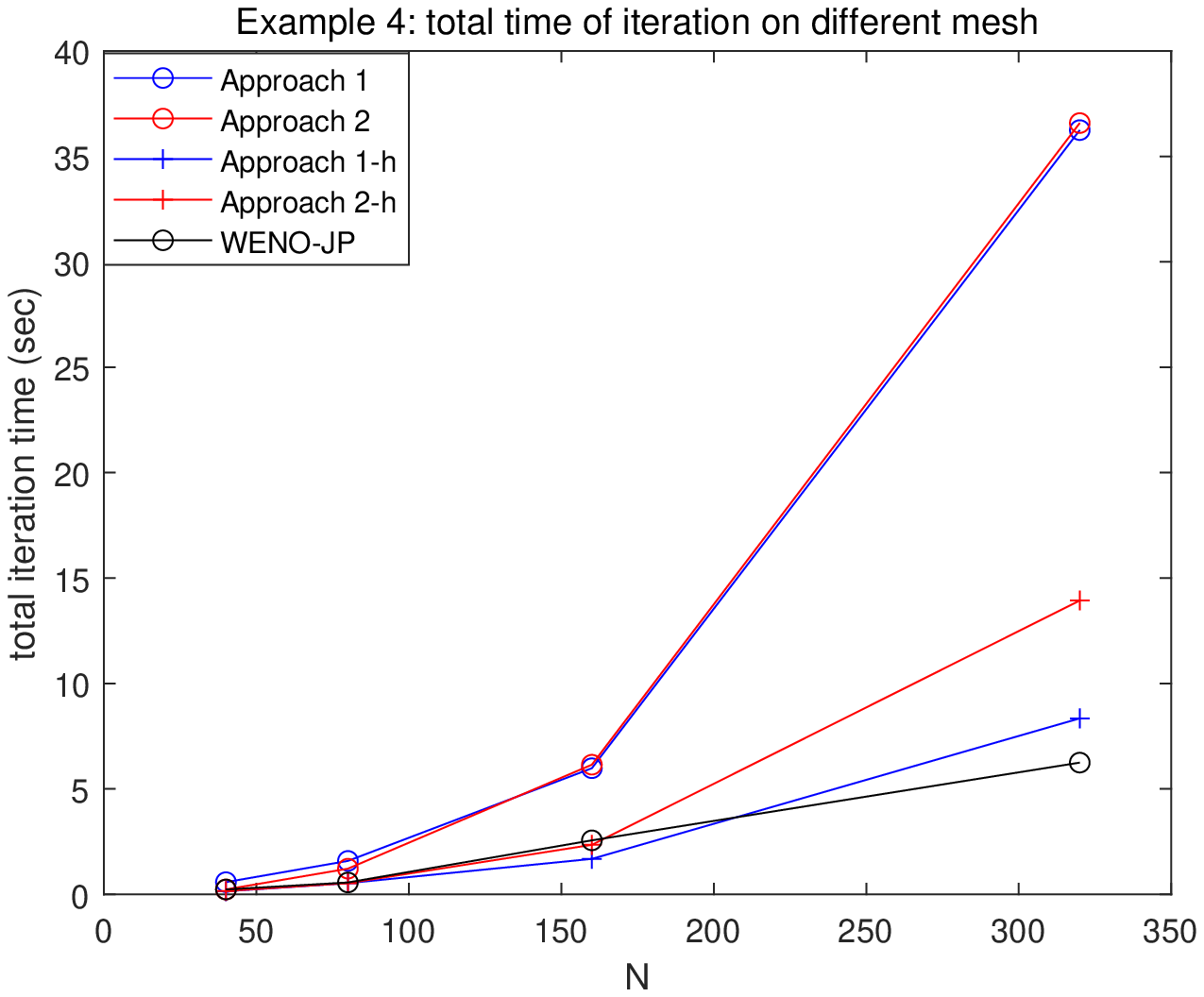}
	\includegraphics[width=5.4cm]{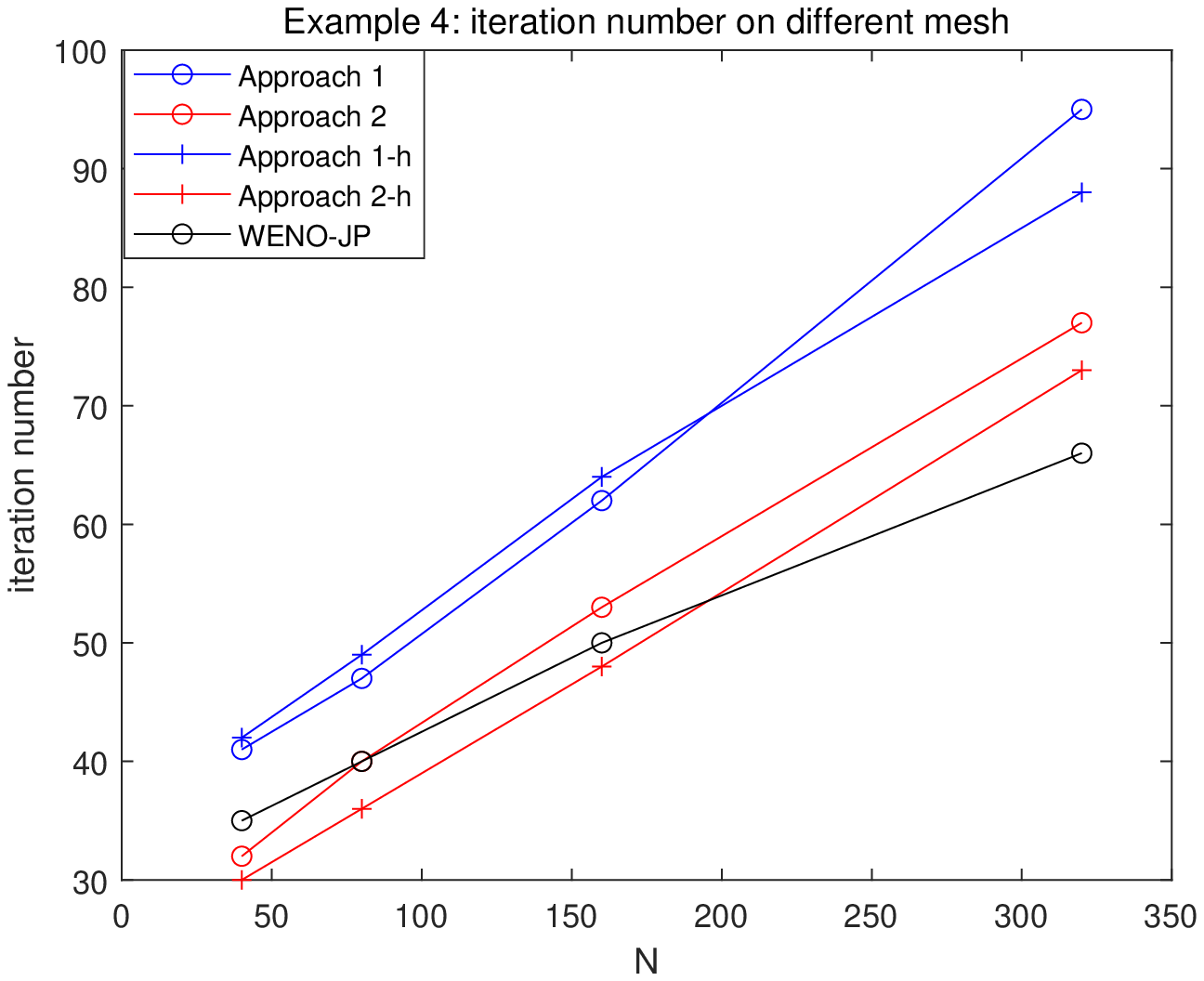}
	\includegraphics[width=5.4cm]{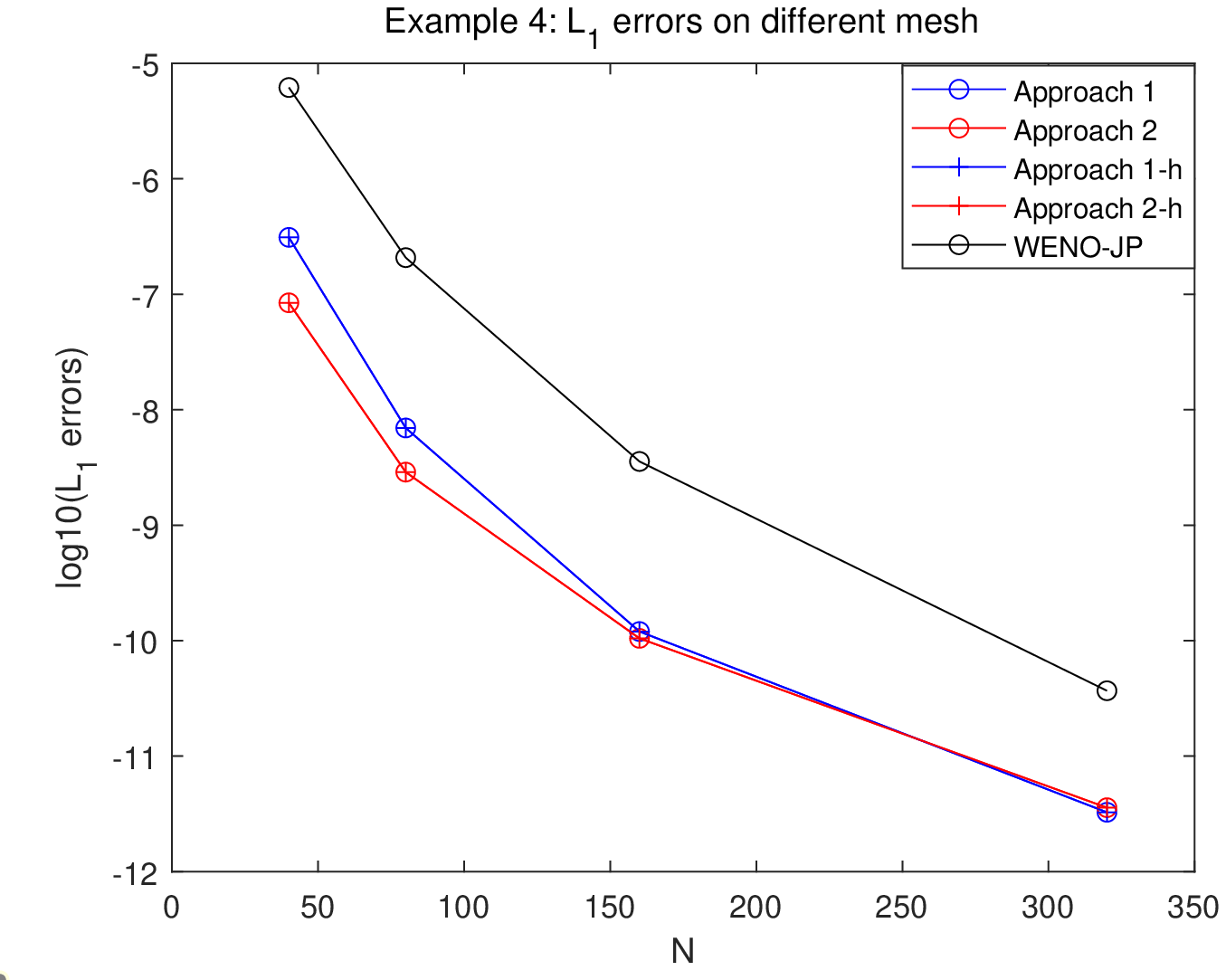}
\end{center}
\caption{Example 4 with Approach 1, 2, the hybrid Approach 1, 2 and WENO-JP. Left: mesh number $N$ vs CPU time; Middle: mesh number $N$ vs number of iterations; Right: mesh number N vs $L_{1}$ error.}\label{fig1h4}
\end{figure}
\begin{figure}[htbp!]
\begin{center}
	\includegraphics[width=5.4cm]{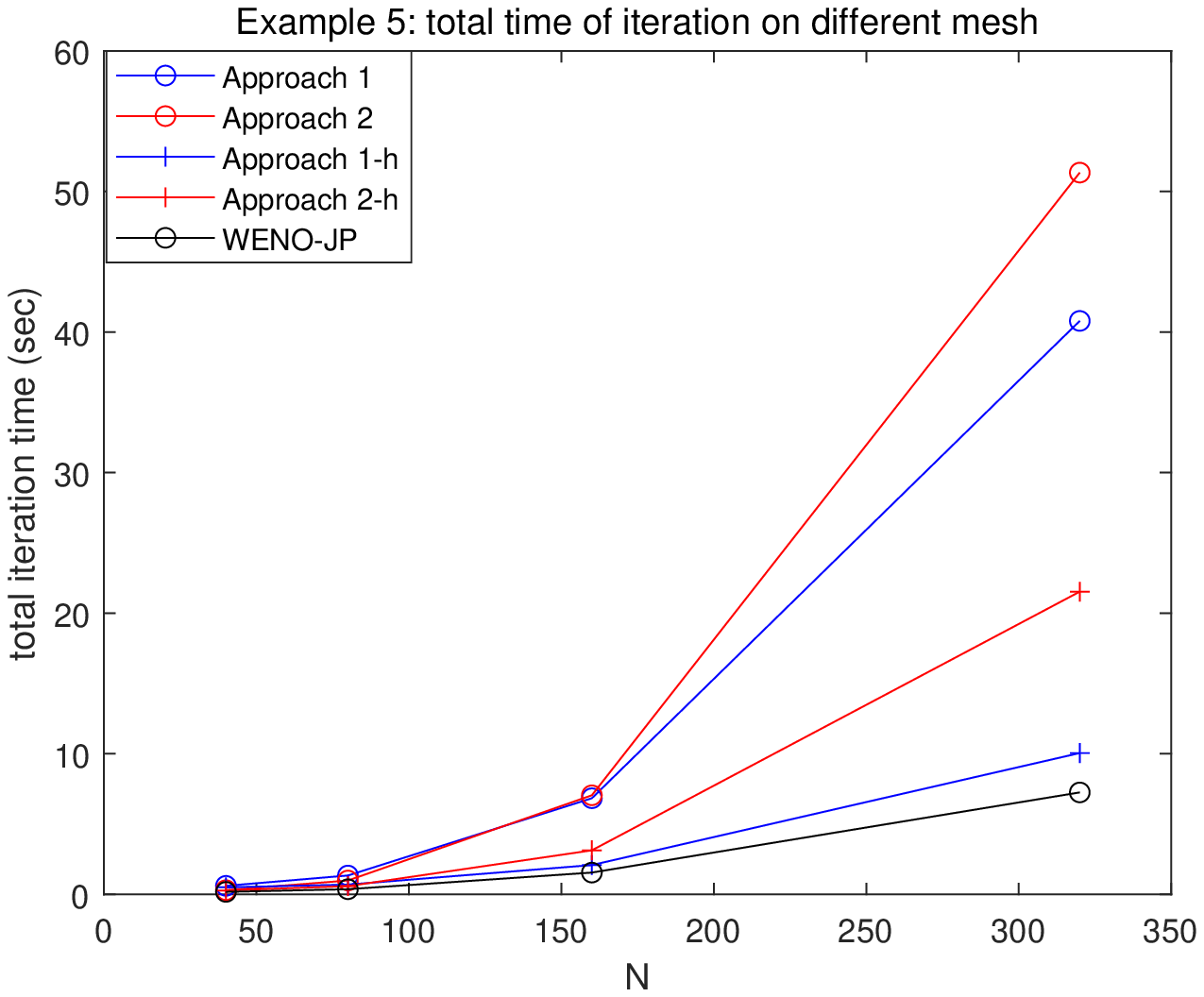}
	\includegraphics[width=5.4cm]{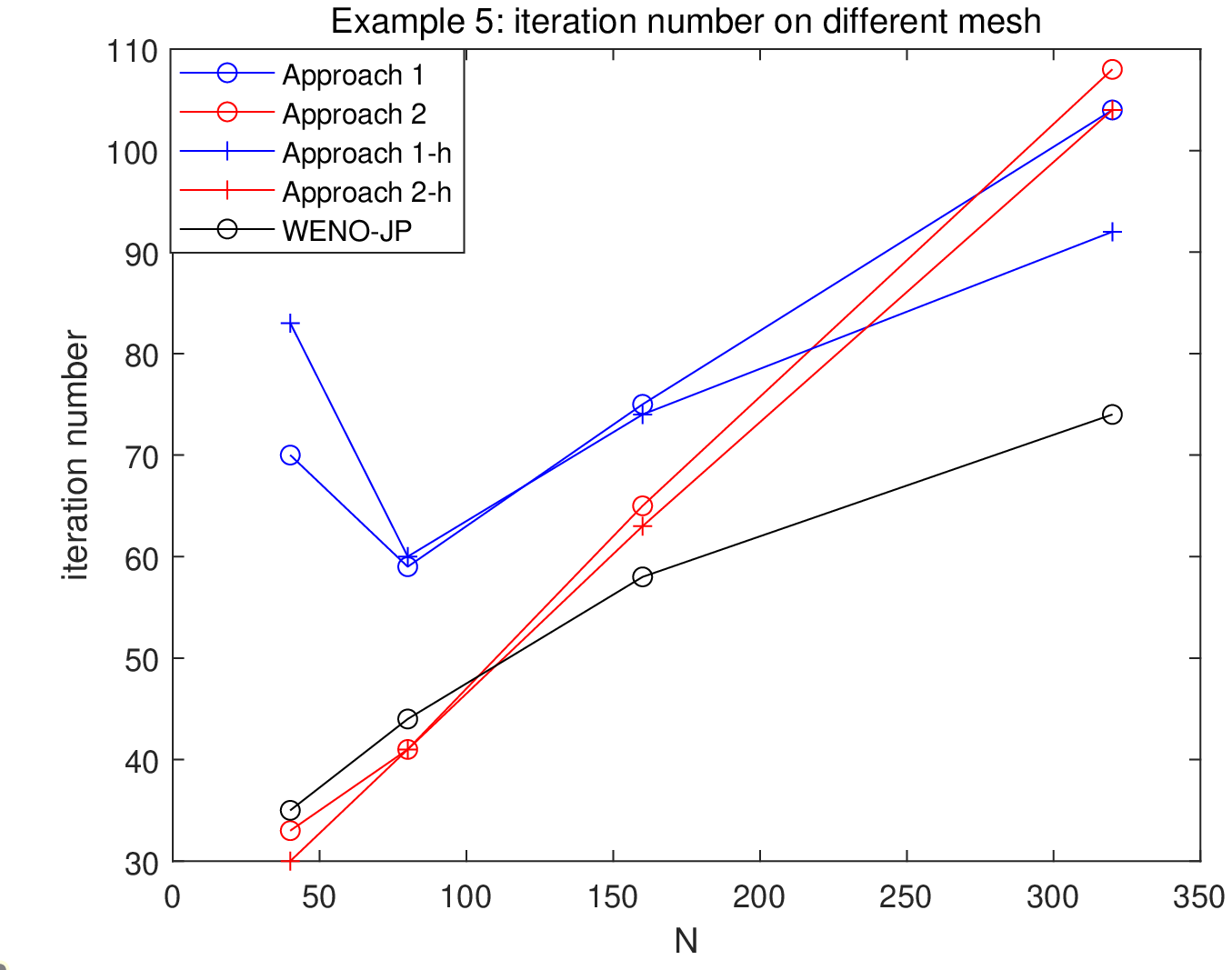}
	\includegraphics[width=5.4cm]{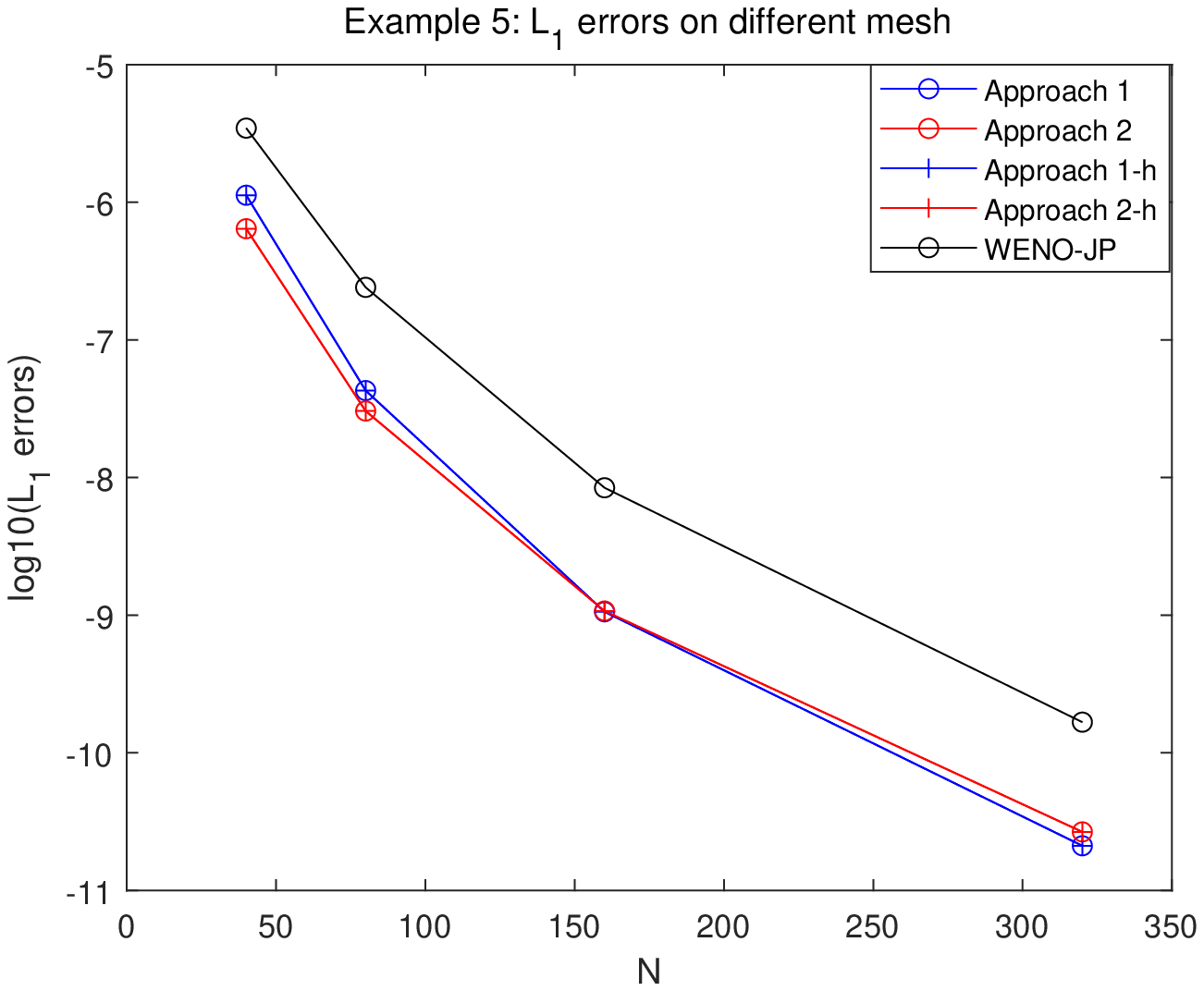}
\end{center}
\caption{Example 5 with Approach 1, 2, the hybrid Approach 1, 2 and WENO-JP. Left: mesh number $N$ vs CPU time; Middle: mesh number $N$ vs number of iterations; Right: mesh number N vs $L_{1}$ error.}\label{fig1h5}
\end{figure}
\begin{figure}[htbp!]
\begin{center}
	\includegraphics[width=5.4cm]{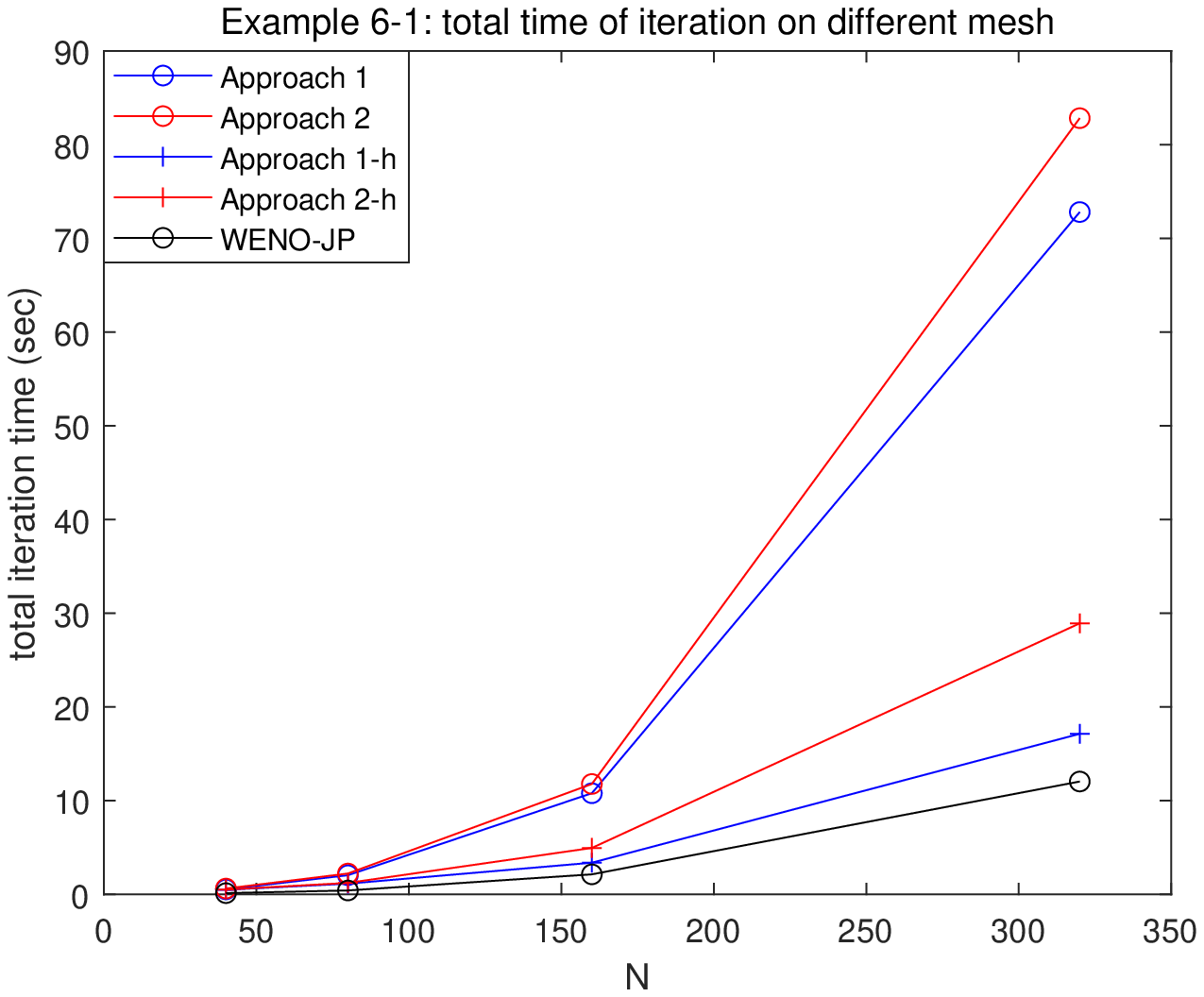}
	\includegraphics[width=5.4cm]{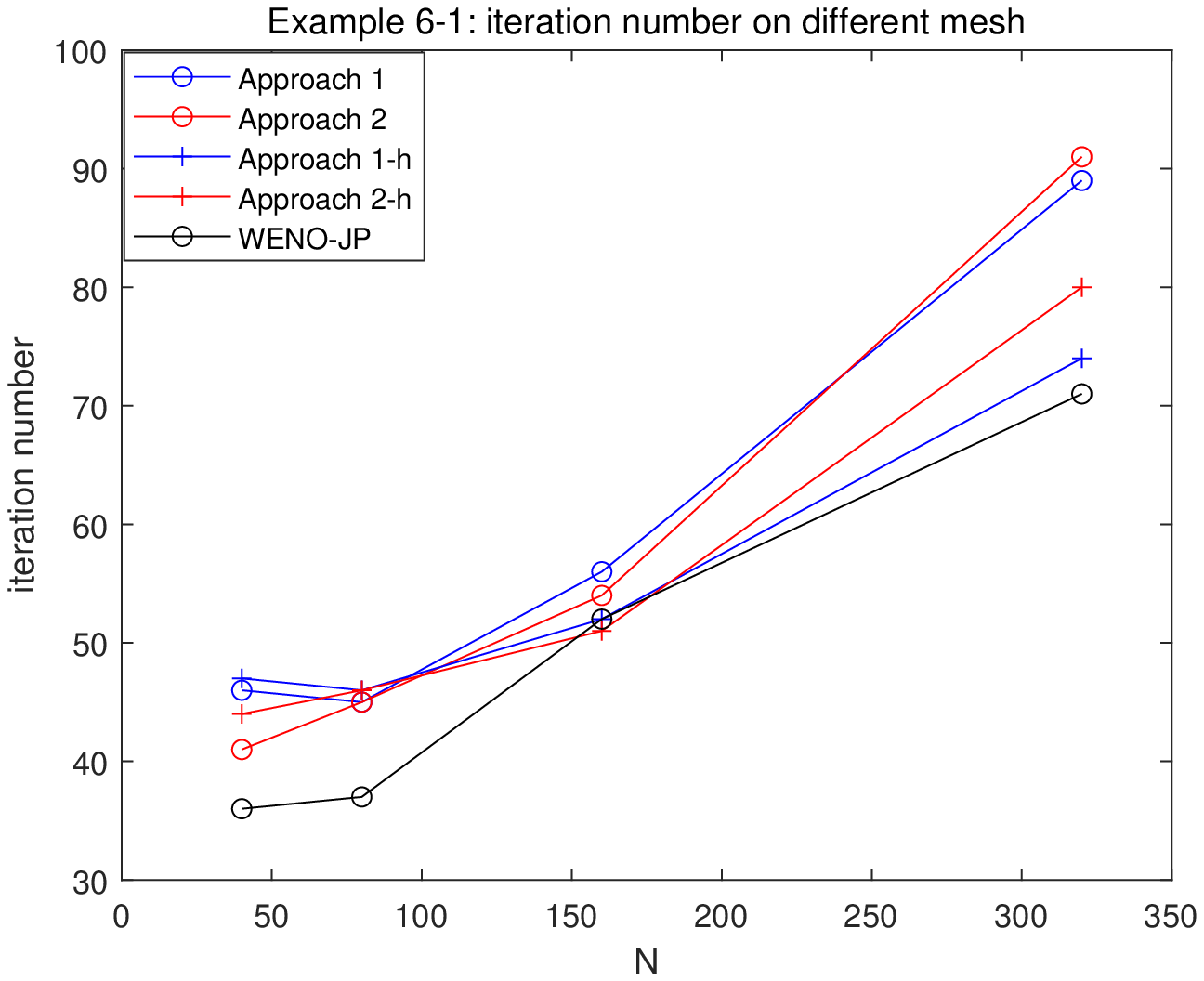}
	\includegraphics[width=5.4cm]{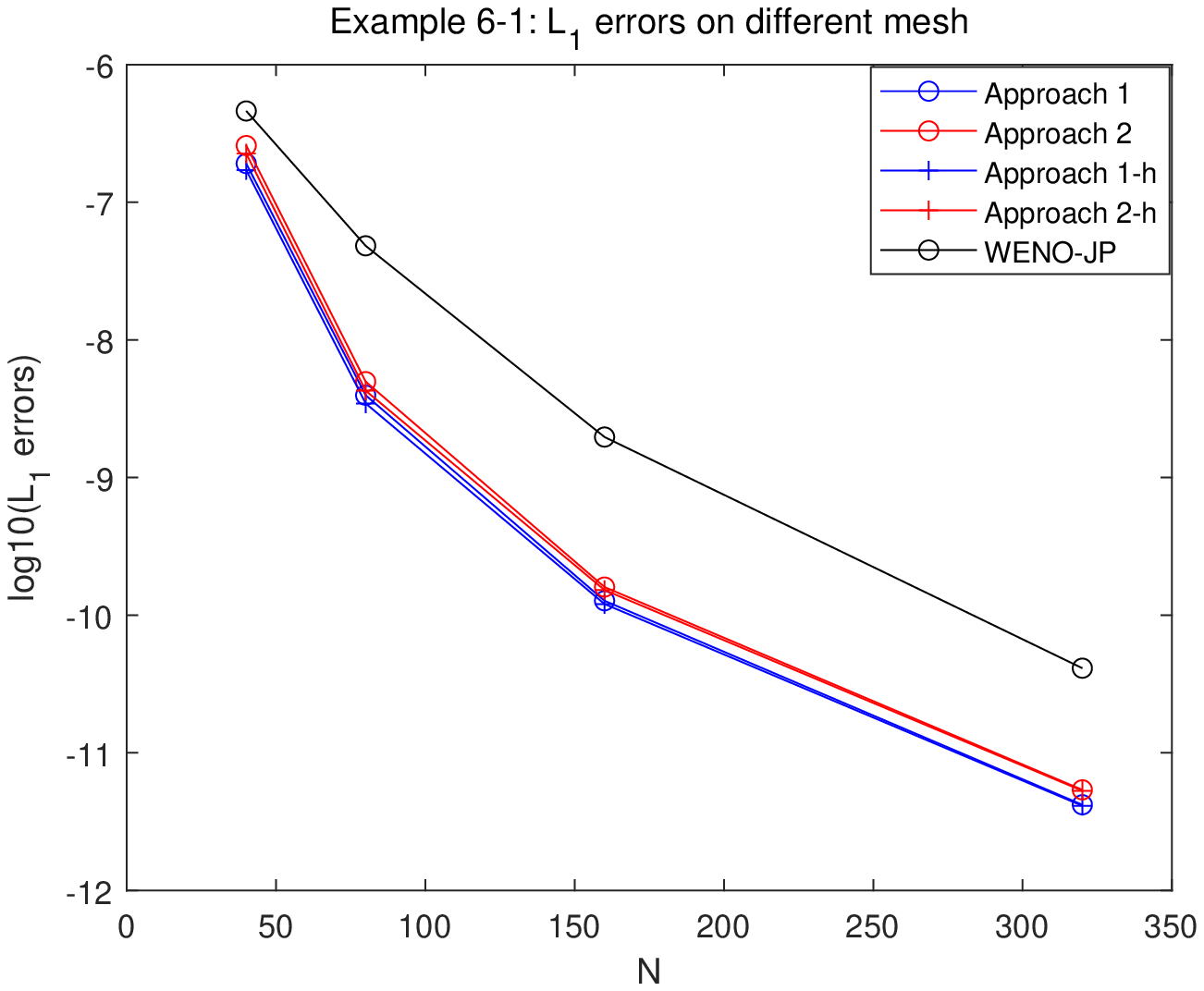}
\end{center}
\caption{Example 6-a with Approach 1, 2, the hybrid Approach 1, 2 and WENO-JP. Left: mesh number $N$ vs CPU time; Middle: mesh number $N$ vs number of iterations; Right: mesh number N vs $L_{1}$ error.}\label{fig1h6}
\end{figure}
\begin{figure}[htbp!]
\begin{center}
	\includegraphics[width=5.4cm]{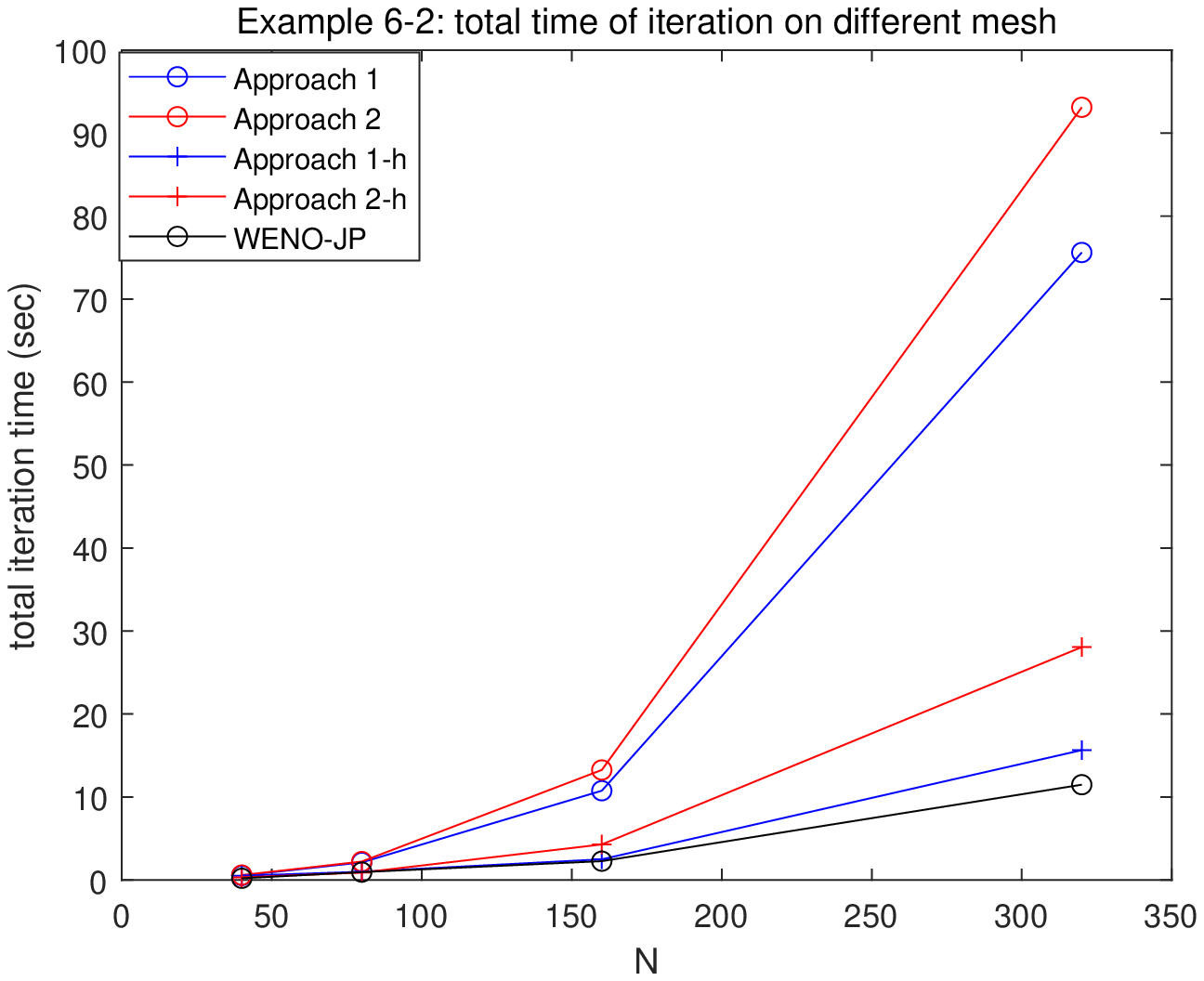}
	\includegraphics[width=5.4cm]{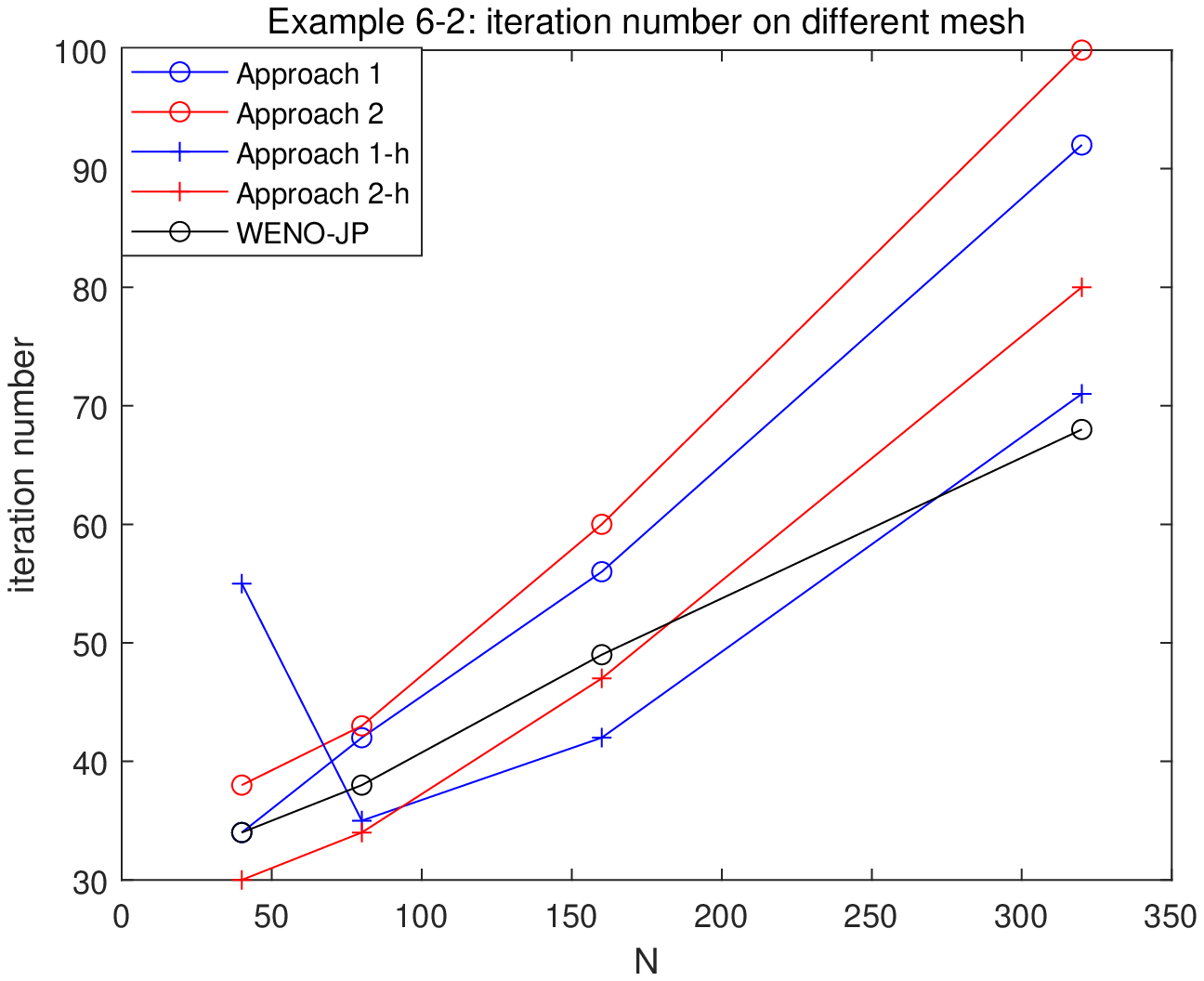}
	\includegraphics[width=5.4cm]{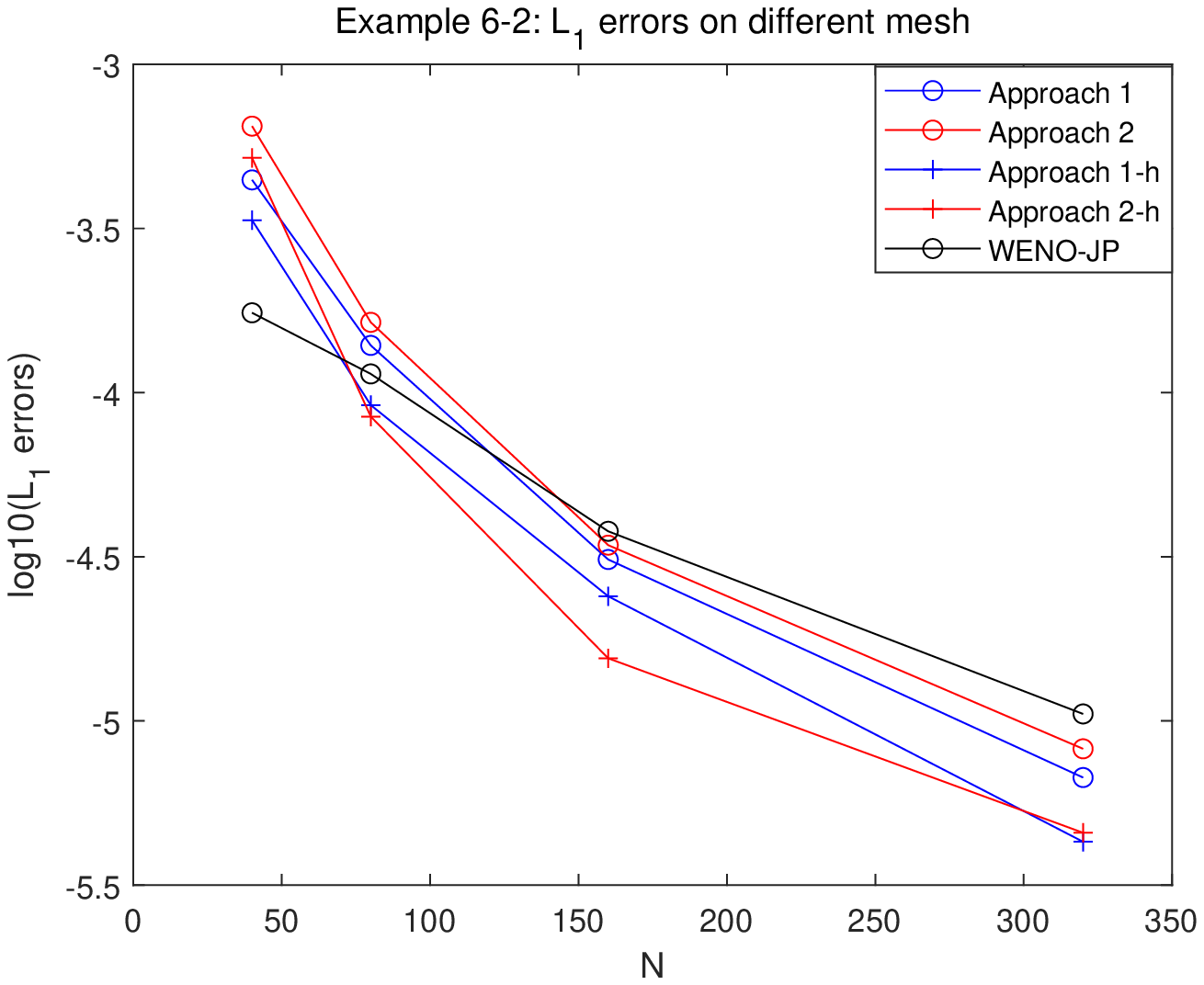}
\end{center}
\caption{Example 6-b with Approach 1, 2, the hybrid Approach 1, 2 and WENO-JP. Left: mesh number $N$ vs CPU time; Middle: mesh number $N$ vs number of iterations; Right: mesh number N vs $L_{1}$ error.}\label{fig1h7}
\end{figure}

\begin{figure}[htbp!]
\begin{center}
	\includegraphics[width=5.4cm]{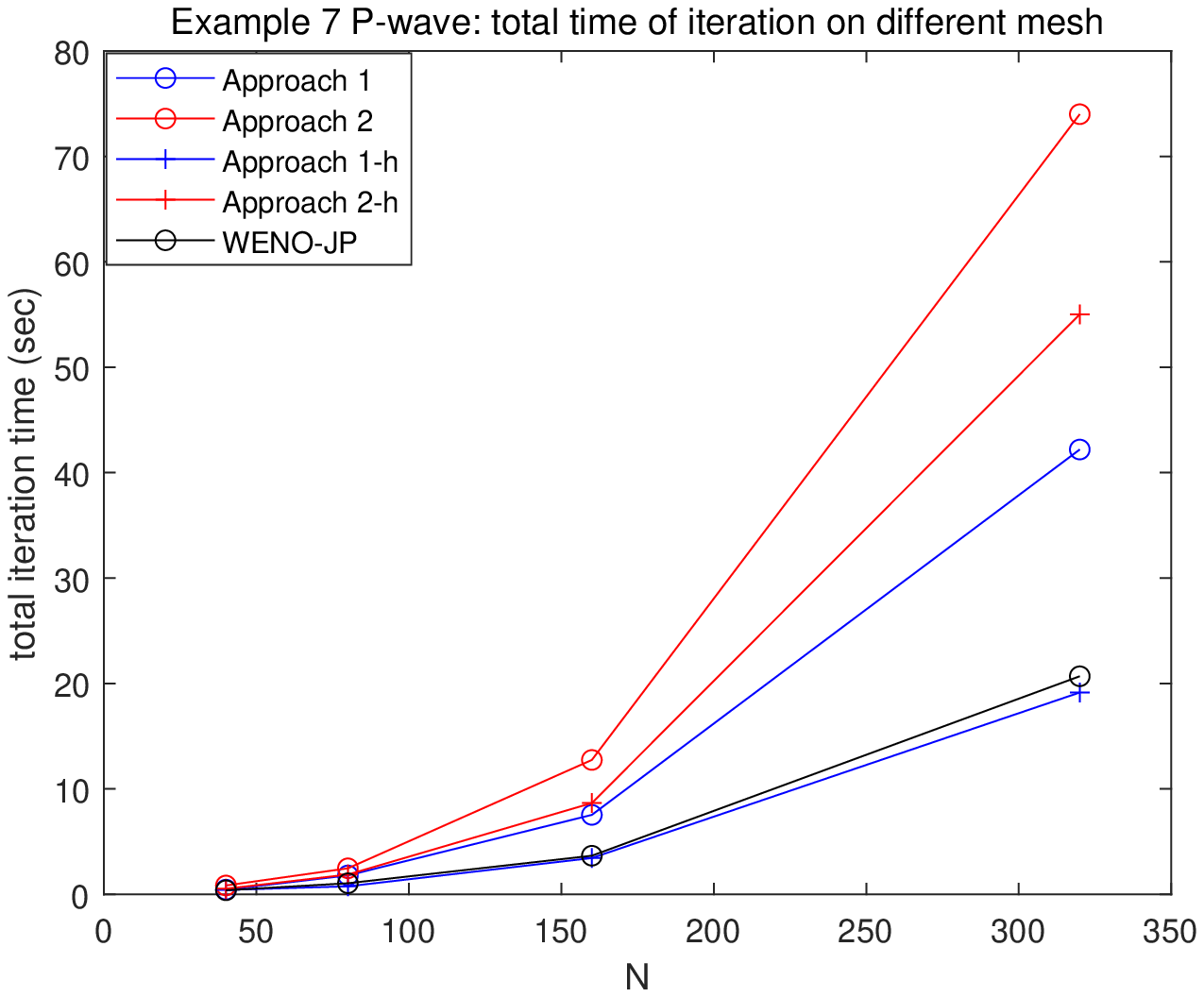}
	\includegraphics[width=5.4cm]{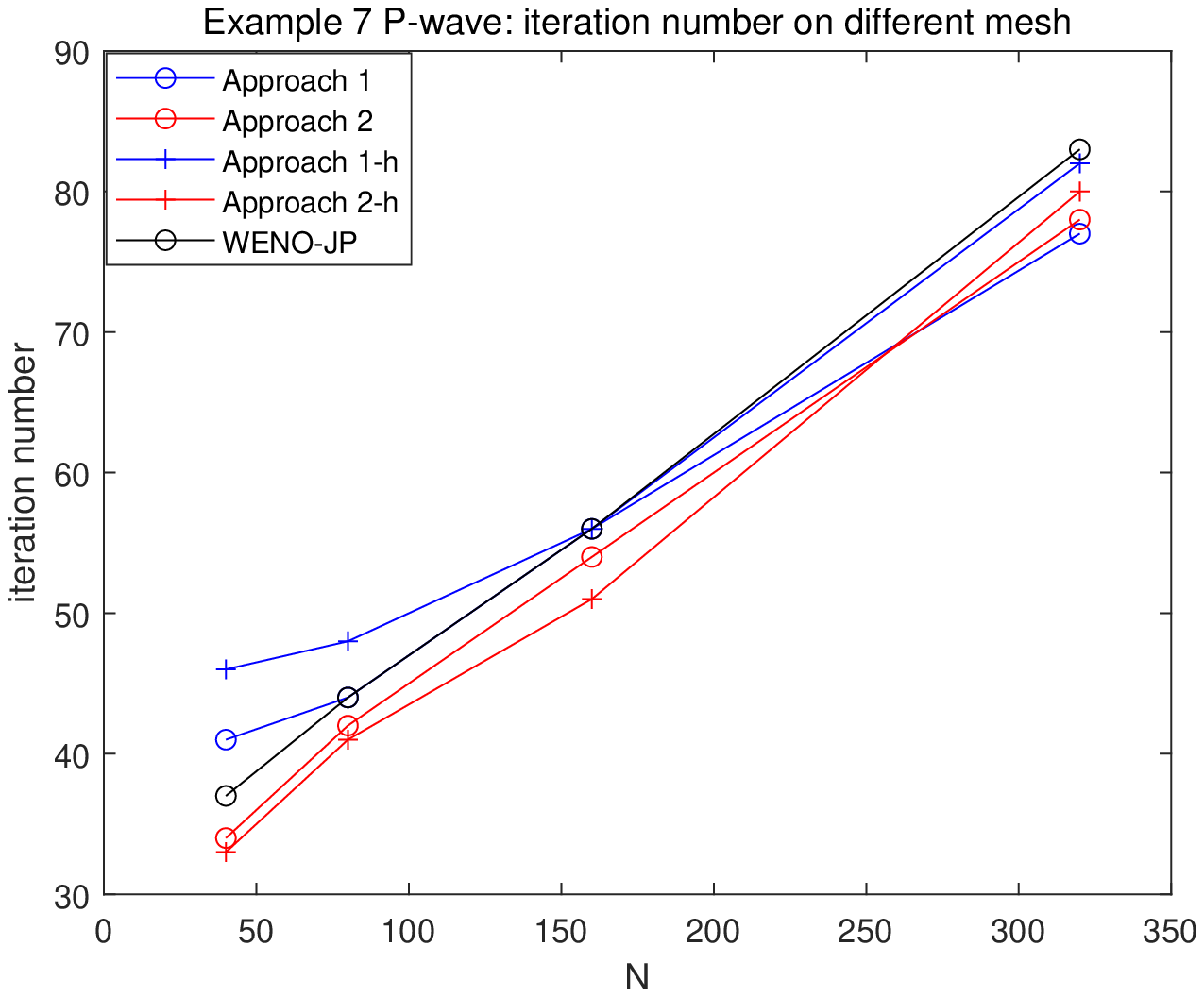}
	\includegraphics[width=5.4cm]{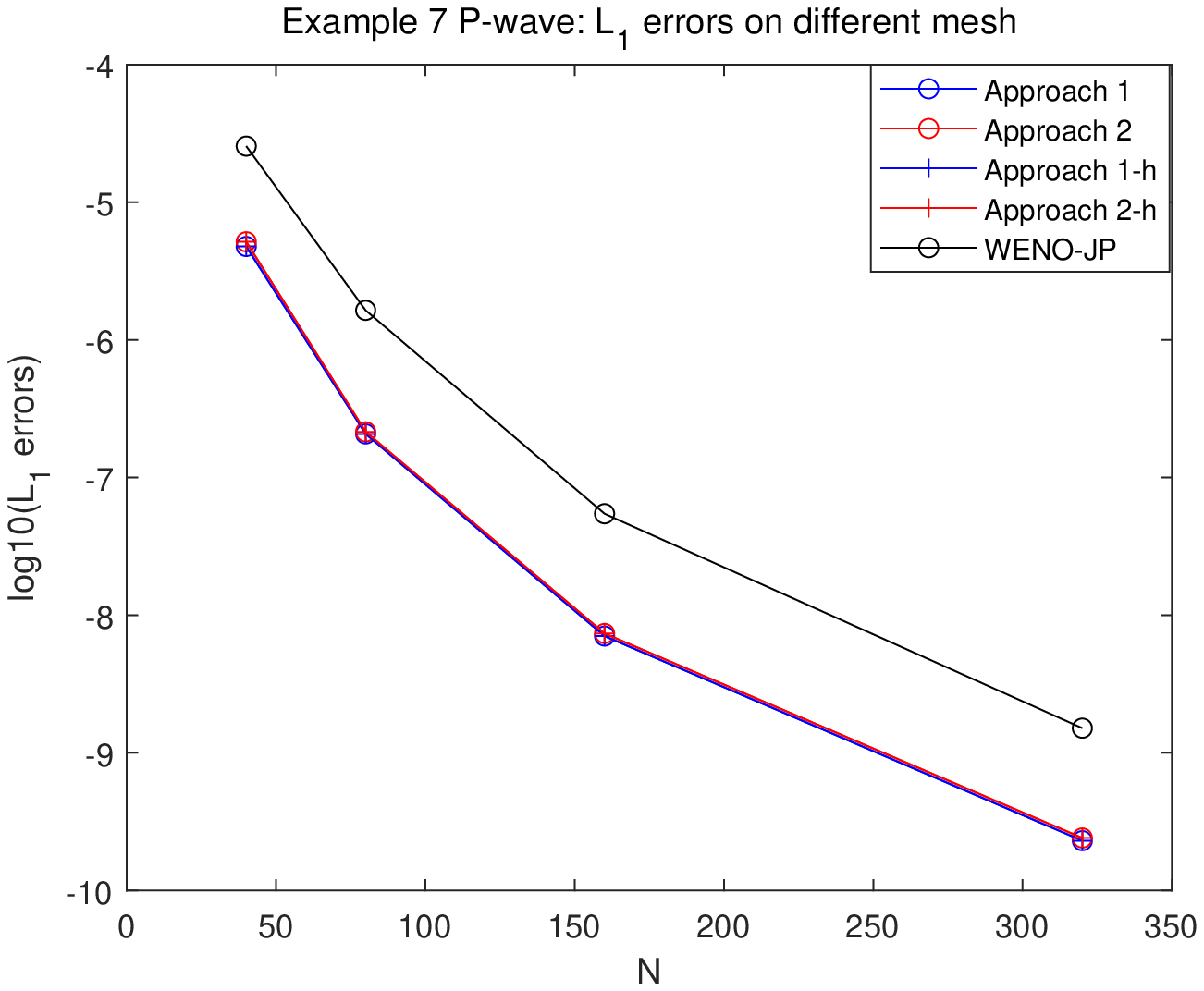}
\end{center}
\caption{Example 7 P-wave with Approach 1, 2, the hybrid Approach 1, 2 and WENO-JP. Left: mesh number $N$ vs CPU time; Middle: mesh number $N$ vs number of iterations; Right: mesh number N vs $L_{1}$ error.}\label{fig1h8}
\end{figure}
\begin{figure}[htbp!]
\begin{center}
	\includegraphics[width=5.4cm]{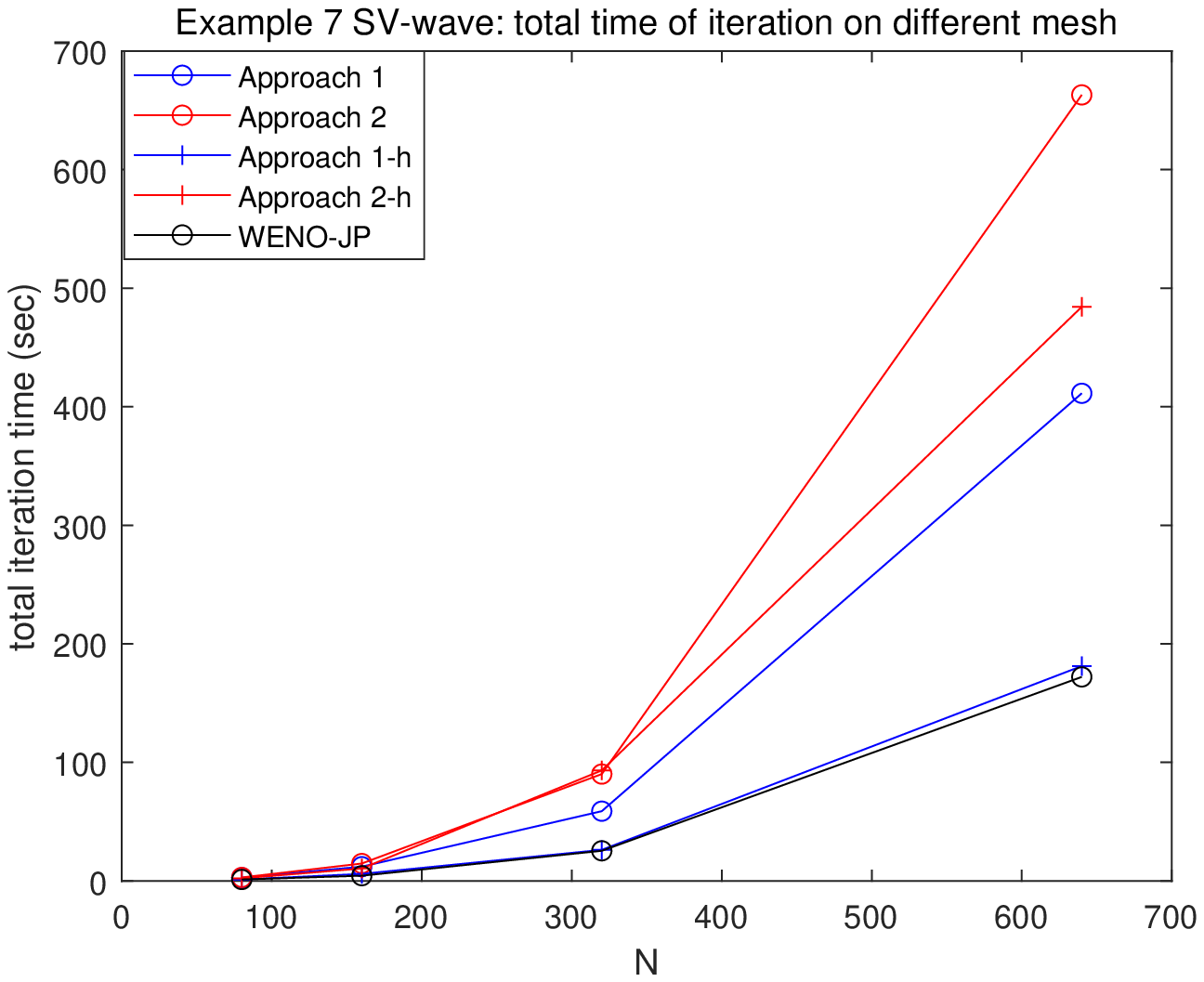}
	\includegraphics[width=5.4cm]{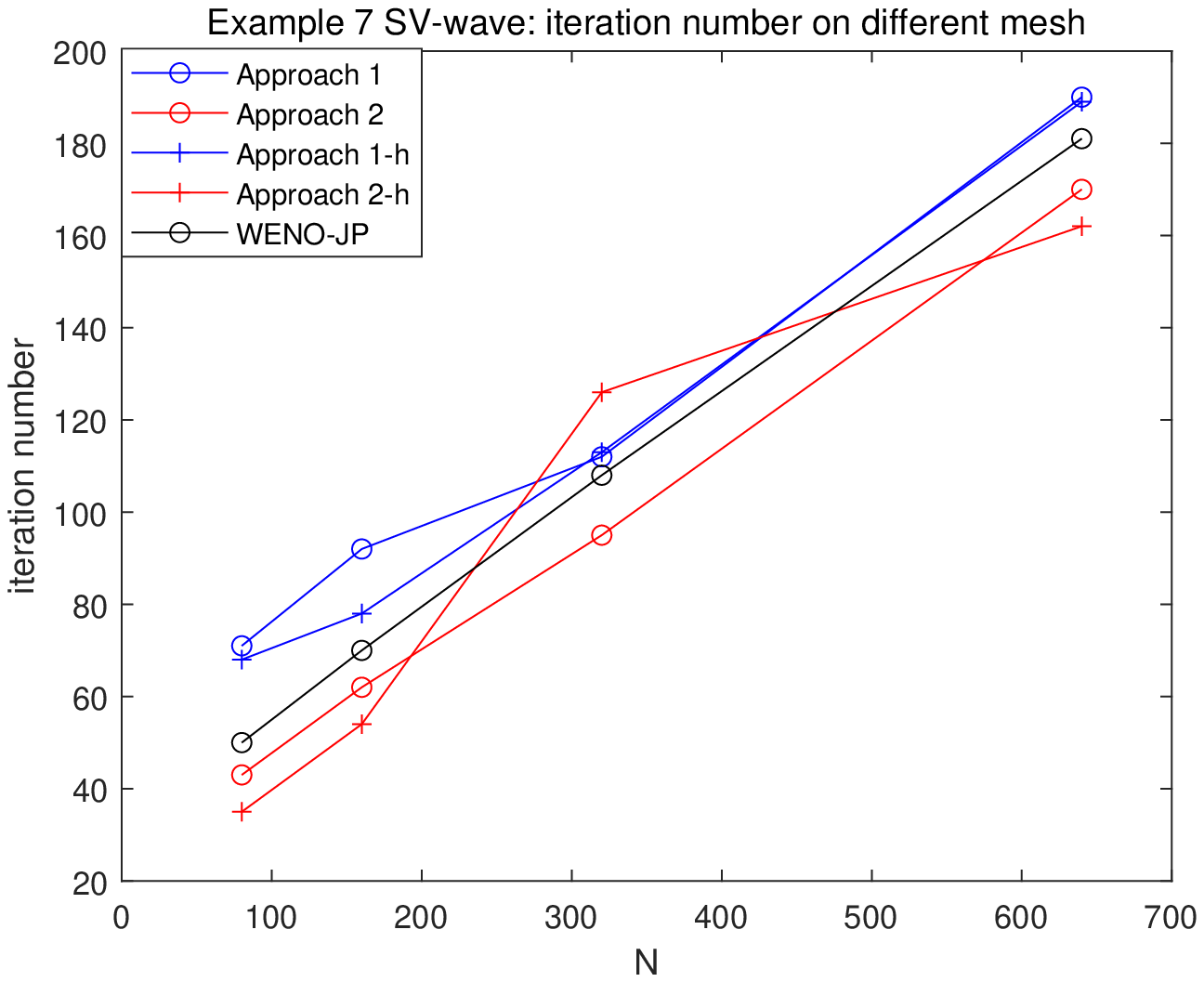}
	\includegraphics[width=5.4cm]{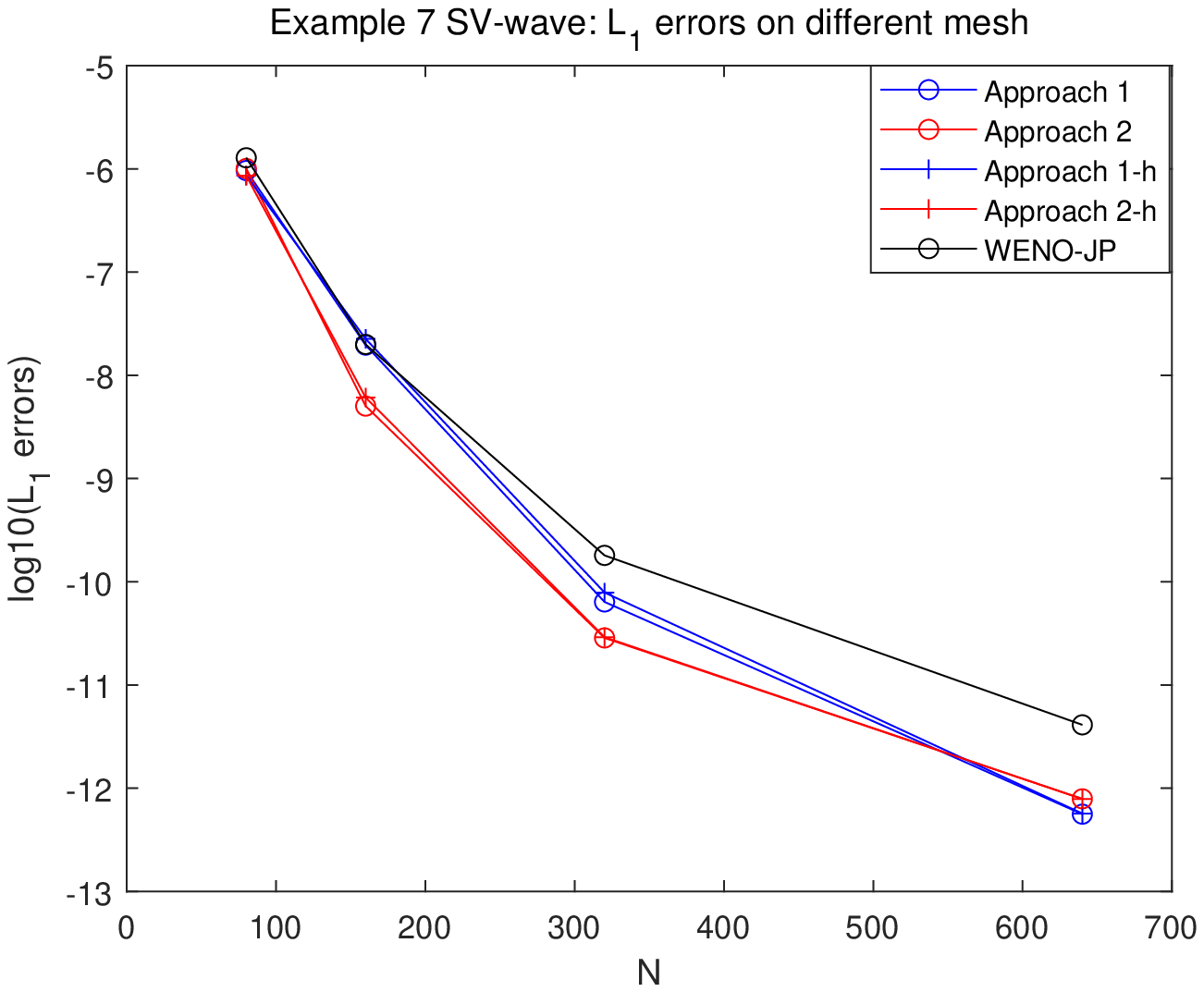}
\end{center}
\caption{Example 7 SV-wave with Approach 1, 2, the hybrid Approach 1, 2 and WENO-JP. Left: mesh number $N$ vs CPU time; Middle: mesh number $N$ vs number of iterations; Right: mesh number N vs $L_{1}$ error.}\label{fig1h9}
\end{figure}
\newpage
\begin{remark}
At the end of this section, we propose a small trick to further reduce the computational cost, namely to freeze the nonlinear weights. Since the targeting problem is steady state and iterative method is used to solve it, the nonlinear weights changes only slightly when the iteration is close to converge.
Therefore, when $L_{1}$ error of the nonlinear weights between two iteration steps is less than the given threshold, we can fix these nonlinear weights in following iterations until convergence.
Numerical tests demonstrate that the value of $\epsilon_{1}$ can be taken as $10^{-4}$ or $10^{-5}$, and this could effectively reduce the computational costs.
\end{remark}

\section{Conclusion Remark}\label{section7}
In this work, we have combined the fifth order HWENO-ZZQ scheme with the fast sweeping idea, to design efficient algorithms for directly solving static Hamilton-Jacobi equations. The novel approach seeks to use the updated $\phi$ to directly approximate the spatial derivatives of $\phi$, and there is no need to introduce and solve additional equations.  As a comparison, the second approach is based on the traditional HWENO idea, with additional equations governing the spatial derivatives of $\phi$. The first approach has great savings in both computational time and storage, which improves the computational efficiency of the traditional HWENO scheme. Extensive numerical experiments demonstrate that these two methods perform well numerically and lead to smaller numerical errors when compared with WENO methods. A hybrid strategy which combines both linear and HWENO-ZZQ reconstruction is also proposed and tested, which yields additional savings in computational time. Especially, the hybrid version of the proposed novel HWENO method enjoys more savings in computational time than the traditional HWENO approach.

\end{document}